\pdfoutput=1
\documentclass[11pt]{article}

\usepackage[final]{acl}

\usepackage{times}
\usepackage{latexsym}
\usepackage{hyperref}
\usepackage[T1]{fontenc}
\usepackage[utf8]{inputenc}

\usepackage{microtype}

\usepackage{inconsolata}
\usepackage{enumitem}

\usepackage{booktabs}
\definecolor{azure(colorwheel)}{rgb}{0.0, 0.5, 1.0}

\usepackage{graphicx}
\usepackage{amsmath}
\usepackage{multirow}
\usepackage{pifont}
\usepackage{amssymb} % 用于数学符号
\DeclareMathOperator*{\argmin}{arg\,min}
\usepackage{tabularx}
\usepackage{listings}
\usepackage{xcolor}
\usepackage{graphicx}
\usepackage{subcaption}

\lstdefinestyle{markdownstyle}{
    backgroundcolor=\color{lightgray},
    basicstyle=\ttfamily,
    breaklines=true,
    frame=single,
    keywordstyle=\color{blue},
    commentstyle=\color{green!60!black},
    morekeywords={Problem, Setup, Examples, Task}, % 自定义高亮关键词
    columns=flexible
}

\title{Exploiting Edited Large Language Models as General Scientific Optimizers}

\author{
    Qitan Lv\textsuperscript{1},
    Tianyu Liu\textsuperscript{1}\thanks{Equal Contributions.},
    Hong Wang\textsuperscript{1}\thanks{The Corresponding author.},
    \\
    \textsuperscript{1} University of Science and Technology of China
    \\
    \texttt{\{qitanlv, tianyu\_liu, wanghong1700\}@mail.ustc.edu.cn}
}

\begin{document}
\maketitle

\begin{abstract}

Large language models (LLMs) have been widely adopted in mathematical optimization in scientific scenarios for their extensive knowledge and advanced reasoning capabilities. 
Existing methods mainly focus on utilizing LLMs to solve optimization problems in a prompt-based manner, which takes observational feedback as additional textual descriptions. 
However, due to LLM's \textbf{high sensitivity to the prompts} and \textbf{tendency to get lost in lengthy prompts}, these methods struggle to effectively utilize the {observational} feedback from each optimization step, which severely hinders the applications for real-world scenarios. 
To address these challenges, we propose a conceptually simple and general {bi-level} optimization method, namely \textbf{G}eneral \textbf{S}cientific \textbf{O}ptimizers (GSO).
Specifically, GSO first utilizes inner-level simulators as experimental platforms to evaluate the current solution and provide observational feedback. Then, LLMs serve as knowledgeable and versatile scientists, generating new solutions by refining potential errors from the feedback as the outer-level optimization.
Finally, simulations together with the expert knowledge in LLMs are jointly updated with bi-level interactions via model editing.
Extensive experiments show that GSO consistently outperforms existing state-of-the-art methods using \textit{six} different LLM backbones on \textit{seven} different tasks, demonstrating the effectiveness and a wide range of applications. 

\end{abstract}
\section{Introduction} \label{sec:intro}
Optimization is ubiquitous across a wide range of scientific domains, spanning mathematics, physics, chemistry, pharmacology, etc \cite{back_sgd, math_opt, op_chem}. Various research endeavors innovate within its field, creating methods tailored to its specific challenges and nuances to automate and accelerate the process of mathematical optimization in scientific scenarios \cite{ai_age}. However, the need for customization of optimization algorithms to address specific challenges highlights the absence of a universally applicable philosophy \cite{logic_sci, sos}, which is crucial for establishing a standardized optimization framework and enhancing the efficiency of scientific research. We aim to transcend specific domains and provide a generalized approach to boost mathematical optimization in scientific scenarios.

\begin{figure}[t]
    \centering 
    \includegraphics[width=1\columnwidth]{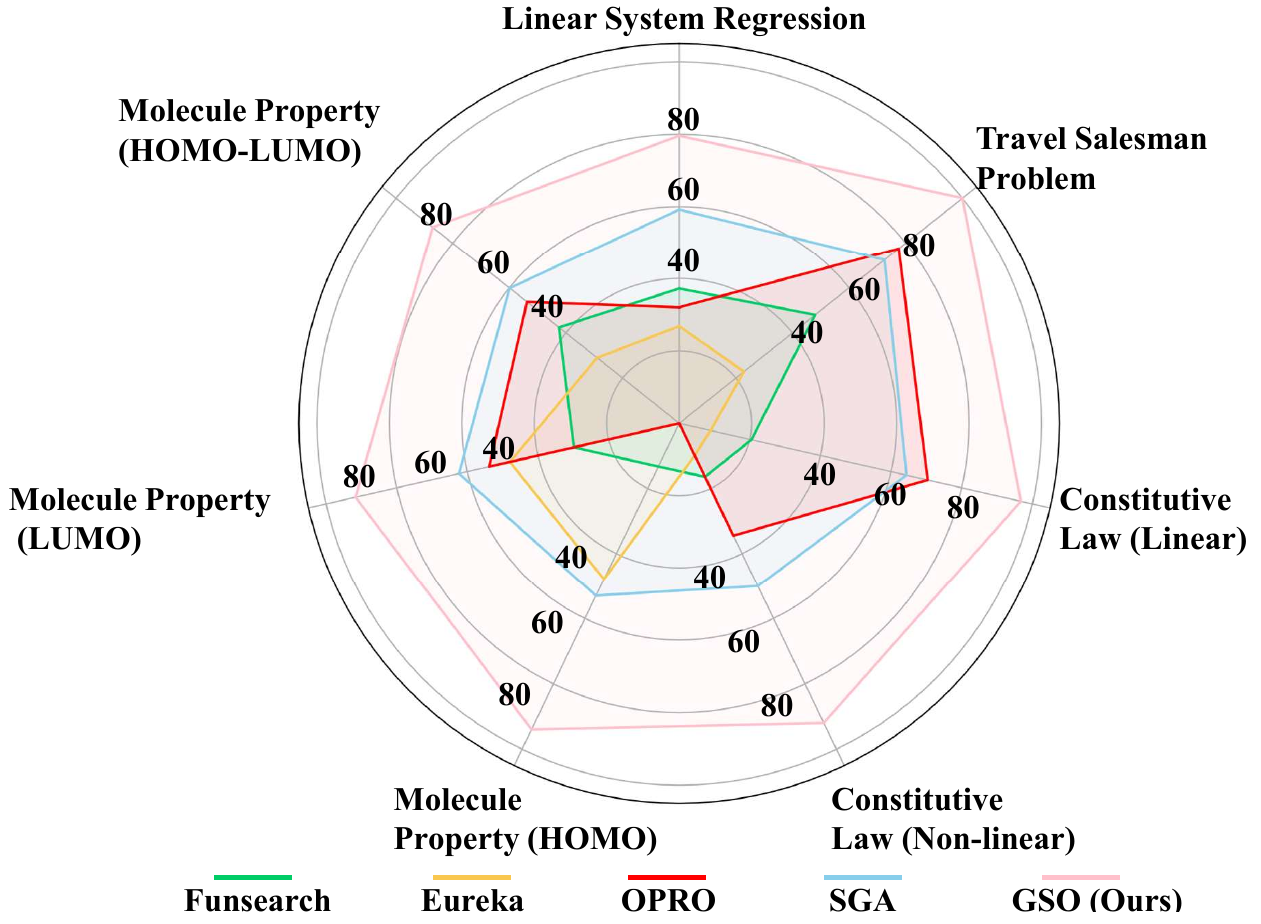}
    \caption{GSO achieves state-of-the-art performance on a broad
range of scientific optimization tasks compared with existing methods, using LLama 3 8B \cite{llama3} as the backbone. Results of other five LLMs are in Figures \ref{fig:leida_qiansan} and \ref{fig:leida_housan}. }
    \label{fig:score_overview}
    % \vspace{-5mm}
\end{figure}

Standing out as versatile tools with vast knowledge repositories, large language models (LLMs) have recently risen to prominence in optimization across scientific domains for their expansive knowledge bases, advanced reasoning capabilities, and human-friendly natural language interface \cite{ai4sci}. Extensive research efforts have been devoted to boosting general mathematical optimization in scientific scenarios. Canonical methods mainly focus on fine-tuning LLMs using domain-specific data to align natural language with scientific information, such as chemical structures \cite{chem_2024, chem_berta} or drug structures \cite{ai_drug_2021}. 

However, these approaches are constrained to specific domains and require substantial amounts of data and extensive computation resources for broader applicability. Recently, prompt-based iterative optimization methods---which enhance the inherent capabilities of pre-trained LLMs by incorporating the optimization feedback to LLMs---have emerged as a promising approach for advancing scientific optimization \cite{opro}. 

Extensive researches have explored leveraging LLMs as optimizers or agents \cite{robots1} for tasks such as mathematical problem-solving \cite{nature_math, opro}, conducting chemical experiments \cite{nature_chem}, advancing physical scientific discovery \cite{bi_level}, molecular discovery \cite{module_chem}, and drug discovery \cite{gpt_drug}.

Albeit with multiple benefits of the prompt-based methods, they confront one significant challenge that severely hinders their general applications---struggling to effectively utilize observational feedback. This challenge primarily stems from two limitations inherent in existing prompt-based methods: \textbf{(i) }\textbf{LLMs are shown to be sensitive to the prompt format} \cite{sen_1, sen_2, sen_3}. In particular, semantically similar prompts can yield drastically different performance \cite{dra2,dra3,opro}, and the optimal prompt formats may be model-specific and task-specific, which severely limits the generalizability across different scientific tasks.  

{\textbf{(ii) }\textbf{LLMs may get lost in lengthy prompt} \cite{coft_2024}. In multi-round iterative optimization, the input prompt can become increasingly lengthy to trace the optimization trajectory, which may distract the LLMs' reasoning and result in the \textit{lost in the middle issue} \cite{distract}. Despite LLMs' ability to process long contexts, performances significantly decrease as the input grows longer, even for models explicitly designed for long contexts \cite{coft_2024}.}

{
To address this challenge, we propose a conceptually simple, flexible,
and general method, namely \textbf{G}eneral \textbf{S}cientific \textbf{O}ptimizers (GSO). 
Specifically, GSO is a bi-level optimization method involving \textbf{inner-level} optimization and \textbf{outer-level} optimization, together with \textbf{bi-level} interaction between them.
For a given optimization task, GSO will iteratively conduct the following process: (i) the inner-level optimization {first} employs simulators serving as an experimental platform to provide observational feedback for reasoning and refining; (ii) the \textbf{outer-level} optimization {then} utilizes knowledgeable LLMs as reasoning agents to generate hypotheses, reason with observational feedback, and refine the previous hypothetical solutions. It also devises a dynamic exploit-and-explore strategy to adaptively adjust the LLM reasoning trajectory based on observational feedback; (iii) finally, the \textbf{bi-level} optimization jointly updates the simulations together with the knowledge embedded within LLMs via model editing.

}

GSO is a novel framework to boost general mathematical optimization in scientific scenarios. 

As shown in Figure \ref{fig:score_overview}, extensive experiments demonstrate the effectiveness and generalization of our GSO, leading significant and consistent superiority using \textit{six} different LLM backbone on \textit{seven} different tasks than existing state-of-the-art methods.

\begin{figure*}[t]
    \centering 
    \includegraphics[width=2\columnwidth]{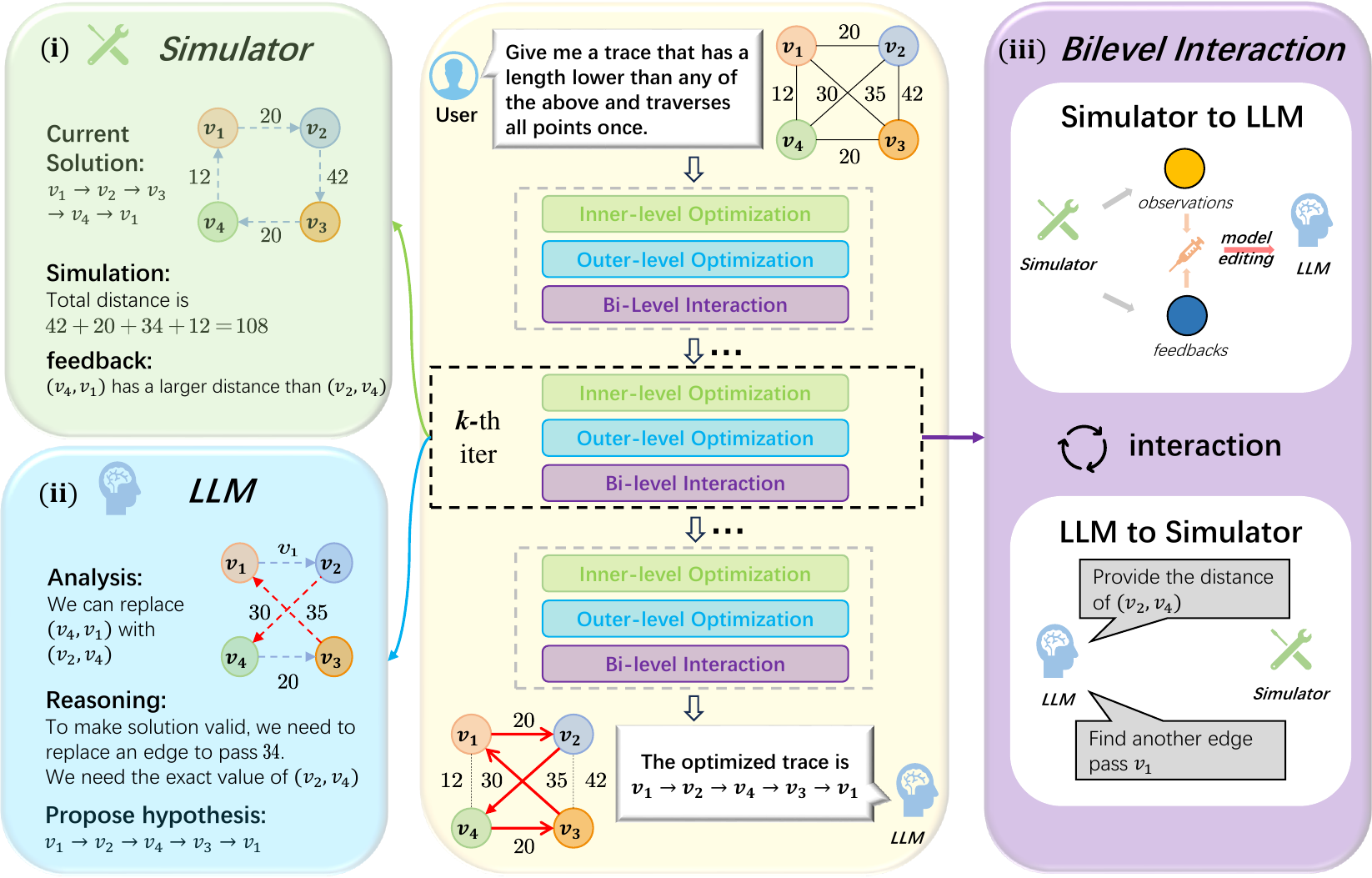}
    \caption{The overview of GSO. For a given optimization task, GSO iteratively conducts the \textcolor[HTML]{addb74}{\textbf{inner-level}} optimization, \textcolor[HTML]{1cb9ef}{\textbf{outer-level}} optimization, and \textcolor[HTML]{9567b1}{\textbf{bi-level}} interaction sequentially. The workflow is as follows: (i) the \textcolor[HTML]{addb74}{\textbf{inner-level}} simulator $\Phi$ conducts numerical simulations based on the current step's hypothetical solution $s_k$ ($v_1 \to v_2 \to v_3 \to v_4 \to v_1$) and returns observational feedback $f_k, \mathcal{L}_k$ (the edge ($v_4, v_1$) has a larger distance than the edge ($v_2, v_4$), current total distance: $108$); (ii) the \textcolor[HTML]{1cb9ef}{\textbf{outer-level}} LLM $\mathcal{M}_{\theta_k}$ generates new hypothetical solutions $s_{k+1}$ ($v_1 \to v_2 \to v_4 \to v_3 \to v_1$) based on the observational feedback $f_k, \mathcal{L}_k$; (iii) the \textcolor[HTML]{9567b1}{\textbf{bi-level}} interaction jointly updates simulations in conjunction with the expert knowledge within the LLMs through model editing. }
    \label{fig:overview}
\end{figure*}

\section{Preliminaries}

% \paragraph{Notations} 

\paragraph{Model Editing} Model editing mainly focuses on altering the internal knowledge within LLMs. It aims to update a range of intricate learned concepts, such as logical reasoning, spatial awareness, and numerical understanding, to make tailored adjustments to the model's behavior. In this paper, we follow \cite{assess, ME_survey} to update factual knowledge represented as (subject \emph{s}, relation \emph{r}, object \emph{o}). An LLM is expected to retrieve a memory corresponding to $o$ when given a natural language prompt. Editing a fact involves replacing the existing knowledge triple $(s, r, o)$ with a new one $(s, r, o^{*})$. For simplicity, an edit is denoted as $e = (s, r, o, o^{*})$. Given a set of edits $\mathcal{E}=\left\{e_{1}, e_{2}, \ldots\right\}$ and an original model $f_{\theta_0}$, sequential model editing applies each edit consecutively, i.e., $\mathcal{F}(f_{\theta_{n-1}}, e_n) = f_{\theta_n}$, where $f_{\theta_n}$ refers to the model after $n$ edits.

\section{Related Work} \label{sec:realted_work}
% We present related works on LLMs for scientific optimization in the main text. More related works on LLMs and model editing are in Appendix \ref{app:realted_work}.

% 1014lty model editing这部分的参数引入的很奇怪
\paragraph{Model Editing} Model editing involves modifying the memorized knowledge contained in LLMs. Various kinds of complex learned beliefs such as logical, spatial, or numerical knowledge are expected to be edited.
% Based on whether model parameters are modified, existing editing methods can be categorized into parameter-modifying \cite{edit2,edit3,edit4} and parameter-preserving \cite{edit5,pre_edit2} approaches.
Many studies explore the role of MLP layers in Transformers, revealing that these layers store knowledge, which can be localized to specific neurons and edited \cite{ana1,ana2,ana3}. KE \cite{ke} and MEND \cite{edit4} train a hypernetwork to compute gradient adjustments for updating model parameters. ROME \cite{edit2} and MEMIT \cite{edit3} apply the Locate-Then-Edit strategy, which first locates MLP storing factual knowledge, and then edits such knowledge by injecting a new key-value pair in the MLP module.

%  这一段放在preliminary
% In this paper, we integrate the parameter hypotheses generated by the outer-level LLM optimizer with the observational feedback from inner-level simulation platforms, editing into the LLM accordingly. Specifically,  we follow \cite{assess, ME_survey} to update factual knowledge represented as (subject \emph{s}, relation \emph{r}, object \emph{o}). An LLM is expected to retrieve a memory corresponding to $o$ when given a natural language prompt $p(s, r)$. Editing a fact involves replacing the existing knowledge triple $(s, r, o)$ with a new one $(s, r, o^{*})$. For simplicity, an edit is denoted as $e = (s, r, o, o^{*})$. Given a set of edits $\mathcal{E}=\left\{e_{1}, e_{2}, \ldots\right\}$ and an original model $f_{\theta_0}$, sequential model editing applies each edit consecutively, i.e., $\mathcal{F}(f_{\theta_{n-1}}, e_n) = f_{\theta_n}$, where $f_{\theta_n}$ refers to the model after $n$ edits. 

\paragraph{LLMs for Scientific Optimization}
A recent line of research explores methods to enhance LLM performance by incorporating natural language feedback as prompts to revise model outputs in improving reasoning \cite{reason1,reason2}, code generation \cite{code1,code3,code4}, dialogue applications \cite{diag1,diag2,restliu}, and so on \cite{soon1,soon2}. LLMs have also been demonstrated to optimize complex problems by utilizing tools \cite{LATM}. AlphaGeometry's \cite{geo} success in solving complex geometry problems without human demonstrations underscores the potential of LLMs in automating complex tasks. OPRO \cite{opro} employs LLMs as black-box optimizers for complex reasoning tasks. Eureka \cite{eureka} generates multiple solutions in each step to increase the success rate of produced codes. Funsearch \cite{funsearch} utilizes an evolutionary strategy to avoid local optima by using LLMs along with a systematic evaluator. SGA \cite{bi_level} introduces a bilevel optimization framework to enhance the knowledge-driven capabilities of LLMs by integrating simulations. 

\section{Method}

We first briefly outline the bi-level pipeline of our GSO, describing how the inner-level and the outer-level optimization are jointly conducted. 
Specifically, for a scientific optimization task $y$ (e.g., the traveling salesman problem), GSO aims to predict the optimal solution $\hat{s}$, by which GSO iteratively optimizes the solution $s_k$ from an initial solution $s_0$, where $k$ denotes the iteration step. 

Each iteration consists of three parts: the inner-level simulation platform, outer-level LLM optimization, and bi-level interactions. An overview of GSO is in Figure \ref{fig:overview}. 

The inner-level simulator \textbf{first} serves as an experimental platform, where a simulator $\Phi$ is employed to take a intermediate solution $s_k$ and its scientific expression $\mathcal{E}_k$ (e.g., coordinates of the nodes in the TSP task) as inputs, and outputs the corresponding observational feedback $f_k$ and the optimization objective $\mathcal{L}_k$ (e.g., the total distance of the current route in the TSP task):

\begin{equation*}
    (f_k, \mathcal{L}_k) = \Phi\left(s_k; \mathcal{E}_k\right)
\end{equation*}

The outer-level optimization \textbf{then} utilizes an LLM with its parameter $\theta_{k}$ (denote by $\mathcal{M}_{\theta_{k}}$) as a reasoning agent, propose a new hypothesis to optimize the current solution $s_k$:

\begin{equation*}
% \vspace{-2mm}
    \mathcal{E}_{k+1}, s_{k+1}=\mathcal{M}_{\theta_k}\left(\left\{\mathcal{L}_i, f_i, \mathcal{E}_i, s_i\right\}_{i=0}^{k}; \mathcal{P}\right). 
    % \vspace{-2mm}
\end{equation*}

Here $\mathcal{P}$ denotes the input prompt to LLM, and $\{\mathcal{L}_i, f_i, \mathcal{E}_i, s_i\}_{i=0}^{k}$ represents the historical optimization trajectory.

The bi-level interaction finally utilizes the feedback $f_k$ from inner-level simulation and the solutions $s_k$ provided by the outer-level LLM as a new key-value pair to update the LLM parameters $\theta$ through model editing:
\begin{equation*}
    \mathcal{M}_{\theta_{k+1}} = \mathcal{F}\left(\mathcal{M}_{\theta_{k}}, s_k, f_k \right), 
\end{equation*}

where $\theta_{k+1}$ represents the new parameter matrix obtained by applying a model editing step to the previous parameters $\theta_{k}$, $\mathcal{M}_{\theta_{k+1}}$ represents the according updated LLM, and $\mathcal{F}$ denotes the model editing process. We can then define the overall optimization task:

\begin{subequations}
\label{eq:bilevel}
\begin{align} 
    \min_{s}~&\mathcal{L}\Big(y(\mathcal{E}, s, \Phi)\Big)  \\ 
    \text{s.t.}~& G(\mathcal{E}, s, \Phi) \leq 0 
\end{align}
\end{subequations}

where $G\left(\cdot\right)\leq 0$ indicates that the current solution satisfies the constraints of the given optimization problem (for example, in the TSP task, the solution must visit all points and return to the starting point).

\subsection{Inner-level Simulation Platforms} 

{The inner-level optimization conducts numerical simulations for a hypothetical solution.
 
Simulations can provide domain-specific knowledge, transferring information from the simulation to optimization outputs (e.g., feedback in the form of total distance for the TSP task).

These outputs, paired with the predicted solution ${s}$, are then fed back into the LLMs to iteratively refine the hypothesis. 

The feedback can include the simulation loss relative to the target metric $\mathcal{L}$, as well as auxiliary information during the optimization process, which can provide more insights for improving solutions. 

For instance, if $\mathcal{L}$ represents the total distance of a solution for a TSP task, the auxiliary information may include the specific distance between two given points or the relative magnitude of distances between pairs of points (e.g. simulator can provide auxiliary information that $(v_4, v_1)$ has a larger distance than $(v_2, v_4)$ in  Figure \ref{fig:overview}).

\subsection{Outer-level LLM Optimization}

Many studies have demonstrated that LLMs are capable of discovering high-quality solutions and can match or even outperform hand-designed heuristic algorithms \cite{opro, many_na}. Therefore, we design outer-level optimization to:
\textbf{(i) }analyze the current optimization task and make reasonable hypotheses;
\textbf{(ii)} utilize the feedback from the simulations (including experimental phenomena and loss values, etc.), and update the LLM's knowledge through model editing; \textbf{(iii) } design and propose new potential solutions and input into the simulation platforms.

Inspired by the scientific optimization process by human scientists, we observe that when confronted with a new problem, scientists tend to be \textit{daring adventurers}, utilizing expert knowledge and problem parameters to thoroughly explore the search space, eliminating improbable parameter regions. This may significantly reduce the subsequent solution space, yielding solutions that satisfy optimization constraints, though they may not be optimal \cite{inspired}.

 Then, 
after gaining a deeper understanding of the problem through initial attempts, 
they change to be \textit{cautious followers} to keep the trajectory and trail previous solutions, adhering more closely to previous optimizations. This allows them to continually improve the metrics while satisfying optimization constraints, ultimately approaching the optimal solution. We thus design a\textbf{ dynamic exploitation and exploration strategy} to mimic this process.

Specifically, we calculate the loss change as $\Delta L = \frac{L_{\text{prev}} - L_{\text{cur}}}{L_{\text{prev}}}$. If $\Delta L > 0$, this indicates that the current hypothetical solutions are better than the previous ones. We then apply a lower decoding temperature to allow the LLM to follow the trace, i.e., $T_{\text{cur}} = T_{\text{pre}} \times \left( \frac{1}{1 + \Delta L} \right)$; if $\Delta L < 0$, this indicates that current solutions are worse than the previous ones. We thus apply a higher decoding temperature to be more exploratory and adventurous, i.e., $T_{\text{cur}} = T_{\text{pre}} \times (1 + |\Delta L|)$. We also clip the temperature to $[0, 1]$ in case that $T_{\text{cur}} < 0$ or $T_{\text{cur}} > 1$.

Empirically, we divide solutions into $\mathcal{S}_\text{exploit}$ and $\mathcal{S}_\text{explore}$. We observe that (i) $\mathcal{S}_\text{exploit}$ often contains repetitive solutions from previous iterations, and (ii) $\mathcal{S}_\text{explore}$ tends to yield solutions to be informative for guiding optimization, or otherwise infeasible.

\subsection{Bi-level Interactions}

The key challenge in integrating the two levels of optimization is developing a protocol that facilitates efficient, structured, and flexible communication between them. 

Therefore, we propose an interaction strategy considering the two aspects.
\textbf{On the one hand}, \textbf{ from simulation platforms to LLMs}, drawing upon \cite{edit2}, we incorporate the feedback from simulation platforms as new knowledge edited into the model. We \textbf{first} locate how facts are stored within the parameters of LLMs. 
We begin by analyzing and identifying which specific layers and their parameters ${W}$ have the strongest causal effect on predictions of individual facts through \textbf{causal tracing} \cite{edit2} (Detailed steps and results of the causal tracing for different LLMs are in Appendix \ref{app:casual_trace}). 
\textbf{Next}, we integrate the experimental results from the simulation platforms as new knowledge, formatted as triples, into the model. 
For instance, in the context of the TSP problem, suppose the simulation platform's feedback template includes: (i) the path length for the solution \( v_1 \to v_2 \to v_3 \to v_4 \) is $108$, (ii) the distance between \( v_1 \) and \( v_4 \) is greater than that between \( v_2 \) and \( v_4 \), and etc. These statements are then transformed into triples, such as \( (v_1 \to v_2 \to v_3 \to v_4, \hspace{1mm} \text{has distance of}, 108) \) and \( (v_1, v_4), \hspace{1mm} \text{greater than}, \hspace{1mm} (v_2, v_4)) \), which are subsequently edited into the model. \textbf{Finally}, 
by collecting such triples \( (s, r, o) \), we compile a series of key-value pairs for a set of vector keys \( S = [ s_1, s_2, \dots, s_n] \) and corresponding vector values \( O = [o_1, o_2, \dots, o_n] \) in each iteration step. Therefore, for the \( u \) new \( (s, r, o) \) pairs obtained from the next iteration step, the objective of editing can be given by the following formula:

\begin{small}
\begin{align*}
    W_1 \triangleq \argmin_{\hat{W}} \left( \sum_{i=1}^{n} \left\lVert \hat{W} s_i - o_i \right\rVert^2 + \sum_{i=n+1}^{n+u} \left\lVert \hat{W} s_{i} - o_{i} \right\rVert^2 \right)
\end{align*}
\end{small}

Here, $W$ denotes the original matrix and $W_1$ denotes the edited weight matrix. 

Due to this straightforward algebraic structure, any fact can be inserted directly once $(s, o)$ is determined. We can update the parameters of the LLM by editing the newly obtained \( u \) parameter-feedback pairs all at once, which adaptively utilizes the feedback from each simulation platform interaction into the LLM. It directly updates the LLM’s weight parameters and prevents the input prompt length from increasing with iterations, alleviating the loss-in-the-middle issue during long-text and multi-round optimization processes. 

\textbf{On the other hand, from LLMs to simulation platforms}, we leverage the LLM as a domain expert to guide the setup of each subsequent simulation based on the results of the previous results. This process is akin to an experienced scientist providing tailored guidance for the experimental configuration of the next step simulation based on the current results. The simulation platforms can then continue running experiments based on the guidance of the LLM to provide observational feedback. Specifically, each time the outer-level LLM receives feedback from the previous iteration of simulation platforms, it can simultaneously request intermediate results from the inner-level platforms to aid in its reasoning. Taking the TSP problem as an example, when the LLM provides the solution \( v_1 \to v_2 \to v_3 \to v_4 \to v_1 \), it can also request the platforms to provide the specific distance between two points. This additional feedback can then assist in refining the reasoning for the next iteration.

\begin{table*}[t]
\caption{Results of our GSO against \textbf{6} baselines using GPT-J 6B, Llama3 8B, and Mistral 7B as three representative backbone models (for more results of different backbone models, please see Appendix \ref{app:backbone}). Our experiments encompass \textbf{7} different tasks, which are divided into linear system regression (LSR) \textbf{(a)}, travel salesman problem (TSP) \textbf{(b)}, constitutive law prediction \textbf{(c-d)}, and molecule property prediction \textbf{(e-g)}. 
We report the mean $\pm$ standard error of each optimization result. The symbol N/A indicates that the model cannot provide a feasible solution for the current task. \textbf{A lower value is preferable across all tasks.} The best results are highlighted in \textbf{bold} text.}\label{tab:main_res}

\resizebox{\linewidth}{!}{
\begin{tabular}{llccccccc}
\toprule
\multirow{3}{*}{\textbf{Backbone}} & \multirow{3}{*}{\textbf{Method}} & \textbf{Linear System} & \textbf{Travel Salesman} & \multicolumn{2}{c}{\textbf{Constitutive Law}} & \multicolumn{3}{c}{\textbf{Molecule Property}} \\ 
\cmidrule(lr){3-3} \cmidrule(lr){4-4} \cmidrule(lr){5-6} \cmidrule(lr){7-9}
 & & \textbf{(a)} $\downarrow$ & \textbf{(b)}  $\downarrow$ & \textbf{(c)} $\downarrow$ & \textbf{(d)} $\downarrow$ & \textbf{(e)} $\downarrow$ & \textbf{(f)} $\downarrow$ & \textbf{(g)} $\downarrow$ \\
\midrule
\multirow{7}{*}{\textbf{GPT-J 6B}} & \textbf{Vanilla} &
N/A & 6.0 $\pm$ 2.0 & N/A & N/A & N/A & N/A & N/A \\
 & \textbf{CoT} &
N/A & 4.4 $\pm$ 1.9 & N/A & N/A & N/A & N/A & N/A \\
 & \textbf{Funsearch} &
29.7 $\pm$ 10.2 & 0.9 $\pm$ 0.1 & 77.4 $\pm$ 19.5 & {58.7 $\pm$ 17.2} & 230.8 $\pm$ 30.0 & 301.6 $\pm$ 51.9 & 34.9 $\pm$ 6.8 \\
 & \textbf{Eureka} &
45.1 $\pm$ 19.8 & 1.9 $\pm$ 0.8 & 155.3 $\pm$ 31.2 & 91.7 $\pm$ 18.5 & 390.0 $\pm$ 33.8 & 388.0 $\pm$ 70.3 & 55.3 $\pm$ 8.4 \\
 & \textbf{OPRO} &
41.0 $\pm$ 18.0 & 0.4 $\pm$ 0.2 & 60.9 $\pm$ 22.6 & 81.3 $\pm$ 10.2 & 455.9 $\pm$ 83.3 & 155.8 $\pm$ 18.9 & 23.7 $\pm$ 4.0 \\
 & \textbf{SGA} &
23.9 $\pm$ 5.0 & 0.2 $\pm$ 0.1 & 84.7 $\pm$ 17.9 & 73.2 $\pm$ 21.0 & 238.1 $\pm$ 17.5 & 89.5 $\pm$ 13.0 & 15.5 $\pm$ 3.1 \\
 & \textbf{GSO (ours)} &
\textbf{12.3 $\pm$ 4.8} & \textbf{0.0 $\pm$ 0.0} & \textbf{15.7 $\pm$ 5.0} & \textbf{55.4 $\pm$ 8.0} & \textbf{70.1 $\pm$ 17.7} & \textbf{83.1 $\pm$ 13.9} & \textbf{8.5 $\pm$ 2.0} \\
\bottomrule
\end{tabular}
}

\vspace{2mm}

\resizebox{\linewidth}{!}{
\begin{tabular}{llccccccc}
\toprule
\multirow{3}{*}{\textbf{Backbone}} & \multirow{3}{*}{\textbf{Method}} & \textbf{Linear System} & \textbf{Travel Salesman} & \multicolumn{2}{c}{\textbf{Constitutive Law}} & \multicolumn{3}{c}{\textbf{Molecule Property}} \\ 
\cmidrule(lr){3-3} \cmidrule(lr){4-4} \cmidrule(lr){5-6} \cmidrule(lr){7-9}
 & & \textbf{(a)} $\downarrow$ & \textbf{(b)}  $\downarrow$ & \textbf{(c)} $\downarrow$ & \textbf{(d)} $\downarrow$ & \textbf{(e)} $\downarrow$ & \textbf{(f)} $\downarrow$ & \textbf{(g)} $\downarrow$ \\
\midrule
\multirow{7}{*}{\textbf{Llama3 8B}} & \textbf{Vanilla} &
N/A & 6.0 $\pm$ 2.0 & N/A & N/A & N/A  & N/A  & N/A  \\
 & \textbf{CoT} &
N/A & 4.4 $\pm$ 1.9 & 3397.0 $\pm$ 298.9 & N/A  & N/A  & N/A  & N/A  \\
 & \textbf{Funsearch} &
15.7 $\pm$ 8.3 & 1.3 $\pm$ 0.5 & 198.9 $\pm$ 30.8 & 335.9 $\pm$ 112.1 & 433.7 $\pm$ 25.5 & 175.5 $\pm$ 28.8 & 91.4 $\pm$ 19.8 \\
 & \textbf{Eureka} &
18.3 $\pm$ 5.1 & 2.1 $\pm$ 0.2 & 255.6 $\pm$ 80.9 & 401.3 $\pm$ 98.0 & 260.7 $\pm$ 38.9 & 130.1 $\pm$ 15.9 & 155.1 $\pm$ 28.3 \\
 & \textbf{OPRO} &
17.0 $\pm$ 4.7 & 0.3 $\pm$ 0.1 & 74.1 $\pm$ 10.3 & 227.8 $\pm$ 33.2 & 537.1 $\pm$ 69.7 & 115.5 $\pm$ 35.3 & 51.9 $\pm$ 10.4 \\
 & \textbf{SGA} &
10.2 $\pm$ 3.7 & 0.5 $\pm$ 0.1 & 89.1 $\pm$ 5.9 & 150.0 $\pm$ 40.3 & 235.1 $\pm$ 55.7 & 94.1 $\pm$ 22.0 & 30.5 $\pm$ 15.3 \\
 & \textbf{GSO (ours)} &
\textbf{5.1 $\pm$ 1.0} & \textbf{0.0 $\pm$ 0.0} & \textbf{8.1 $\pm$ 2.1} & \textbf{20.1 $\pm$ 3.9} & \textbf{30.1 $\pm$ 14.9} & \textbf{20.9 $\pm$ 9.3} & \textbf{9.7 $\pm$ 3.6} \\
\bottomrule
\end{tabular}
}

\vspace{2mm}

\resizebox{\linewidth}{!}{
\begin{tabular}{llccccccc}
\toprule
\multirow{3}{*}{\textbf{Backbone}} & \multirow{3}{*}{\textbf{Method}} & \textbf{Linear System} & \textbf{Travel Salesman} & \multicolumn{2}{c}{\textbf{Constitutive Law}} & \multicolumn{3}{c}{\textbf{Molecule Property}} \\ 
\cmidrule(lr){3-3} \cmidrule(lr){4-4} \cmidrule(lr){5-6} \cmidrule(lr){7-9}
 & & \textbf{(a)} $\downarrow$ & \textbf{(b)}  $\downarrow$ & \textbf{(c)} $\downarrow$ & \textbf{(d)} $\downarrow$ & \textbf{(e)} $\downarrow$ & \textbf{(f)} $\downarrow$ & \textbf{(g)} $\downarrow$ \\
\midrule
\multirow{7}{*}{\textbf{Mistral 7B}} & \textbf{Vanilla} &
N/A & 6.0 $\pm$ 2.0 & N/A & N/A & N/A & N/A & N/A \\
 & \textbf{CoT} &
N/A & 4.4 $\pm$ 1.9 & N/A & N/A & N/A & N/A & N/A \\
 & \textbf{Funsearch} &
67.3 $\pm$ 14.7 & 1.3 $\pm$ 0.3 & 317.9 $\pm$ 43.0 & 199.4 $\pm$ 31.0 & 335.9 $\pm$ 61.3 & 310.4 $\pm$ 33.1 & 86.8 $\pm$ 17.4 \\
 & \textbf{Eureka} &
37.9 $\pm$ 13.1 & 2.5 $\pm$ 0.8 & 337.3 $\pm$ 51.1 & 379.0 $\pm$ 139.1 & 201.3 $\pm$ 33.7 & 331.7 $\pm$ 25.4 & 110.4 $\pm$ 19.5 \\
 & \textbf{OPRO} &
N/A & 0.3 $\pm$ 0.1 & 88.6 $\pm$ 12.1 & 275.7 $\pm$ 54.1 & 55.7 $\pm$ 13.9 & 252.4 $\pm$ 31.7 & 12.5 $\pm$ 3.9 \\
 & \textbf{SGA} &
27.9 $\pm$ 5.8 & 0.3 $\pm$ 0.1 & 145.1 $\pm$ 11.7 & 227.3 $\pm$ 19.9 & 130.5 $\pm$ 33.1 & 55.7 $\pm$ 12.1 & 55.4 $\pm$ 18.4 \\
 & \textbf{GSO (ours)} &
\textbf{5.6 $\pm$ 3.0} & \textbf{0.1 $\pm$ 0.0} & \textbf{3.2 $\pm$ 1.7} & \textbf{13.0 $\pm$ 2.1} & \textbf{23.9 $\pm$ 5.1} & \textbf{2.8 $\pm$ 1.5} & \textbf{1.7 $\pm$ 0.4} \\
\bottomrule
\end{tabular}
}

\end{table*}

\section{Experiments}

\subsection{Problem Descriptions} 
We provide a brief definition of each task. More detailed definitions are in Appendix \ref{app:def}. For the \textbf{linear system regression} task, {the objective is to estimate the linear coefficients to model relationships between input and corresponding output} ~\cite{fisher1922}. We use \textbf{"\#steps"} (optimization steps for successfully finding the optimal solution) \textbf{as the metric}. We follow the same dataset setting as OPRO \cite{opro}. For the \textbf{TSP} task, the objective is that given a set of \(n\) nodes with known coordinates, it seeks to find the shortest possible route that visits each node once and returns to the starting point. We follow the same dataset setting as OPRO \cite{opro} and use the Gurobi solver \cite{gurobi} to construct the oracle solutions and compute the \textbf{optimality gap as the metric} (the difference between the distance in the solution by the evaluated approach and by the oracle solution, divided by the distance of the oracle solution). For the \textbf{constitutive law prediction} task, we consider both fitting linear and non-linear materials and follow the same dataset setting as SGA \cite{bi_level}.  We use \textbf{mean square error (MSE) as the metric}. For the \textbf{molecular property prediction} task, we consider three tasks: predict a molecule's highest occupied molecular orbital (HOMO), lowest unoccupied molecular orbital (LUMO), and the HOMO-LUMO gap based on their conformations and quantum mechanical properties. We follow SGA to use the QM9 dataset \cite{qm9} for experiments. We also use \textbf{MSE as the metric}.

\subsection{Experiment Setups}

\begin{table*}[t]
\caption{The results of the ablation study of our GSO on the seven scientific optimization tasks, using Llama3 8B as the backbone model (We provide more results for the other \textit{five} backbone models, in  Appendix \ref{app:abla}).}\label{tab:abla}

\resizebox{\linewidth}{!}{
\begin{tabular}{lccccccc}
\toprule
\multirow{2}{*}{\textbf{Method}} & \textbf{Linear System} & \textbf{Travel Salesman} & \multicolumn{2}{c}{\textbf{Constitutive Law}} & \multicolumn{3}{c}{\textbf{Molecule Property}} \\ \cmidrule(lr){2-2} \cmidrule(lr){3-3} \cmidrule(lr){4-5} \cmidrule(lr){6-8}
 & \textbf{(a)} $\downarrow$ & \textbf{(b)}  $\downarrow$ & \textbf{(c)} $\downarrow$ & \textbf{(d)} $\downarrow$ & \textbf{(e)} $\downarrow$ & \textbf{(f)} $\downarrow$ & \textbf{(g)} $\downarrow$ \\
\midrule
\textbf{GSO$_{w/o \hspace{1mm} edit}$} &
33.7 $\pm$ 10.1 & 0.7 $\pm$ 0.3 & 116.1 $\pm$ 30.1 & 141.5 $\pm$ 36.9 & 307.7 $\pm$ 31.4 & 331.2 $\pm$ 89.2 & 164.3 $\pm$ 38.9 \\
\textbf{GSO$_{w/o \hspace{1mm} dynamic}$} &
6.5 $\pm$ 3.5 & 0.0 $\pm$ 0.0 & 23.7 $\pm$ 9.9 & 76.9 $\pm$ 11.4 & 35.0 $\pm$ 17.1 & 45.0 $\pm$ 19.3 & 21.9 $\pm$ 7.3 \\
\midrule
\textbf{GSO} &
\textbf{5.1 $\pm$ 1.0} & \textbf{0.0 $\pm$ 0.0} & \textbf{8.1 $\pm$ 2.1} & \textbf{20.1 $\pm$ 3.9} & \textbf{30.1 $\pm$ 14.9} & \textbf{20.9 $\pm$ 9.3} & \textbf{9.7 $\pm$ 3.6} \\
\bottomrule
\end{tabular}
}
\end{table*}

\noindent \textbf{Implementation Details. } We apply LLMs including GPT-J 6B \cite{gpt6b}, Llama3 8B \cite{llama3}, Mistral 7B \cite{mistral}, Llama2 13B \cite{llama2}, Yi9b \cite{yi9b}, and Internlm 7B \cite{internlm}. We follow \cite{eureka} to conduct all experiments five times using different random seeds to guarantee stable and reproducible results. More details of implementation details are in Appendix \ref{app:def}.

We consider \textbf{six} strong baselines for evaluation: (i) vanilla LLMs without additional modules. Vanilla LLMs represent the original capabilities of LLMs. (ii) \textbf{CoT} prompting~\cite{wei2022chain} solves the problem by looking at step-by-step solutions from examples. We provide $5$ examples with explanations as the initial solutions. (iii) \textbf{Funsearch}~\cite{funsearch} utilizes evolutionary strategy to avoid local optimum. (iv) \textbf{Eureka}~\citep{eureka} generates multiple solutions in each iteration to improve the success rate. (v) \textbf{OPRO}~\cite{opro} highlights the advantages of involving a sorted optimization trajectory. (vi) \textbf{SGA}~\cite{bi_level} utilizes a top-k heap to generate more diverse solutions.

\begin{figure}[t]
    \centering 
    \includegraphics[width=1\columnwidth]{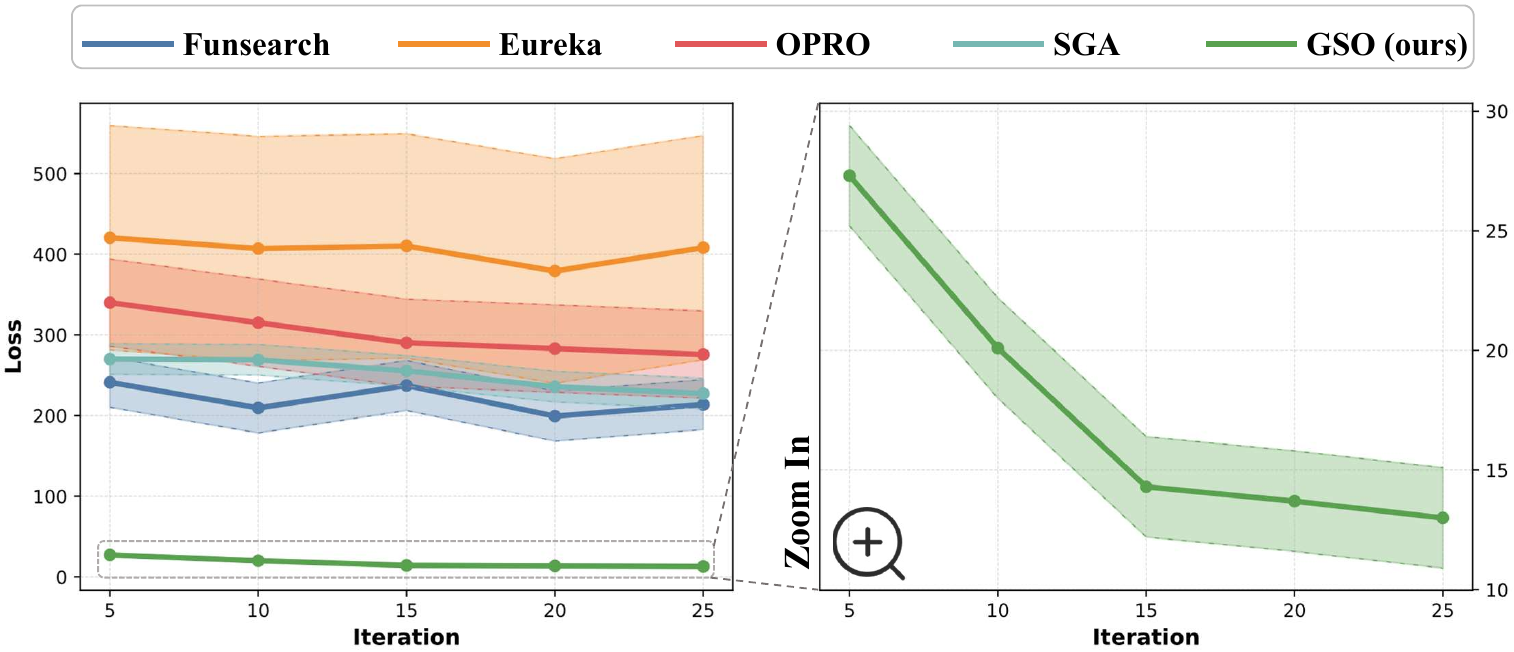}
    \caption{We visualize the average MSE loss values of each method for the non-linear Constitutive Law task (d) across five random seeds at the same optimization steps using Mistral 7B as the backbone model, with shading representing the standard deviation.}
    \label{fig:loss_curve}
    % \vspace{-5mm}
\end{figure}

\begin{table*}[t]
\caption{We evaluate our GSO against other baseline methods on five prompts that vary in format but are semantically similar, and report their average results of each method, using Mistral 7B as the representative backbone model. The symbol N/A indicates that the model cannot provide a feasible solution for the current task. }\label{tab:prompt_robust}

\resizebox{\linewidth}{!}{
\begin{tabular}{llccccccc}
\toprule
\multirow{3}{*}{\textbf{Backbone}} & \multirow{3}{*}{\textbf{Method}} & \textbf{Linear System} & \textbf{Travel Salesman} & \multicolumn{2}{c}{\textbf{Constitutive Law}} & \multicolumn{3}{c}{\textbf{Molecule Property}} \\ 
\cmidrule(lr){3-3} \cmidrule(lr){4-4} \cmidrule(lr){5-6} \cmidrule(lr){7-9}
 & & \textbf{(a)} $\downarrow$ & \textbf{(b)}  $\downarrow$ & \textbf{(c)} $\downarrow$ & \textbf{(d)} $\downarrow$ & \textbf{(e)} $\downarrow$ & \textbf{(f)} $\downarrow$ & \textbf{(g)} $\downarrow$ \\
\midrule

% & \textbf{Vanilla} & N/A & --- & N/A & N/A & N/A & N/A & N/A \\
%  & \textbf{CoT}  & N/A & ---  & ---  & ---  & ---  & --- & --- \\
 \multirow{5}{*}{\textbf{Mistral 7B}}  & \textbf{Funsearch}  & 52.2 $\pm$ 29.4 & 2.4 $\pm$ 1.5  & 411.0 $\pm$ 193.4  & 254.9 $\pm$ 109.3  & 473.3 $\pm$ 130.2 & 403.1 $\pm$ 135.0 & 173.0 $\pm$ 66.3 \\
 & \textbf{Eureka}  & 60.25 $\pm$ 29.2 & 3.5 $\pm$ 2.6  & 303.7 $\pm$ 180.3  & 179.0 $\pm$ 65.8  & 219.0 $\pm$ 40.6  & 287.9 $\pm$ 99.0 & 148.9 $\pm$ 55.0 \\
 & \textbf{OPRO}  & N/A & 0.4 $\pm$ 0.3  & 107.6 $\pm$ 61.4  & 224.8 $\pm$ 34.2  & 94.3 $\pm$ 50.6  & 593.7 $\pm$ 266.8 & 40.5 $\pm$ 28.0\\
 & \textbf{SGA}  & 38.7 $\pm$ 20.3  & 0.3 $\pm$ 0.2  & 55.4 $\pm$ 43.9  & 203.2 $\pm$ 22.3  & 211.7 $\pm$ 98.0  & 73.0 $\pm$ 18.8 & 62.9 $\pm$ 33.2 \\
 & \textbf{GSO (ours)}  & \textbf{8.3 $\pm$ 3.3} & \textbf{0.1 $\pm$ 0.0}  & \textbf{5.5 $\pm$ 3.1}  & \textbf{20.9 $\pm$ 8.3} & \textbf{18.8 $\pm$ 7.0}  & \textbf{3.9 $\pm$ 2.1} & \textbf{2.2 $\pm$ 0.4} \\
\bottomrule
\end{tabular}
}
\end{table*}

\begin{table*}[t]
\caption{Results of our GSO against popular closed-source LLMs, including GPT-4o and Claude-3.5. We present the results of our GSO Mistral 7B as the backbone model for comparison. }\label{tab:comp_with_close}

\resizebox{\linewidth}{!}{
\begin{tabular}{lccccccc}
\toprule
\multirow{2}{*}{\textbf{Method}} & \textbf{Linear System} & \textbf{Travel Salesman} & \multicolumn{2}{c}{\textbf{Constitutive Law}} & \multicolumn{3}{c}{\textbf{Molecule Property}} \\ \cmidrule(lr){2-2} \cmidrule(lr){3-3} \cmidrule(lr){4-5} \cmidrule(lr){6-8}
 & \textbf{(a)} $\downarrow$ & \textbf{(b)}  $\downarrow$ & \textbf{(c)} $\downarrow$ & \textbf{(d)} $\downarrow$ & \textbf{(e)} $\downarrow$ & \textbf{(f)} $\downarrow$ & \textbf{(g)} $\downarrow$ \\
\midrule
\textbf{GPT-4o} & 6.4 $\pm$ 1.5 & 0.2 $\pm$ 0.0 & 59.7 $\pm$ 7.4  & 45.5 $\pm$ 9.1 & 181.4 $\pm$ 33.7 & 313.5 $\pm$ 45.4 & 67.8 $\pm$ 10.9  \\
% \midrule
\textbf{Claude-3} & 12.1 $\pm$ 2.0 & 0.2 $\pm$ 0.1 & 128.1 $\pm$ 13.5 & 103.0 $\pm$ 5.2 & 162.3 $\pm$ 21.2 & 519.3 $\pm$ 37.4 & 37.8 $\pm$ 8.3  \\
\midrule
\textbf{GSO} & \textbf{5.6 $\pm$ 3.0} & \textbf{0.1 $\pm$ 0.0} & \textbf{3.2 $\pm$ 1.7} & \textbf{13.0 $\pm$ 2.1} & \textbf{23.9 $\pm$ 5.1} & \textbf{2.8 $\pm$ 1.5} & \textbf{1.7 $\pm$ 0.4} \\
\bottomrule
\end{tabular}
}
\end{table*}

\subsection{Main Results} \label{sec:main_res}

We conduct our experiments on the $7$ designed tasks using GPT-J 6B, Llama3 8B, and Mistral 7B as three representative backbone models in Table \ref{tab:main_res}. We also provide more results for other \textit{three} different backbone models in Table \ref{tab:main_res_housange} in Appendix \ref{app:backbone} to demonstrate the versatility of our GSO across different backbones. GSO enables tasks, which are challenging to effectively optimize with traditional Vanilla and CoT methods, to become feasible. We also observe that GSO significantly and consistently outperforms existing methods on the scientific optimization tasks, which demonstrates the effectiveness of our GSO. Notably, for the molecule property prediction task on predicting the HOMO-LUMO gap, GSO achieves a maximum precision improvement of over 32.6$\times$ when utilizing Mistral 7B as the backbone model. These results underscore the importance of effectively utilizing the observational feedback to adaptively adjust its optimization directions. The universality of \textit{six} popular open-source backbone models in Appendix \ref{app:backbone} also suggests the generalizability of our GSO.

\subsection{Ablation Study}\label{sec:abla}
To further investigate the contribution of each component within GSO, we conduct a series of ablation experiments on the entire framework. Specifically, we denote GSO without \textit{edit} as GSO$_{w/o \hspace{1mm} edit}$, GSO without the dynamic strategy as GSO$_{w/o \hspace{1mm} dynamic}$, respectively. We use Llama3 8B as the backbone model, the results of other backbone models are in Appendix \ref{app:abla}. We present the ablation results of GSO using Llama 3 8B as the backbone model in Table \ref{tab:abla}. More Results using the other five different backbone models are in Appendix \ref{app:abla}.
As shown in Table \ref{tab:abla}, the absence of any component within GSO results in a performance
degradation of the entire framework. Notably, GSO$_{w/o \hspace{1mm} edit}$ exhibits more significant impacts on the performance of GSO,  which demonstrates the importance of effectively utilizing the observational feedback to adaptively adjust the optimization direction.

\subsection{Case Study} \label{sec:case_study}

\paragraph{Robustness of the Prompt} As mentioned in Section \ref{sec:intro}, one appealing feature of our GSO compared to other methods is its robustness to prompts. Prompts with similar semantics do not require meticulous crafting to yield consistently promising results. To provide more insights into our GSO, we manually generated two prompts and used OpenAI o1 to rewrite an additional three based on the original task prompts for each task (details of the generated prompts are in Appendix \ref{augmented_prompt}). As shown in Table \ref{tab:prompt_robust}, we observe that prompt-based methods tend to be sensitive to the prompts, with semantically similar prompts often leading to fluctuating results across many tasks. This necessitates significant effort from users to perform prompt 'tuning' during application.
In contrast, \textbf{our GSO consistently produces robust results across different prompts, with minimal variation between semantically similar prompts.} This reduces the need for extensive prompt "tuning," highlighting its potential for broader real-world applications.

\paragraph{Loss Curve and Decoding Temperature} We also investigate the impact of the number of lengthy input prompts as optimization iterations increase on each method. We present the feedback loss trend for the non-linear constitutive law task (d) in Figure \ref{fig:loss_curve} as an example. We observe that GSO achieves significantly lower loss and exhibits a clear convergence trend compared to existing baselines. This demonstrates that GSO can effectively leverage feedback from each iteration to achieve stable and consistent improvements in scenarios involving lengthy prompts. Notably, some methods show improvement at the beginning (i.e., loss reduction) when the initial prompt is short and the optimization space is large, but as the number of optimization steps increases, the feedback loss fluctuates, making further optimization difficult. In contrast, GSO effectively leverages each step feedback to adaptively adjust its optimization direction, leading to a stable and consistent loss decrease. These results demonstrate that \textbf{GSO can effectively observational feedback from lengthy prompts to facilitate further optimization, thereby alleviating the loss in the middle issue.}

\paragraph{Comparison with Advanced Closed-source LLMs} We also compared our GSO framework on open-source models with the current state-of-the-art closed-source model, GPT-4o\footnote{https://openai.com/o1/} and Claude-3-Sonnet\footnote{https://www.anthropic.com/news/claude-3-5-sonnet}. As shown in Table \ref{tab:comp_with_close}, our GSO utilizing Mistral 7B can consistently outperform GPT-4o and Claude-3-Sonnet, demonstrating the effectiveness of our approach. Note that our method is orthogonal to the choice of backbone model, making it a versatile plug-and-play module that can be directly applied to more advanced LLMs to further achieve enhanced results.
\section{Conclusion}
In this paper, we propose a novel \textbf{G}eneral \textbf{S}cientific \textbf{O}ptimizers method, effectively enabling LLMs to utilize observational feedback from each optimization step.
Specifically, GSO consists of a bi-level optimization framework: outer-level LLMs function as knowledgeable and versatile scientists, generating new hypotheses to optimize experimental hyperparameters; inner-level simulations function as experimental platforms to perform numerical simulations to these hypotheses and provide observational feedback; a bi-level interaction then update the simulators together with the expert knowledge within LLMs via model editing. GSO effectively guides LLMs to derive a more precise and stable optimization direction, yielding superior optimization results. Extensive experiments on  
 \textit{six} different open-source and \textit{seven} different scientific tasks demonstrate the superiority of our GSO, delivering consistent, robust, generalizable, and nearly monotonic improvement. We view our GSO as a trailblazer, establishing a new paradigm for utilizing LLMs and simulations as bi-level optimization to further advancements in scientific optimizations.\footnote{More discussions on GSO are in Appendix \ref{app:more_discuss}.}
 
 % \clearpage

 \section*{Acknowledgements}

    The authors would like to thank all the anonymous reviewers for their insightful comments.
 
\section{Limitations}
We consider a few limitations and future directions. (i) The content generated by LLMs exhibits a certain degree of randomness, and the optimization process cannot guarantee interpretability or transparency. (ii) LLM inference requires large computational resources and thus increases expense.  It paves the way for research on LLM inference acceleration to expedite our GSO \cite{sd,sd2}.
(iii) GSO requires access to model weights, which limits its applicability to closed-source models like GPT-4 and Claude 3.5. We believe that as the community progresses and the performance gap between open-source and closed-source models narrows, our GSO will be able to demonstrate its capabilities more effectively.

% \clearpage
% \bibliography{custom}

\appendix
\clearpage
\section{More Related Works} \label{app:realted_work}

\subsection{Classical Language Model Methods for Scientific Optimization} 
% wh：LMs？
Besides the methods using LLMs for scientific optimization in Section \ref{sec:realted_work}, these are several classical studies that have fine-tuned language models (LMs) to function as mutation and crossover operators within evolutionary algorithms \cite{tune1, tune2, tune3}. LMX \cite{LMX} employs language models guided by few-shot exemplars to generate evolutionary crossovers in tasks like image and code generation, enhancing the model's ability to adapt and innovate across diverse domains. ELM \cite{ELM} trains an LLM to generate code diffs, which serve as the mutation operator, and introduces a fine-tuning method to enhance performance in the Sodarace domain for robotic simulation. EvoPrompting \cite{EvoPrompting} leverages large language models to evolve neural network architectures by integrating evolutionary search with soft prompt tuning. OptFormer \cite{OptFormer} incorporates trajectories as input for optimization by training a transformer model on extensive hyperparameter optimization datasets. These methods primarily focus on fine-tuning LLMs with domain-specific data to align natural language with scientific information.  However, these approaches are domain-bound and demand substantial data which also limits their broader applicability.

\subsection{Large Language Models}
Language models such as GPT \cite{gpt1}, BERT \cite{bert}, RoBERTa \cite{roberta}, and  Megatron-LM \cite{megantron} have led to a learning paradigm shift in natural language processing (NLP). Models are first pre-trained on extensive volumes of unlabeled text corpora with language modeling objectives and then fine-tuned on downstream tasks. Recently, large language models (LLMs) including LLama \cite{llama3, llama2} ChatGPT \cite{chatgpt} GPT4 \cite{gpt4}, PaLM \cite{palm}, Gemini \cite{gemini}, and Claude3 \cite{claude3} have shown great performance in both few-shot and even zero-shot scenarios \cite{few-shot}. 
To further enhance the interpretability of these LLMs, some research endeavors explain LLMs through attribution analysis \cite{104interpretability, greater, finding, sac-kg}. 
Another line of work aims to retrieve the knowledge explicitly from LLMs as the basis for interpreting them, including the reasoning task \cite{distract} and the QA task \cite{bertnet, qa1, qa2}.

\begin{table*}[t]
\caption{Results of our GSO against \textbf{6} baselines using Llama2 13B, Yi 9B, and Internlm 7B as the backbone models. Our experiments encompass \textbf{7} different tasks, which are divided into linear system regression (LSR) \textbf{(a)}, travel salesman problem (TSP) \textbf{(b)}, constitutive law prediction \textbf{(c-d)}, and molecule property prediction \textbf{(e-g)}. For the LSR task, we use the number of steps of successfully finding
the optimal solution as the metric. For the TSP task, we use the optimality gap as the metric. For the rest five
tasks, we use MSE loss as the metric. We calculate the mean $\pm$ standard error of each optimization result. The symbol N/A indicates that the model is unable to provide a feasible solution for the current task. A lower value is preferable across all tasks. The best results are highlighted in \textbf{bold} text.}\label{tab:main_res_housange}

\resizebox{\linewidth}{!}{
\begin{tabular}{llccccccc}
\toprule
\multirow{3}{*}{\textbf{Backbone}} & \multirow{3}{*}{\textbf{Method}} & \textbf{Linear System} & \textbf{Travel Salesman} & \multicolumn{2}{c}{\textbf{Constitutive Law}} & \multicolumn{3}{c}{\textbf{Molecule Property}} \\ 
\cmidrule(lr){3-3} \cmidrule(lr){4-4} \cmidrule(lr){5-6} \cmidrule(lr){7-9}
 & & \textbf{(a)} $\downarrow$ & \textbf{(b)}  $\downarrow$ & \textbf{(c)} $\downarrow$ & \textbf{(d)} $\downarrow$ & \textbf{(e)} $\downarrow$ & \textbf{(f)} $\downarrow$ & \textbf{(g)} $\downarrow$ \\
\midrule
\multirow{7}{*}{\textbf{Llama2 13B}} & \textbf{Vanilla} &
N/A & 6.0 $\pm$ 2.0 & N/A & N/A & N/A & N/A & N/A \\
 & \textbf{CoT} &
N/A & 4.4 $\pm$ 1.9 & N/A & N/A & N/A & N/A & N/A \\
 & \textbf{Funsearch} &
N/A & 2.5 $\pm$ 1.0 & 170.9 $\pm$ 20.1 & 255.1 $\pm$ 31.1 & 108.5 $\pm$ 13.0 & 115.5 $\pm$ 20.7 & 59.7 $\pm$ 8.3 \\
 & \textbf{Eureka} &
N/A & 2.7 $\pm$ 1.5 & 211.7 $\pm$ 50.4 & 131.7 $\pm$ 21.9 & 98.2 $\pm$ 11.3 & 220.1 $\pm$ 30.1 & 39.3 $\pm$ 13.9 \\
 & \textbf{OPRO} &
28.7 $\pm$ 10.8 & 0.9 $\pm$ 0.3 & 55.8 $\pm$ 10.4 & 165.8 $\pm$ 52.9 & 168.5 $\pm$ 31.7 & 355.4 $\pm$ 43.7 & 35.9 $\pm$ 10.7 \\
 & \textbf{SGA} &
30.0 $\pm$ 10.2 & 0.2 $\pm$ 0.1 & 31.7 $\pm$ 10.1 & 55.4 $\pm$ 15.8 & 87.1 $\pm$ 17.0 & 89.5 $\pm$ 17.7 & 27.6 $\pm$ 5.4 \\
 & \textbf{GSO (ours)} &
\textbf{5.0 $\pm$ 1.6} & \textbf{0.1 $\pm$ 0.1} & \textbf{2.9 $\pm$ 1.3} & \textbf{7.4 $\pm$ 3.1} & \textbf{48.1 $\pm$ 12.0} & \textbf{79.1 $\pm$ 13.3} & \textbf{8.3 $\pm$ 3.2} \\
\bottomrule
\end{tabular}
}

\vspace{2mm}

\resizebox{\linewidth}{!}{
\begin{tabular}{llccccccc}
\toprule
\multirow{3}{*}{\textbf{Backbone}} & \multirow{3}{*}{\textbf{Method}} & \textbf{Linear System} & \textbf{Travel Salesman} & \multicolumn{2}{c}{\textbf{Constitutive Law}} & \multicolumn{3}{c}{\textbf{Molecule Property}} \\ 
\cmidrule(lr){3-3} \cmidrule(lr){4-4} \cmidrule(lr){5-6} \cmidrule(lr){7-9}
 & & \textbf{(a)} $\downarrow$ & \textbf{(b)}  $\downarrow$ & \textbf{(c)} $\downarrow$ & \textbf{(d)} $\downarrow$ & \textbf{(e)} $\downarrow$ & \textbf{(f)} $\downarrow$ & \textbf{(g)} $\downarrow$ \\
\midrule
\multirow{7}{*}{\textbf{Yi 9B}} & \textbf{Vanilla} &
N/A & 6.0 $\pm$ 2.0 & N/A & N/A & N/A & N/A & N/A \\
 & \textbf{CoT} &
N/A & 4.4 $\pm$ 1.9 & N/A & N/A & N/A & N/A & N/A \\
 & \textbf{Funsearch} &
19.6 $\pm$ 6.0 & 2.2 $\pm$ 0.2 & 201.4 $\pm$ 10.3 & 229.5 $\pm$ 31.3 & 201.0 $\pm$ 55.3 & 198.9 $\pm$ 40.1 & 135.0 $\pm$ 20.7 \\
 & \textbf{Eureka} &
8.6 $\pm$ 4.0 & 1.8 $\pm$ 0.6 & 130.1 $\pm$ 22.1 & 381.1 $\pm$ 98.8 & 1559.1 $\pm$ 100.7 & 301.9 $\pm$ 38.7 & 98.3 $\pm$ 10.1 \\
 & \textbf{OPRO} &
9.5 $\pm$ 4.4 & 0.3 $\pm$ 0.2 & 94.0 $\pm$ 23.3 & 163.2 $\pm$ 28.9 & 850.7 $\pm$ 94.0 &  746.9 $\pm$ 31.0 & 129.1 $\pm$ 23.5 \\
 & \textbf{SGA} &
22.5 $\pm$ 5.1 & 0.7 $\pm$ 0.4 & 50.3 $\pm$ 7.9 & 104.4 $\pm$ 19.3 & 133.9 $\pm$ 23.1& \textbf{168.4 $\pm$ 59.0} & 39.6 $\pm$ 5.3 \\
 & \textbf{GSO (ours)} &
\textbf{3.0 $\pm$ 0.8} & \textbf{0.0 $\pm$ 0.0} & \textbf{5.9 $\pm$ 2.1} & \textbf{89.1 $\pm$ 33.9} & \textbf{67.9 $\pm$ 23.1} & {172.9 $\pm$ 43.1} & \textbf{5.5 $\pm$ 2.1} \\
\bottomrule
\end{tabular}
}

\vspace{2mm}

\resizebox{\linewidth}{!}{
\begin{tabular}{llccccccc}
\toprule
\multirow{3}{*}{\textbf{Backbone}} & \multirow{3}{*}{\textbf{Method}} & \textbf{Linear System} & \textbf{Travel Salesman} & \multicolumn{2}{c}{\textbf{Constitutive Law}} & \multicolumn{3}{c}{\textbf{Molecule Property}} \\ 
\cmidrule(lr){3-3} \cmidrule(lr){4-4} \cmidrule(lr){5-6} \cmidrule(lr){7-9}
 & & \textbf{(a)} $\downarrow$ & \textbf{(b)}  $\downarrow$ & \textbf{(c)} $\downarrow$ & \textbf{(d)} $\downarrow$ & \textbf{(e)} $\downarrow$ & \textbf{(f)} $\downarrow$ & \textbf{(g)} $\downarrow$ \\
\midrule
\multirow{7}{*}{\textbf{Internlm 7B}} & \textbf{Vanilla} &
N/A & 6.0 $\pm$ 2.0 & N/A & N/A & N/A & N/A & N/A \\
 & \textbf{CoT} &
N/A & 4.4 $\pm$ 1.9 & N/A & N/A & N/A & N/A & N/A \\
 & \textbf{Funsearch} &
19.1 $\pm$ 1.3 & 0.8 $\pm$ 0.2 & 193.8 $\pm$ 30.9 & 318.3 $\pm$ 87.1 & 110.4 $\pm$ 19.3 & 150.0 $\pm$ 23.7 & 70.3 $\pm$ 10.9 \\
 & \textbf{Eureka} &
29.5 $\pm$ 4.2 & 1.3 $\pm$ 0.4 & 211.0 $\pm$ 30.3 & 230.5 $\pm$ 39.9 & 211.1 $\pm$ 40.1 & 139.6 $\pm$ 38.2 & 50.8 $\pm$ 9.4 \\
 & \textbf{OPRO} &
27.6 $\pm$ 3.8 & 0.2 $\pm$ 0.1 & 59.6 $\pm$ 16.2 & 130.8 $\pm$ 19.1 & 257.9 $\pm$ 93.0 & 195.1 $\pm$ 33.9 & 45.3 $\pm$ 9.7 \\
 & \textbf{SGA} &
14.6 $\pm$ 5.5 & \textbf{0.0 $\pm$ 0.0} & 133.5 $\pm$ 29.3 & 191.4 $\pm$ 15.0 & 93.7 $\pm$ 25.4 & 94.8 $\pm$ 19.3 & 17.9 $\pm$ 4.5 \\
 & \textbf{GSO (ours)} &
\textbf{10.0 $\pm$ 2.1} & \textbf{0.0 $\pm$ 0.0} & \textbf{10.4 $\pm$ 3.1} & \textbf{37.0 $\pm$ 12.9} & \textbf{47.2 $\pm$ 10.8} & \textbf{73.5 $\pm$ 19.6} & \textbf{4.2 $\pm$ 2.3} \\
\bottomrule
\end{tabular}
}

\end{table*}

\begin{table}[ht]
\centering
\caption{We consider an imaginary constitutive law setting to prevent the LLM from cheating by memorization and report the results using Mistral 7B as the representative backbone model.}\label{tab:memory}
\resizebox{\linewidth}{!}{
\begin{tabular}{ccccc}
\toprule
  \textbf{Funsearch} & \textbf{Eureka} & \textbf{OPRO} & \textbf{SGA} & \textbf{GSO}  \\
\midrule
201.1 $\pm$ 39.0 & 291.0 $\pm$ 20.8 & 190.0 $\pm$ 10.3  & 59.7 $\pm$ 5.9 & \textbf{17.1 $\pm$ 4.0}  \\
\bottomrule
\end{tabular}}

\end{table}

\section{Results of the Causal Tracing for Different LLMs} \label{app:casual_trace}
Causal tracing has emerged as a pivotal methodology for dissecting and understanding the internal mechanisms of model \cite{cg_01}. This technique facilitates the identification and modification of specific factual associations within a model without necessitating comprehensive retraining \cite{casual_graph}. In the context of model editing, causal tracing enables precise interventions by isolating the neural correlates responsible for particular behaviors or outputs. We follow \cite{edit2} to build a \textbf{causal graph} \cite{casual_graph} to describe dependencies between the hidden variables. This graph illustrates numerous pathways from the input on the left to the output (next-word prediction) at the lower right. Our aim is to \textbf{determine whether specific hidden state variables are more important than others in the process of recalling a fact}. 

To quantify each state's contribution to a correct factual prediction, we analyze all of LLM's internal activations across three runs: \textbf{a clean run} that accurately predicts the fact,\textbf{ a corrupted run} where the prediction is impaired, and \textbf{a corrupted-with-restoration run} that evaluates the ability of a single state to restore the correct prediction.

Let $\mathbb{P}[o]$, $\mathbb{P}_*[o]$, and $\mathbb{P}_{*,\smash{\text{ clean } {h}_{i}^{l}}}[o]$ denote the probability of emitting the given entity $o$ under the clean, corrupted, and corrupted-with-restoration runs, respectively; dependence on the input $x$ is omitted for notational simplicity.

The \textbf{total effect} (TE) is defined as the difference between two probabilities: 

\[
\text{TE} = \mathbb{P}[o] - \mathbb{P}_*[o].
\]

The \textbf{indirect effect} (IE) of a specific mediating state \( \hat{h_i}^{l} \) is defined as the difference between the probability of \( o \) under the corrupted condition and the probability of \( o \) when that state is restored to its clean version, while the subject remains in a corrupted state:

\[
\text{IE} = \mathbb{P}_{*, \text{clean } {h_i}^{l}}[o] - \mathbb{P}_*[o].
\]

By averaging over a sample of statements, we can derive the \textbf{average indirect effect (AIE)} for each hidden state variable. Specifically, when restoring hidden states from the original run, we substitute the values computed originally for the corresponding layer and token, allowing subsequent computations to proceed without further modification. Taking Llama3 8B as an example, as shown in Figure \ref{fig:casual_llama3}, we observe that the $18\text{-}22$th layers for the last subject token demonstrate the most causality in the AIE metric. Therefore, we use the $18\text{-}22$th MLP layers as the layers where model edits are applied to update knowledge by modifying the parameters. We also visualize other mentioned LLMs in Figures \ref{fig:casual_gpt6b}, \ref{fig:casual_llama2}, \ref{fig:casual_yi9b}, \ref{fig:casual_internlm}, and \ref{fig:casual_mistral} to provide a more comprehensive results.

\begin{figure*}[t]
    \centering
    \begin{minipage}{\textwidth}
        \includegraphics[width=0.33\textwidth, height=0.2\textwidth]{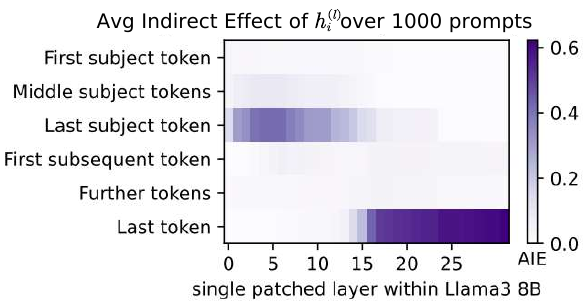}
        \includegraphics[width=0.33\textwidth, height=0.2\textwidth]{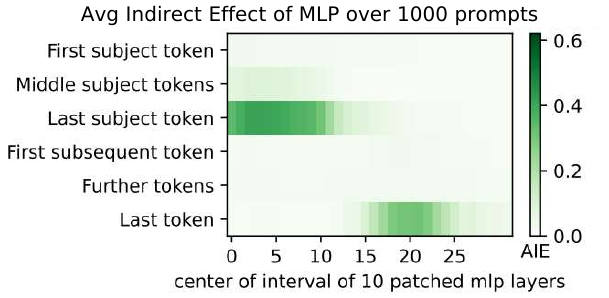}
        \includegraphics[width=0.33\textwidth, height=0.2\textwidth]{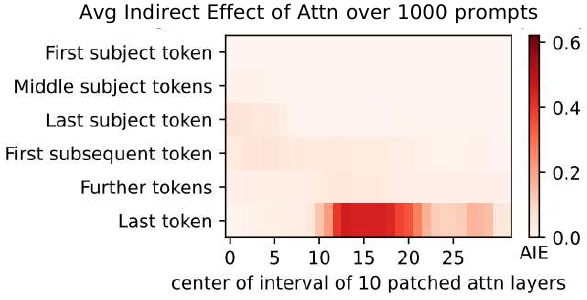}
        % \vspace{-1.5mm}
        \includegraphics[width=\textwidth, height=0.25\textwidth]{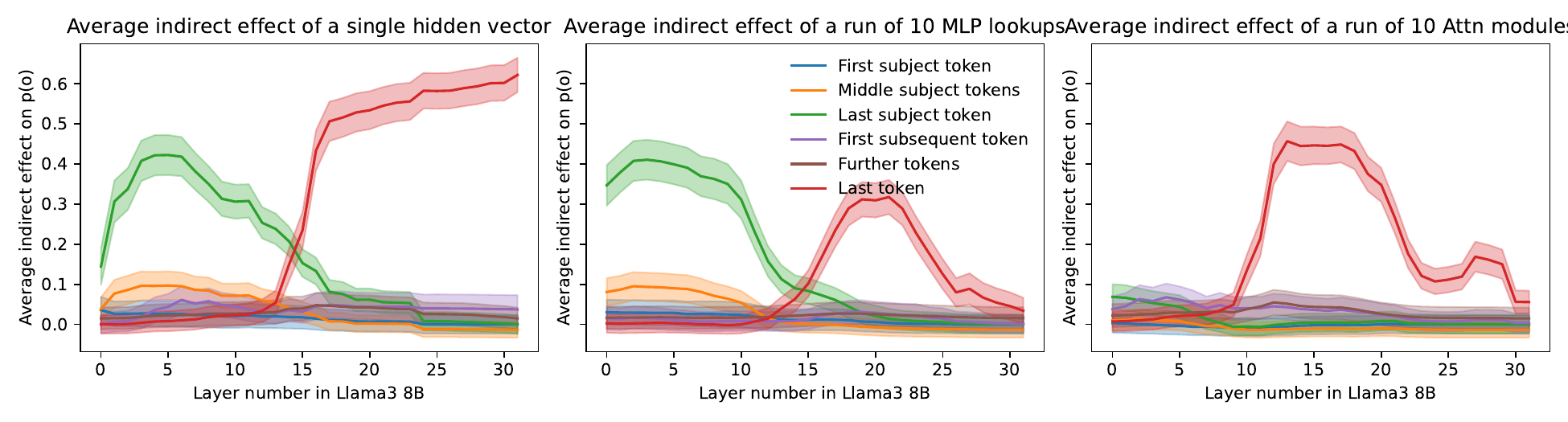}
        % \vspace{-1.5mm}
        \caption{Causal tracing visualization results for Llama3 8B. The causal impact on output probability is mapped for (a) the effect of each hidden state on the prediction, (b) the effect of MLP activations alone, and (c) the effect of attention activations alone. We also give according to mean causal traces of over a sample of 1000 factual statements, shown as a line plot with 95\% confidence intervals, which is below the first three figures. The confidence intervals confirm that the distinctions between peak and non-peak causal effects at both early and late sites are significant.}
        \label{fig:casual_llama3}
    \end{minipage}
\end{figure*}

\begin{figure*}[t]
    \centering
    \begin{minipage}{2\columnwidth}
        \includegraphics[width=0.33\textwidth, height=0.2\textwidth]{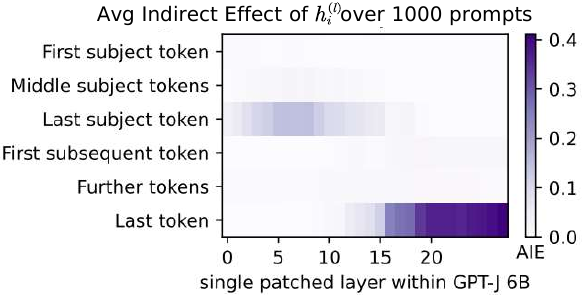}
        \includegraphics[width=0.33\textwidth, height=0.2\textwidth]{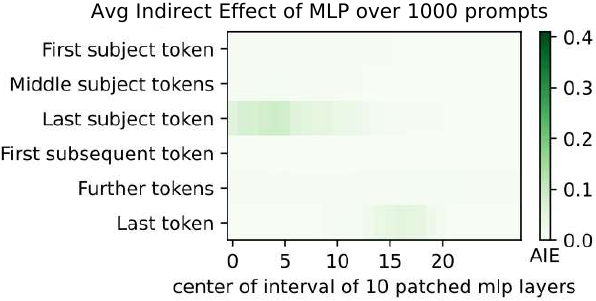}
        \includegraphics[width=0.33\textwidth, height=0.2\textwidth]{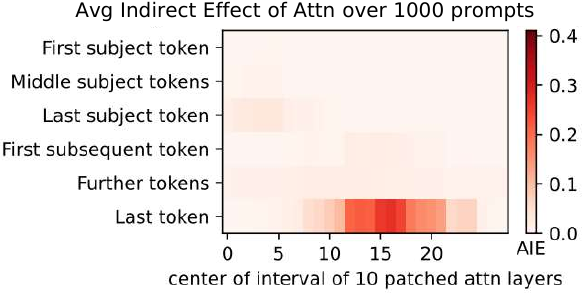}
        \includegraphics[width=\textwidth, height=0.25\textwidth]{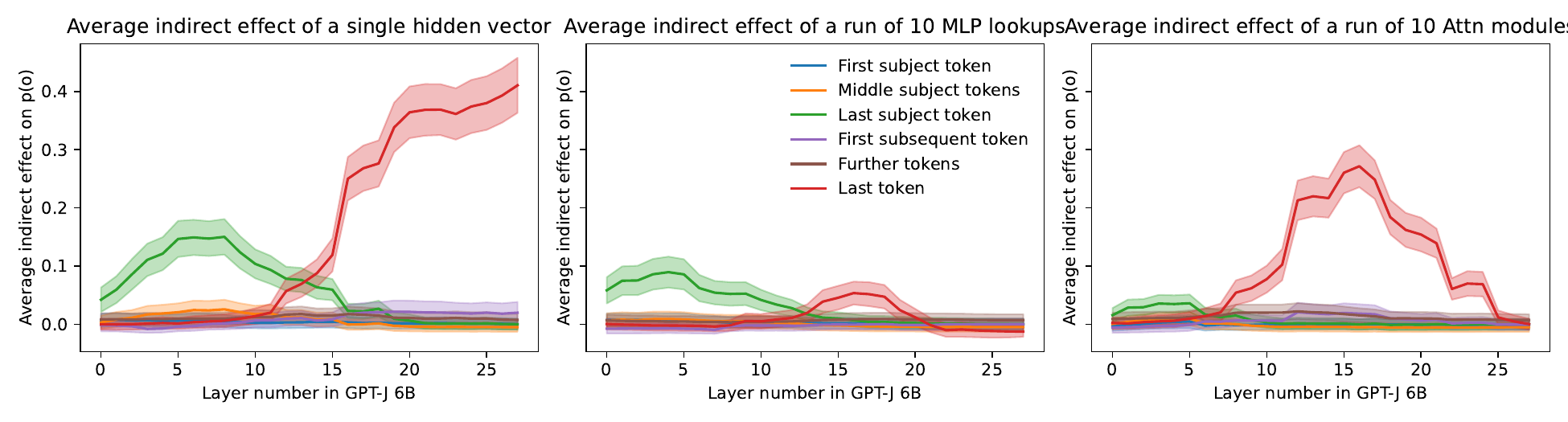}
        % \vspace{-1.5mm}
        \caption{Causal tracing visualization results for GPT-J 6B. The causal impact on output probability is mapped for (a) the effect of each hidden state on the prediction, (b) the effect of MLP activations alone, and (c) the effect of attention activations alone. We also give according to mean causal traces of over a sample of 1000 factual statements, shown as a line plot with 95\% confidence intervals, which is below the first three figures. The confidence intervals confirm that the distinctions between peak and non-peak causal effects at both early and late sites are significant.}
        \label{fig:casual_gpt6b}
    \end{minipage}
\end{figure*}

\begin{figure*}[t]
    \centering
    \begin{minipage}{2\columnwidth}
        \includegraphics[width=0.33\textwidth, height=0.2\textwidth]{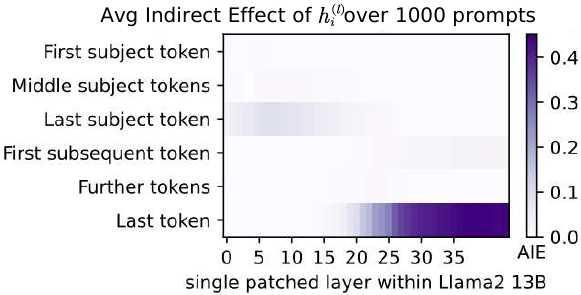}
        \includegraphics[width=0.33\textwidth, height=0.2\textwidth]{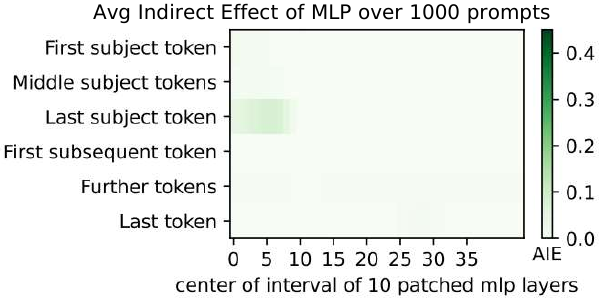}
        \includegraphics[width=0.33\textwidth, height=0.2\textwidth]{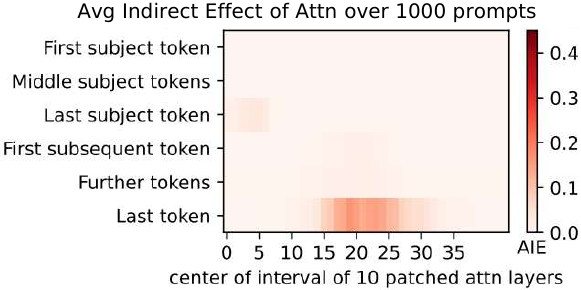}
        \includegraphics[width=\textwidth, height=0.25\textwidth]{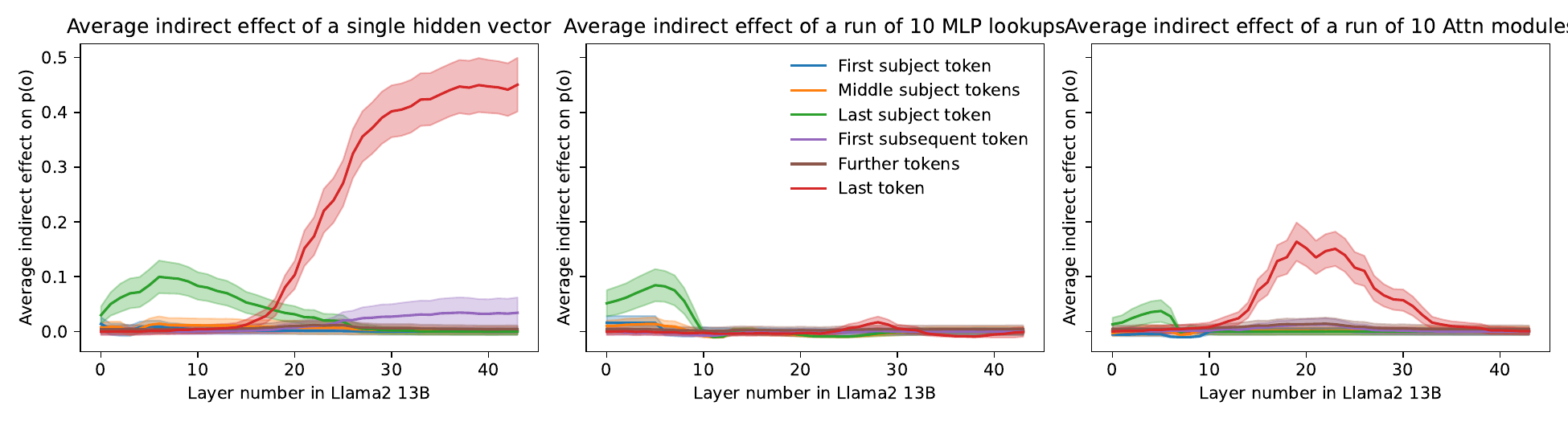}
        % \vspace{-1.5mm}
        \caption{Causal tracing visualization results for Llama2 13B. The causal impact on output probability is mapped for (a) the effect of each hidden state on the prediction, (b) the effect of MLP activations alone, and (c) the effect of attention activations alone. We also give according to mean causal traces of over a sample of 1000 factual statements, shown as a line plot with 95\% confidence intervals, which is below the first three figures. The confidence intervals confirm that the distinctions between peak and non-peak causal effects at both early and late sites are significant.}
        \label{fig:casual_llama2}
    \end{minipage}
\end{figure*}

\begin{figure*}[t]
    \centering
    \begin{minipage}{2\columnwidth}
        \includegraphics[width=0.33\textwidth, height=0.2\textwidth]{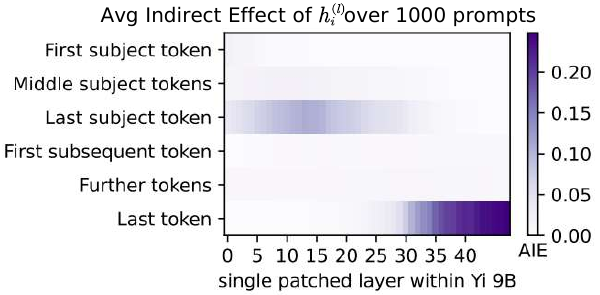}
        \includegraphics[width=0.33\textwidth, height=0.2\textwidth]{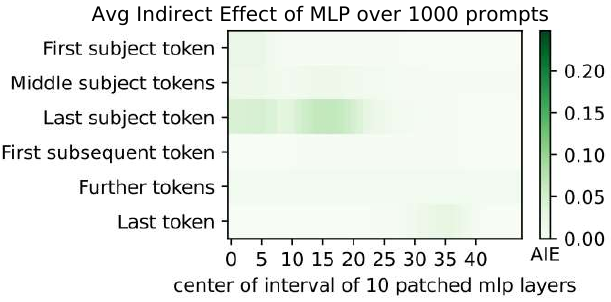}
        \includegraphics[width=0.33\textwidth, height=0.2\textwidth]{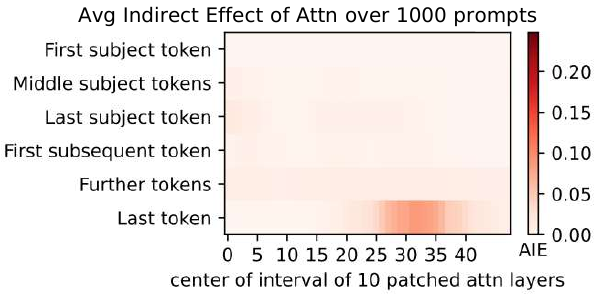}
        \includegraphics[width=\textwidth, height=0.25\textwidth]{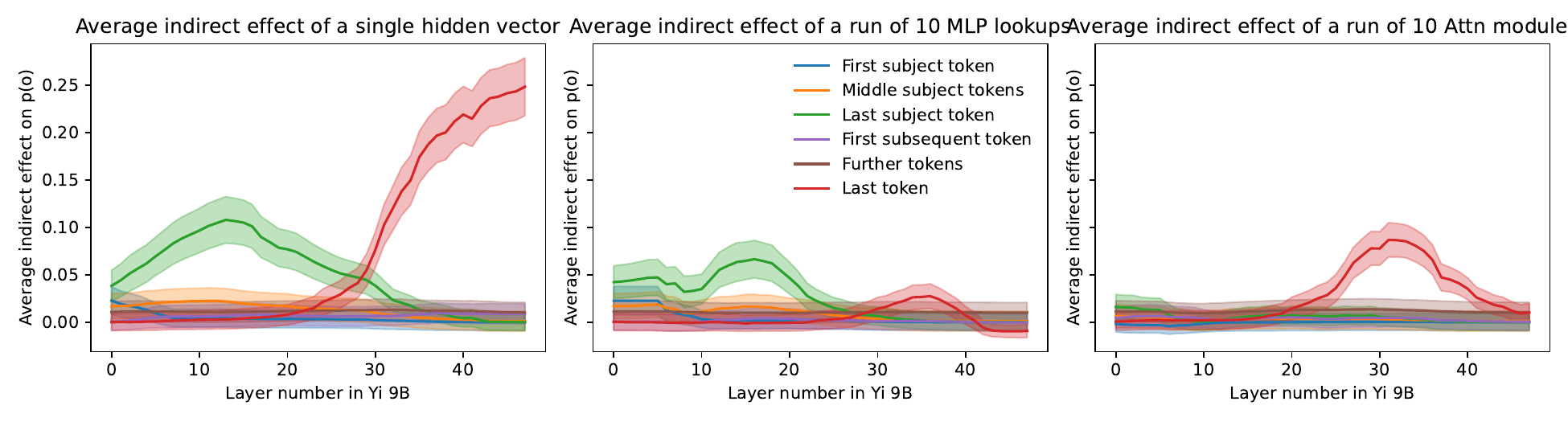}
        % \vspace{-1.5mm}
        \caption{Causal tracing visualization results for Yi 9B. The causal impact on output probability is mapped for (a) the effect of each hidden state on the prediction, (b) the effect of MLP activations alone, and (c) the effect of attention activations alone. We also give according to mean causal traces of over a sample of 1000 factual statements, shown as a line plot with 95\% confidence intervals, which is below the first three figures. The confidence intervals confirm that the distinctions between peak and non-peak causal effects at both early and late sites are significant.}
        \label{fig:casual_yi9b}
    \end{minipage}
\end{figure*}

\begin{figure*}[t]
    \centering
    \begin{minipage}{2\columnwidth}
        \includegraphics[width=0.33\textwidth, height=0.2\textwidth]{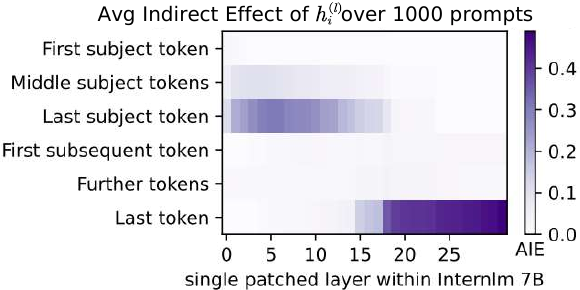}
        \includegraphics[width=0.33\textwidth, height=0.2\textwidth]{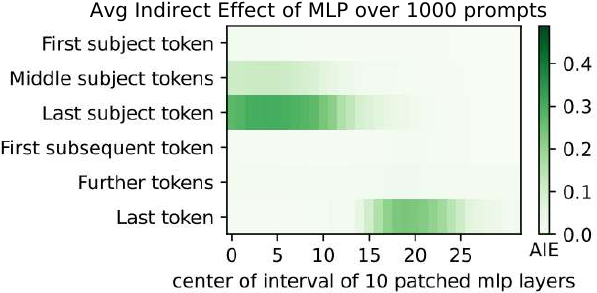}
        \includegraphics[width=0.33\textwidth, height=0.2\textwidth]{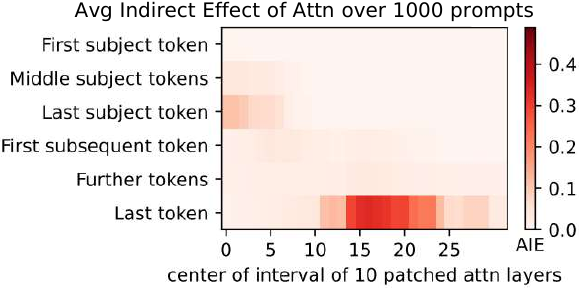}
        \includegraphics[width=\textwidth, height=0.25\textwidth]{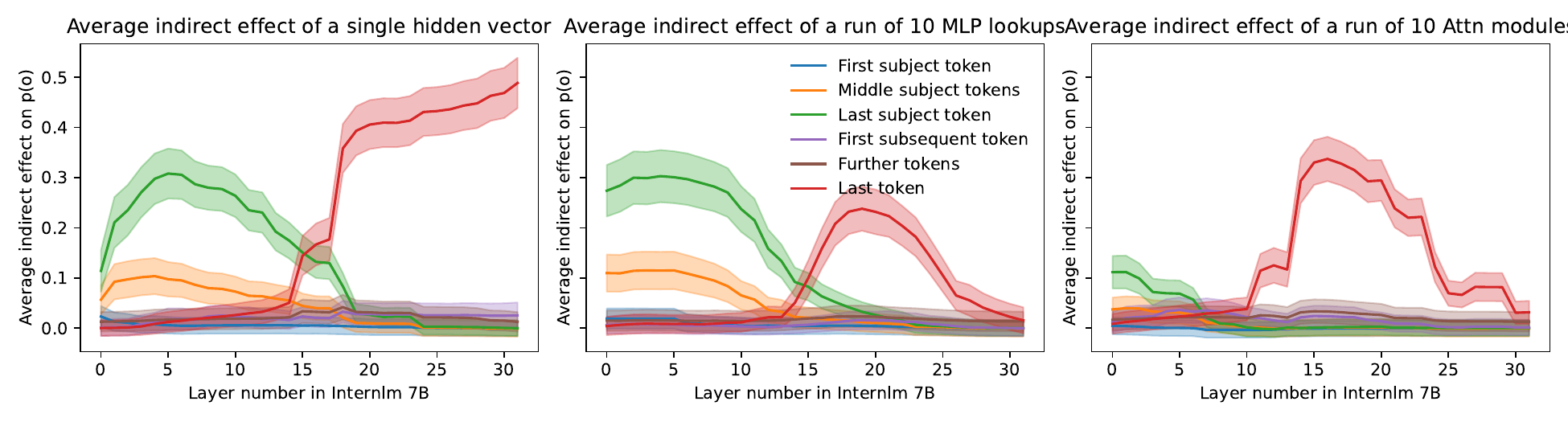}
        % \vspace{-1.5mm}
        \caption{Causal tracing visualization results for Internlm 7B. The causal impact on output probability is mapped for (a) the effect of each hidden state on the prediction, (b) the effect of MLP activations alone, and (c) the effect of attention activations alone. We also give according to mean causal traces of over a sample of 1000 factual statements, shown as a line plot with 95\% confidence intervals, which is below the first three figures. The confidence intervals confirm that the distinctions between peak and non-peak causal effects at both early and late sites are significant.}
        \label{fig:casual_internlm}
    \end{minipage}
\end{figure*}

\begin{figure*}[t]
    \centering
    \begin{minipage}{2\columnwidth}
        \includegraphics[width=0.33\textwidth, height=0.2\textwidth]{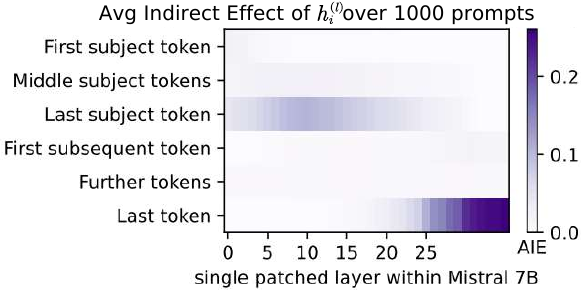}
        \includegraphics[width=0.33\textwidth, height=0.2\textwidth]{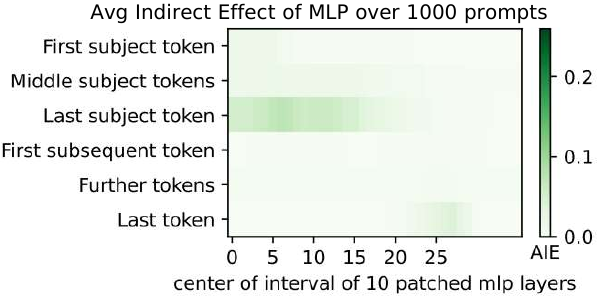}
        \includegraphics[width=0.33\textwidth, height=0.2\textwidth]{figs/mistral_2.pdf}
        % \vspace{-1.5mm}
        \includegraphics[width=\textwidth, height=0.25\textwidth]{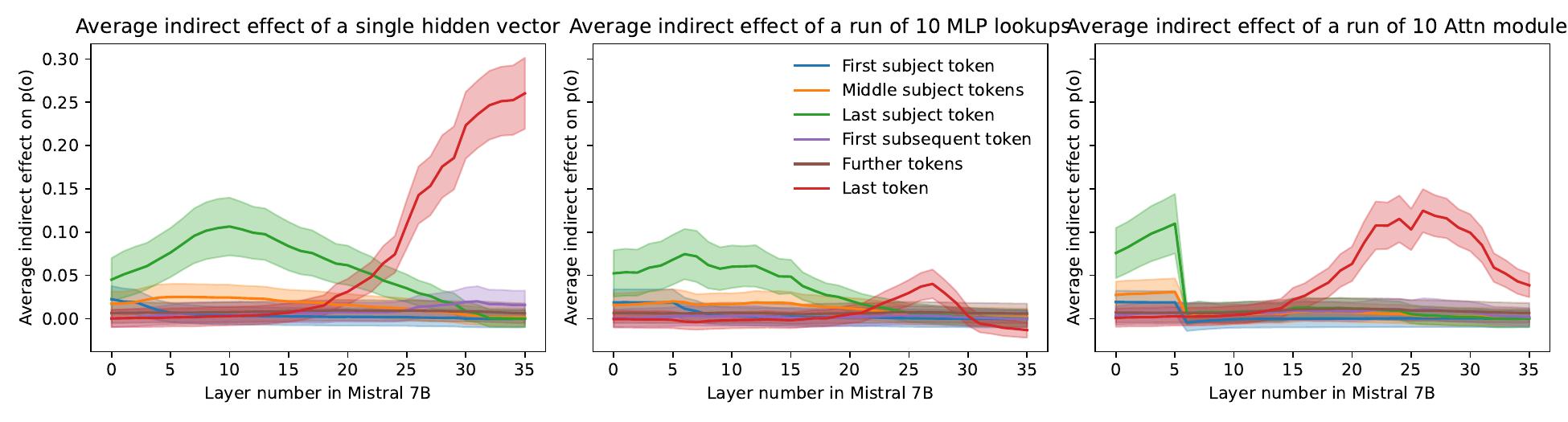}
        \caption{Causal tracing visualization results for mistral 7B. The causal impact on output probability is mapped for (a) the effect of each hidden state on the prediction, (b) the effect of MLP activations alone, and (c) the effect of attention activations alone. We also give according to mean causal traces of over a sample of 1000 factual statements, shown as a line plot with 95\% confidence intervals, which is below the first three figures. The confidence intervals confirm that the distinctions between peak and non-peak causal effects at both early and late sites are significant.}
        \label{fig:casual_mistral}
    \end{minipage}
\end{figure*}

\section{Prompt Templates for Each Task} \label{prompt_temp}
We list the prompt templates for different tasks to offer more visually intuitive results for each task in Table \ref{table:prompt}, respectively. More detailed prompt information for the best performance of each task and dataset can be seen within the code. 

\begin{table*}[h!]
    \centering
    \caption{Prompt templates for each task. In the prompts, "str(Examples)" represents initial solutions for each task. "str(Nodes)" and "str(Properties)" represent some basic properties of the materials/molecules. }\label{table:prompt}
    \begin{tabularx}{\textwidth}{p{0.25\textwidth}|p{0.7\textwidth}}
    \toprule
        \textbf{Optimization Tasks} & \textbf{Prompt Templates}  \\ \midrule
       Linear System Regression  & You will help me minimize a function with two input variables w, b. I have some (w, b) pairs and the function values at those points. The pairs are arranged in descending order based on their function values, where lower values are better.       
       
       Below are some examples: str(Examples)

       Give me a new (w, b) pair that is different from all pairs above, and has a function value lower than any of the examples. 
 \\ \midrule
        Travel Salesman Problem & You are given a list of points with coordinates below: str(Nodes)

        Below are some previous traces and their lengths. The traces are arranged in descending order based on their lengths, where lower values are better: str(Traces)
        
        Give me a new trace that is different from all traces above, and has a length lower than any of the above. The trace should traverse all points exactly once. The trace should start with <trace> and end with </trace>. 
 \\ \midrule
        Constitutive Law Prediction (Linear and Non-Linear) & You are given a list of materials with corresponding properties below: str(Properties)
        
        Below are some previous stress-strain responses for each material, sorted in descending order based on compliance, where lower compliance values indicate better structural integrity: str(Examples)
        
        Provide a new constitutive law that differs from all the laws above, ensuring that it produces a compliance value lower than any of the existing responses. The constitutive law should describe the relationship between stress and strain accurately for the material and should start with <law> and end with </law>. \\ \midrule
        Molecule Property Prediction (HOMO value, LUMO value, and HOMO-LUMO gap) & You are given a list of molecules with their chemical properties below: str(Properties). 

        Your goal is to predict the HOMO (can be adjusted to LUMO or HOMO-LUMO gap depending on different tasks) values for the listed molecules. Below are some examples: str(Examples)
        
        The HOMO value for each molecule should be distinct from any previously reported values, and the predictions should ensure an accurate representation of their electronic properties. Each predicted value should be formatted to start with <value> and end with </value>
        \\ \bottomrule
    \end{tabularx}
\end{table*}

\begin{table*}[h!]
    \centering
    \caption{Augmented prompt templates for linear system regression and travel salesman problem. In the prompts, "str(Examples)" represents initial solutions for each task. "str(Nodes)" represents some basic properties of the materials. }\label{table:aug_prompt1}
    \begin{tabularx}{\textwidth}{p{0.25\textwidth}|p{0.7\textwidth}}
    \toprule
        \textbf{Optimization Tasks} & \textbf{Augmented Prompts}  \\ \midrule
       Linear System Regression  & (i) You will help me minimize a function with two input variables w, b. I have some (w, b) pairs and the function values at those points. The pairs are arranged in descending order based on their function values, where lower values are better.

       Please provide a new (w, b) pair distinct from those above, with a lower function value than any previous pair.

       (ii) You are assisting in minimizing a function with two variables w and b. Provided are some (w, b) pairs along with their function values, sorted in descending order where lower values are better.

        Generate a (w, b) pair not seen in the above list, ensuring it yields a function value lower than any listed.
        
       (iii) Help me find the minimum of a function dependent on variables w and b. Below are (w, b) pairs and their corresponding function values, arranged from highest to lowest (lower is better).

       Give a (w, b) pair not seen in the above list, ensuring it yields a function value lower than any listed.
       
       (iv) I need assistance in minimizing a function with inputs w and b. Here are some (w, b) pairs and their function values, listed in descending order of their function values (lower values indicate better results).

       Give me a better (w, b) pair that is not included above.

       (v) Your task is to help minimize a function of two variables w and b. The following are (w, b) pairs and their function values, sorted from highest to lowest (lower values are preferable).

       Provide a new and better (w, b) pair from those above.
       
 \\ \midrule
        Travel Salesman Problem & (i) Given the coordinates of points: str(Nodes). Here are some previous routes and their lengths, sorted from longest to shortest (shorter is better). Give me a new trace that is different from all traces above, and has a length lower than any of the above. The trace should traverse all points exactly once. The trace should start with <trace> and end with </trace>.

        (ii) Consider the following points with coordinates: str(Nodes). Generate a new trace that differs from all previous traces and is shorter than any of them. Ensure that this trace covers each point exactly once, begins with <trace>, and ends with </trace>.

        (iii) You have a set of points at these coordinates below: str(Nodes). Please provide a unique trace that is distinct from all preceding traces and has a length shorter than any of the prior traces. This trace should start with <trace> and conclude with </trace>, while visiting each point one time.

        (iv) Here are points with their coordinates: str(Nodes). Please craft a new trace that is unique from all previous traces and shorter in length than any listed above. This trace should traverse all points a single time, beginning with <trace> and ending with </trace>.

        (v) You are provided with points with their coordinates: str(Nodes). Produce a trace that does not resemble any of the existing traces and has a length less than the shortest one listed. It should begin with <trace>, finish with </trace>, and cover each point exactly once.

        \\ \bottomrule
    \end{tabularx}
\end{table*}

\begin{table*}[h!]
    \centering
    \caption{Augmented prompt templates for constitutive law prediction and molecule property prediction. In the prompts, "str(Examples)" represents initial solutions for each task. "str(Properties)" represent some basic properties of the molecules. }\label{table:aug_prompt2}
    \begin{tabularx}{\textwidth}{p{0.25\textwidth}|p{0.7\textwidth}}
    \toprule
        \textbf{Optimization Tasks} & \textbf{Augmented Prompts}  \\ \midrule
        Constitutive Law Prediction (Linear and Non-Linear) & (i) You have been provided with a list of materials and their associated properties: str(Properties). You will find previous stress-strain responses for each material, where lower compliance indicates stronger structural integrity: str(Examples) Create a new constitutive law that is distinct from all those listed. The law should begin with <law> and end with </law>.
        
        (ii) Provided below are some materials with their properties: str(Properties). It also includes prior stress-strain responses for each material,  where lower values reflect better structural integrity: str(Examples). Develop a new constitutive law that is different from all previous laws. The law should start with <law> and end with </law>.

        (iii) You have a list of materials and their properties outlined below: str(Properties). Alongside, there are previous stress-strain responses for each material with lower compliance suggesting stronger structural integrity: str(Examples). Formulate a new constitutive law. This law should start with <law> and end with </law>.

        (iv) Below is a list of materials and their corresponding properties: str(Properties) The data is in descending order, where a lower compliance value indicates greater structural integrity: str(Examples). Construct a new constitutive law that is distinct from all the above, ensuring that it yields a compliance value below all existing responses. The law should start with <law> and end with </law>.

        (v) A list of materials and their properties is given below: str(Properties). Below are previous responses, where lower values suggest better structural integrity: str(Examples). Design a new constitutive law. The law should begin with <law> and end with </law>.
        \\ \midrule
        Molecule Property Prediction (HOMO value, LUMO value, and HOMO-LUMO gap) & (i) Provided below is a list of molecules and their properties: str(Properties). Predict the HOMO values for these molecules (or adjust to LUMO/HOMO-LUMO gap as needed). 
        
        Ensure each HOMO value accurately represents electronic properties. Format predictions with <value> at the start and </value> at the end.

        (ii) Below are molecules with their chemical properties: str(Properties). 
        
        Your task is to predict HOMO values (or LUMO/HOMO-LUMO gap) for each. Make sure predictions reflect electronic properties accurately, formatted as <value>...</value>.

        (iii) You have a list of molecules with properties: str(Properties). 
        
        Predict the HOMO (or LUMO/HOMO-LUMO gap) values. Use <value>...</value> to format predictions.

        (iv) Here are molecules and properties: str(Properties). 
        
        Predict HOMO values (or LUMO/HOMO-LUMO gap) that capture properties accurately. Prediction should begin with <value> and end with </value>.

        (v) Provided is a list of molecules: str(Properties). 
        
        Predict HOMO values (or LUMO/HOMO-LUMO gap), ensuring accurate electronic property representation. Format each as <value>...</value>.
        \\ \bottomrule
    \end{tabularx}
\end{table*}

\section{Details of Augmented Prompts} \label{augmented_prompt}

As mentioned in Section \ref{sec:case_study}, we discussed the robustness of the GSO method compared to traditional prompt-based approaches when handling semantically similar prompts. Specifically, we manually generated five prompts and used OpenAI o1 to rewrite an additional five prompts based on the original task prompts. We now list these augmented prompts in Tables \ref{table:aug_prompt1} and \ref{table:aug_prompt2} to provide more insights into our GSO. 

\begin{table*}[ht]
\centering
\caption{The inference time (per sample on average) on the TSP task for the baseline methods, including Vanilla, CoT, OPRO, Eureka, Funsearch, SGA, GSO, respectively. We report the results using Llama3 8B, GPT-J 6B, Llama2 13B, Yi9B, Internlm 7B, and Mistral 7B as backbone models, respectively.  (Unit: seconds)}

\label{tab:inference_time}
\resizebox{1.9\columnwidth}{!}{ 
\begin{tabular}{l c c c c c c c}
\toprule
Backbone Models & Vanilla & CoT & OPRO & Eureka & Funsearch & SGA &GSO (ours) \\
\midrule
Llama3 8B & 21.9  & 21.9 & 34.0 & 32.1 & 40.9 & 37.3 & 28.9 \\
GPT-J 6B & 20.1  & 21.3 & 30.8 & 27.5 & 44.5 & 35.8 & 29.0 \\
Llama2 13B & 40.2  & 41.4 & 53.3 & 47.6 & 61.0 & 55.3 & 51.0 \\
Yi 9B & 28.8  & 28.9 & 38.8 & 41.0 & 53.3 & 50.5 & 35.9 \\
Internlm 7B & 16.9  & 17.4 & 25.4 & 21.0 & 36.9 & 37.0 & 23.8 \\
Mistral 7B & 18.2  & 18.5 & 27.2 & 20.4 & 33.9 & 37.3 & 24.0 \\
\bottomrule
\end{tabular}
}

\end{table*}

\section{More Results of Different Edited Models} \label{app:backbone}
As mentioned in Section \ref{sec:main_res}, we select GPT-J 6B, Llama3 8B, and Mistral 7B as representative models in Table \ref{tab:main_res}. In this section, to further demonstrate the generalization and versatility of GSO, we also conducted experiments on several popular open-source LLMs at different scales, including Llama2 13B, Yi-9B, and Internlm 7B. As shown in Table \ref{tab:main_res_housange}, we can still observe that our GSO method significantly outperforms existing baselines. This further demonstrates the effectiveness of our GSO approach.
We also present radar charts for each model to provide a more intuitive performance comparison in Figures \ref{fig:leida_qiansan} and \ref{fig:leida_housan}. We apply a linear mapping to assign a score of $100$ to values with a zero loss. Higher loss values correspond to progressively lower scores, with an inability to answer the question resulting in a score of zero.
For formatting reasons, we also provide the radar chart of Llama3 8B. The effectiveness of our GSO across various popular open-source models further demonstrates its strong versatility and generalization capabilities.

\begin{figure*}[t]
    \centering
    \begin{minipage}{\textwidth}
        \includegraphics[width=0.33\textwidth, height=0.25\textwidth]{figs/leida_llama3.pdf}
        \includegraphics[width=0.33\textwidth, height=0.25\textwidth]{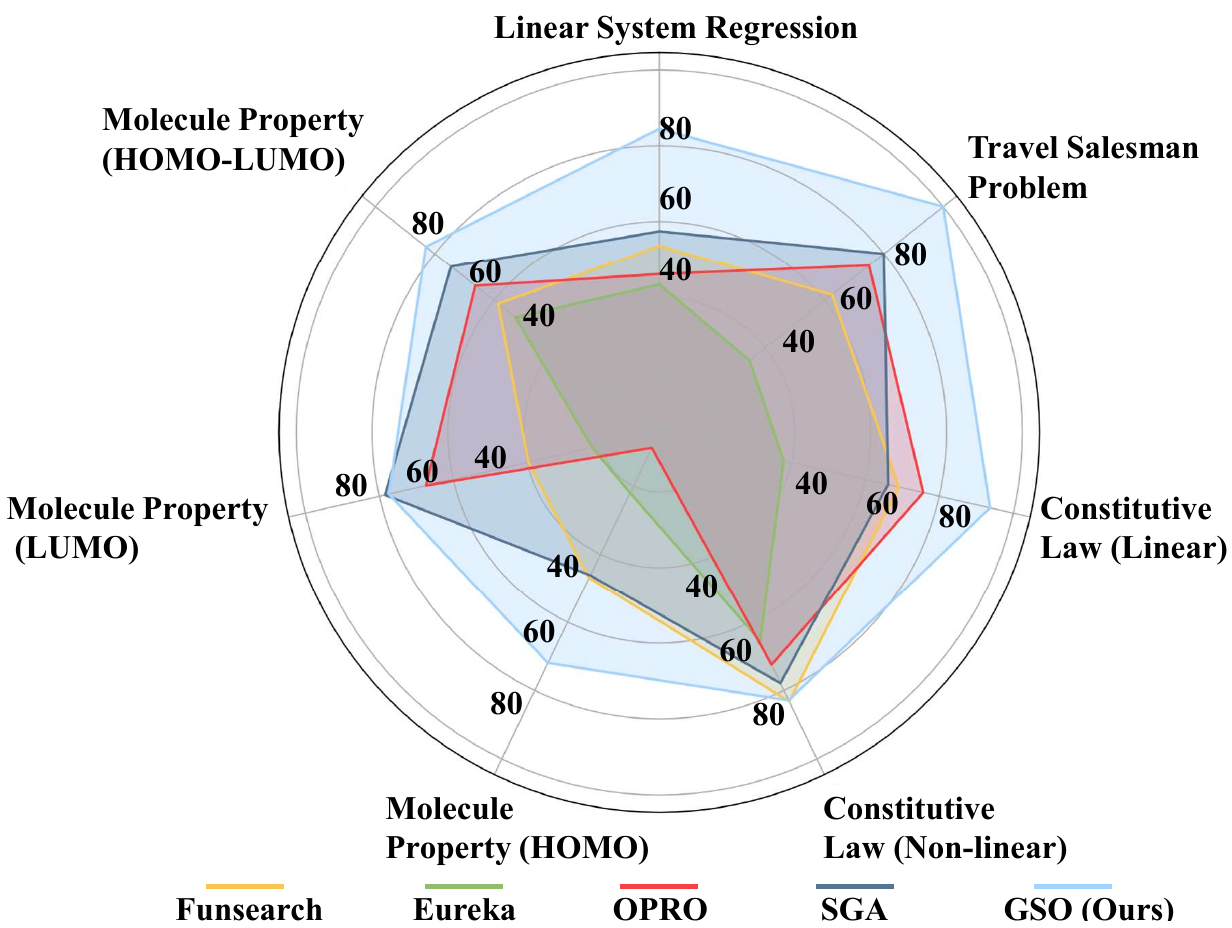}
        \includegraphics[width=0.33\textwidth, height=0.25\textwidth]{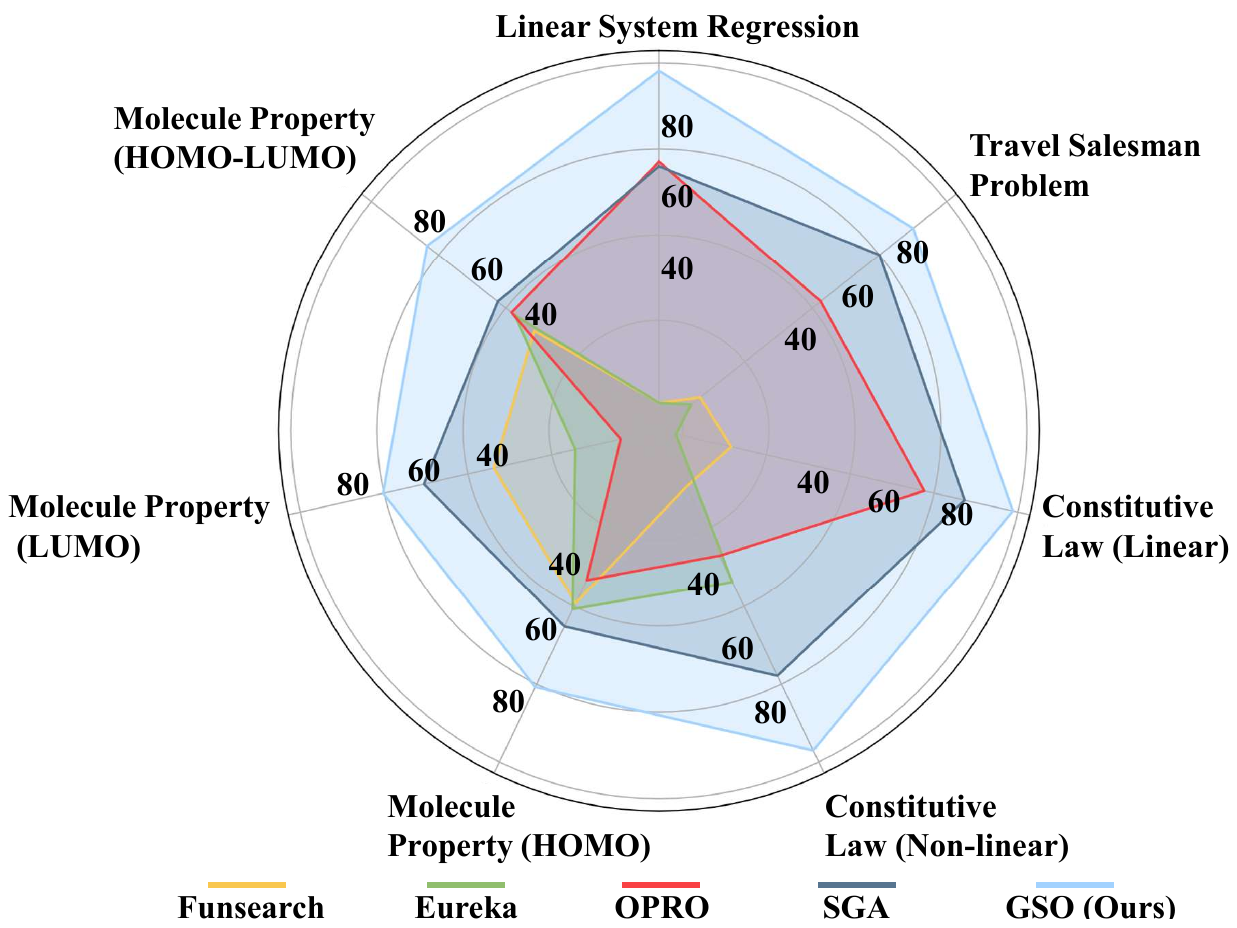}
        % \vspace{-1.5mm}
        \caption{GSO achieves state-of-the-art performance on a broad range of scientific optimization tasks compared with existing methods, using LLama 3 7B, GPT-J 6B, and Llama2 13B as backbone models, respectively.  We linearly map the evaluation metrics to [0, 100] for presentation.}
        \label{fig:leida_qiansan}
    \end{minipage}
\end{figure*}

\begin{figure*}[t]
    \centering
    \begin{minipage}{\textwidth}
        \includegraphics[width=0.33\textwidth, height=0.25\textwidth]{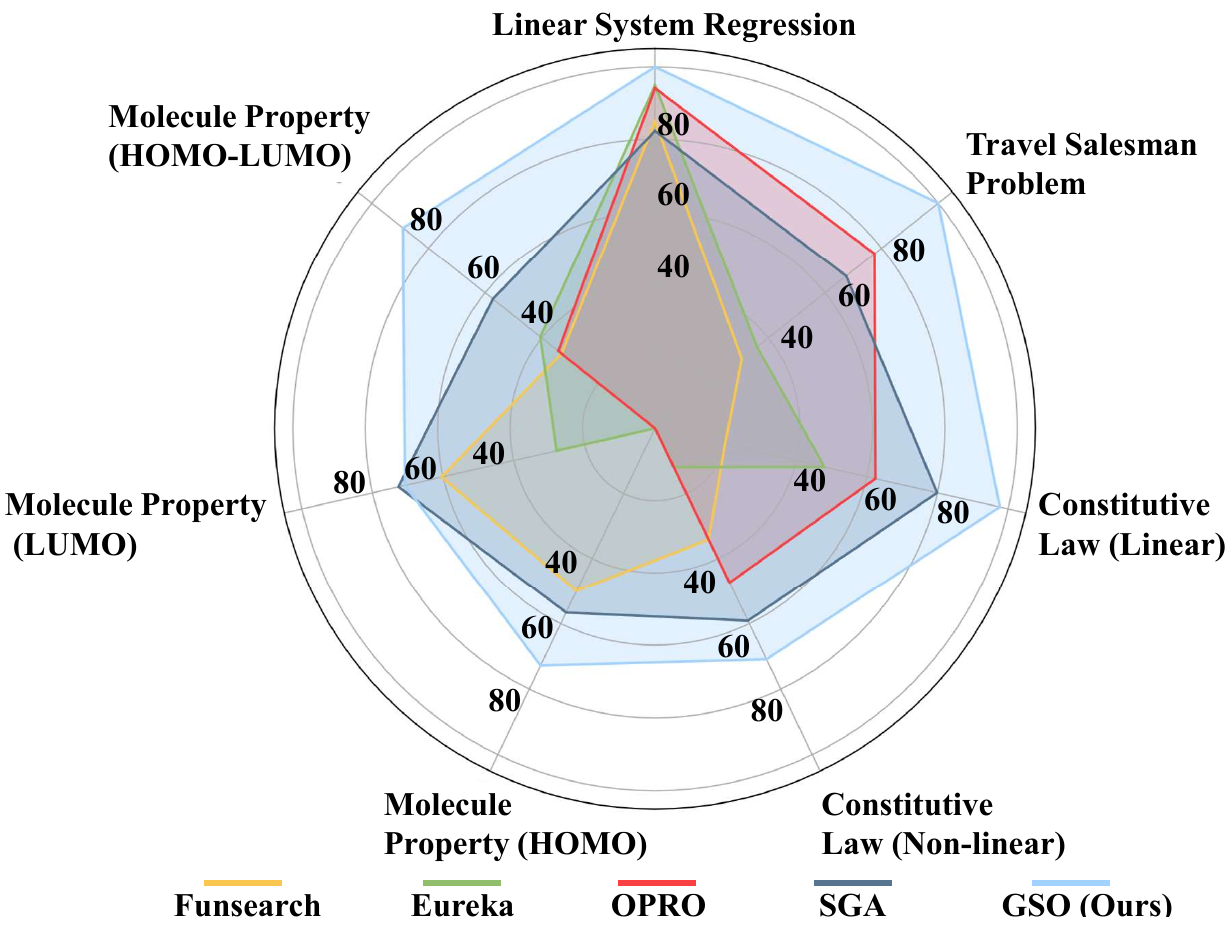}
        \includegraphics[width=0.33\textwidth, height=0.25\textwidth]{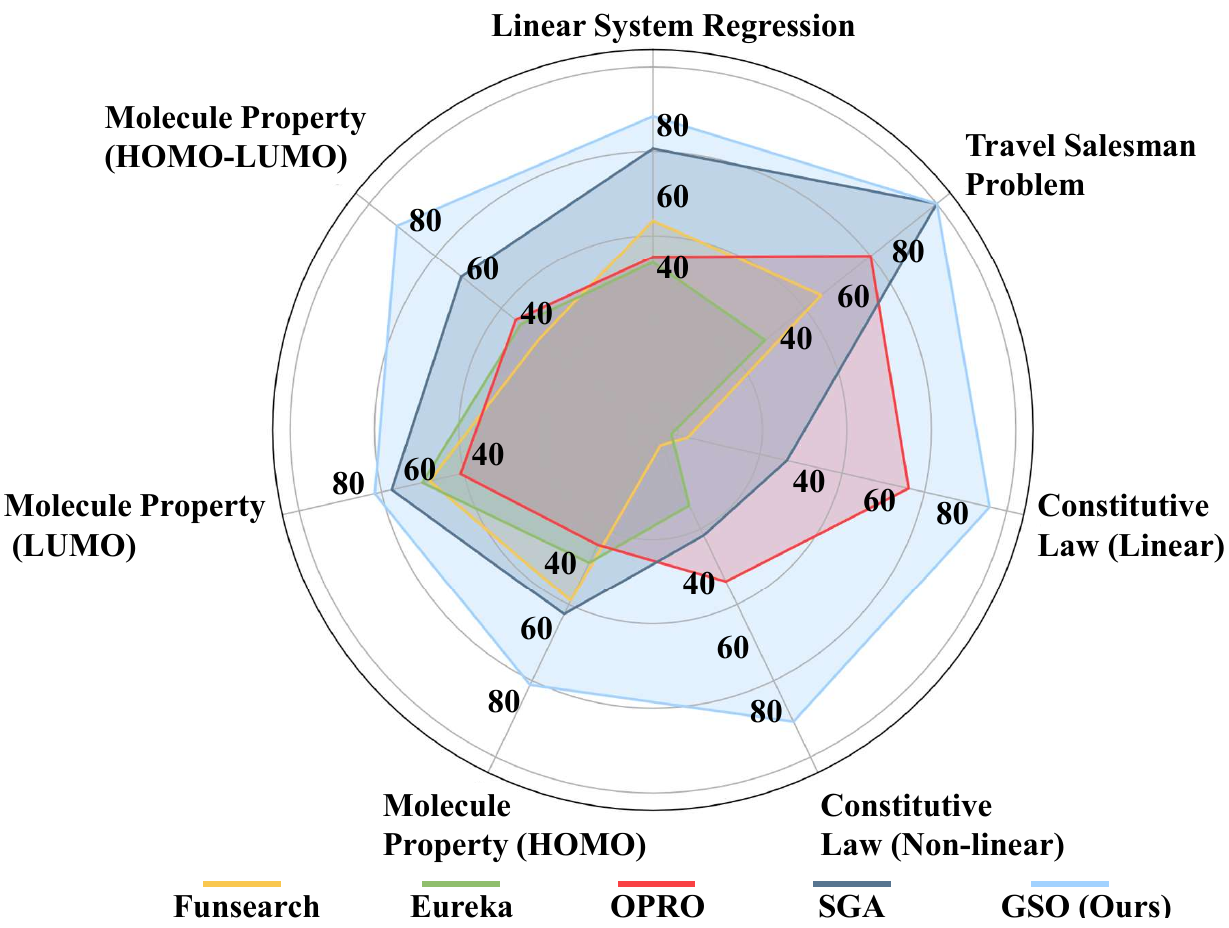}
        \includegraphics[width=0.33\textwidth, height=0.25\textwidth]{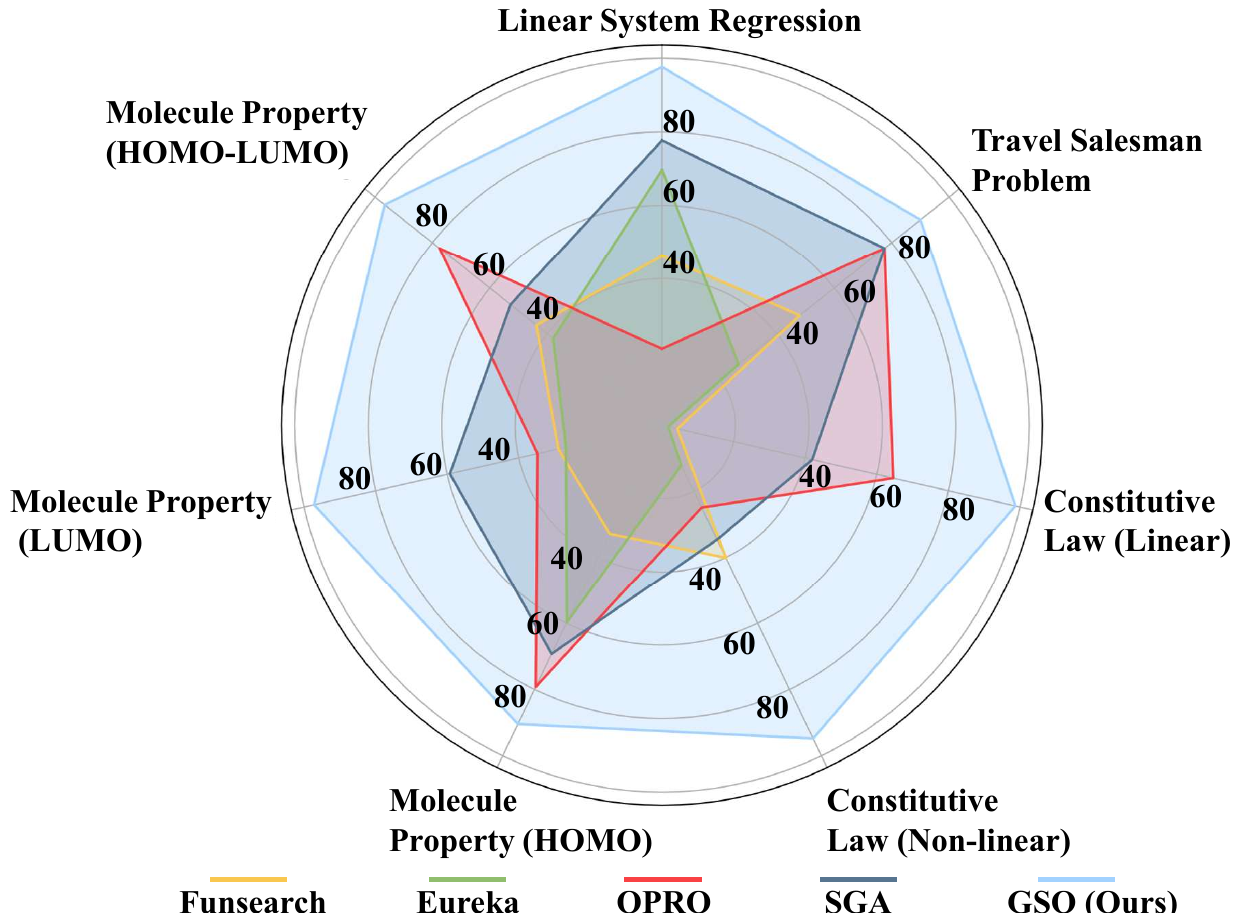}
        % \vspace{-1.5mm}
        \caption{GSO achieves state-of-the-art performance on a broad range of scientific optimization tasks compared with existing methods, using Yi 9B, Internlm 7B, and Mistral 7B as backbone models, respectively.  We linearly map the evaluation metrics to [0, 100] for presentation.}
        \label{fig:leida_housan}
    \end{minipage}
\end{figure*}

\section{More Details of Experiment Setups and Task Definitions} \label{app:def} 
We present more details of experiment setups and task definitions in this section.

\paragraph{Experiment Setups}
For the experimental setup, we apply several representative open-source LLMs: Llama3 8B\footnote{\url{https://huggingface.co/meta-llama/Meta-Llama-3-8B-Instruct}}, GPT-J-6B\footnote{\url{https://huggingface.co/EleutherAI/gpt-j-6b}}, Llama2 13B\footnote{\url{https://huggingface.co/meta-llama/Llama-2-13b-chat-hf}}, Yi 9B\footnote{\url{https://huggingface.co/01-ai/Yi-1.5-9B-Chat}}, InternLM 7B\footnote{\url{https://huggingface.co/internlm/internlm2_5-7b-chat}}, and Mistral 7B\footnote{\url{https://huggingface.co/mistralai/Mistral-7B-Instruct-v0.3}}. Other open-source models can also be replaced based on specific requirements. All experiments were conducted on a single Nvidia A100 GPU (80GB). Our approach is implemented using PyTorch 2.4.1\footnote{\url{https://pytorch.org/}} and Huggingface's Transformers 4.44.2\footnote{\url{https://github.com/huggingface/transformers}}. For each task, we set a maximum of $100$ iterations, with a limit of $2048$ tokens generated per iteration. The optimization process terminates if no performance improvement or accurate result is observed after $10$ consecutive iterations, or if an accurate result is reached earlier. More detailed configurations for the best performance of each task and dataset can be seen within our code.

For task definitions, we provide a more detailed of each task definition and setting. 

\paragraph{Linear System Regression} 
In linear system regression, the objective is to estimate the linear coefficients that most effectively model the relationship between input variables and their corresponding responses, within a probabilistic framework~\cite{fisher1922}. We follow the dataset setting as OPRO \cite{opro} which we include $437$ samples. Linear system regression is a fundamental technique in statistics and machine learning, widely applied in both theoretical and practical settings~\cite{montgomery2021introduction, seber2012linear}.

\paragraph{Traveling Salesman Problem (TSP)}  
The Traveling Salesman Problem (TSP)~\cite{junger1995traveling,gutin2006traveling} is a classical combinatorial optimization problem that has garnered significant attention in the literature, with numerous algorithms proposed, ranging from heuristic methods to exact solvers~\cite{rosenkrantz1977analysis,golden1980approximate,gurobi,helsgaun2017extension}. Recently, the problem has also been approached through the use of deep neural networks~\cite{kool2018attention,deudon2018learning,chen2019learning,nazari2018reinforcement}, underscoring its adaptability to modern machine learning techniques. Formally, given a set of \(n\) nodes with known coordinates, the TSP seeks to determine the shortest possible route that visits each node exactly once and returns to the starting point.We follow the dataset setting as OPRO \cite{opro}, which we include $177$ TSP problems. At each optimization step, a maximum of $8$ new solutions are generated. To evaluate the performance of various methods, we utilize the Gurobi solver~\cite{gurobi} to generate oracle solutions and compute the optimality gap. The optimality gap is defined as the relative difference between the distance of the solution obtained by the tested approach and that of the oracle solution, normalized by the oracle solution's distance.

\paragraph{Constitutive Law Prediction}  
Identifying constitutive laws from elastic material remains one of the most challenging tasks in fields such as physics, materials science, and mechanical engineering. In this paper, we follow the dataset settings as ~\cite{ma2023learning, bi_level}, including $357$ different linear and non-linear materials and  utilizing differentiable Material Point Method (MPM) simulators~\cite{sulsky1995application,jiang2016material}. 
The primary objective of this task is to uncover both the discrete structure and continuous parameters of a constitutive law, specifically identifying the material models $\varphi(\cdot)$ along with their associated material parameters $\theta$, from a ground-truth trajectory of particle positions $\hat{X}_{t \in [1, \dots, T]}$, where $T$ represents the number of time steps. In this context, we consider two classes of constitutive laws: $\varphi_E(\cdot; \theta_E)$ for modeling elastic materials and $\varphi_P(\cdot; \theta_P)$ for modeling plastic materials, both of which contain both linear and non-linear materials and are formally defined as the following:

\begin{subequations}
\begin{align*}
    \varphi_E\left(\mathbf{F};\theta_E\right)&\mapsto\boldsymbol{\tau}\\
    \varphi_P\left(\mathbf{F};\theta_P\right)&\mapsto\mathbf{F}^\text{corrected},
\end{align*}
\end{subequations}
where $\mathbf{F} \in \mathbb{R}^{3\times3}$ is the deformation gradient, ${\tau} \in \mathbb{R}^{3\times3}$ is the Kirchhoff stress tensor, $\mathbf{F}^\text{corrected} \in \mathbb{R}^{3\times3}$ is the deformation gradient after elastic return-mapping correction, and $\theta_{E}$ and $\theta_{P}$ are the continuous material parameters for elastic and elastic constitutive laws respectively. Given a specific constitutive law, we input it to the differentiable simulation and yields a particle position trajectory:
\begin{align*}
    X_{t\in\left[1,\dots,T\right]}=\text{sim}\left(\varphi\left(\cdot;\theta\right)\right),
\end{align*}
and we optimize the constitutive law by fitting the output trajectory to the ground truth $\hat{X}_{t\in\left[1,\dots, T\right]}$.

\paragraph{Molecule Property Prediction}  Predicting the Highest Occupied Molecular Orbital (HOMO), Lowest Unoccupied Molecular Orbital (LUMO), and HOMO-LUMO gap value of molecules is a key challenge in computational chemistry, particularly in the fields of quantum chemistry, material science, and drug design. In this task, we aim to predict the HOMO value based on molecular descriptors, such as molecular weight, Topological Polar Surface Area (TPSA), Octanol-water partition coefficient (logP), etc (maximum to 2k different properties). Although traditional quantum mechanical methods such as Density Functional Theory (DFT) are widely used for HOMO calculations, their computational cost can be prohibitive for large datasets. In this work, we adopt a machine learning approach to predict the HOMO value using these physical-chemical properties as input features. We follow the dataset settings as ~\cite{ma2023learning, bi_level}, including $238$ different molecules.

The primary goal of this task is to construct a model that accurately predicts the HOMO energy level \( \epsilon_\text{HOMO} \) by learning a mapping \( \psi(\cdot) \) between the input molecular properties \( x \in \mathbb{R}^n \), such as molecular weight, TPSA, and logP, and the HOMO energy \( \epsilon_\text{HOMO} \in \mathbb{R} \). The machine learning model is trained on a dataset containing known molecular properties and their corresponding HOMO values, and optimized to minimize the prediction error relative to the true HOMO values \( \hat{\epsilon}_\text{HOMO} \).

Given a set of molecular descriptors, the model outputs the predicted HOMO value:
\[
\epsilon_\text{HOMO} = \psi(x),
\]
and the objective is to minimize the loss function, defined as the difference between the predicted and true HOMO values:
\[
L = \left\| \psi(x) - \hat{\epsilon}_\text{HOMO} \right\|^{2}.
\]
This approach provides a computationally efficient alternative to traditional quantum mechanical methods, offering the potential for high-throughput screening of molecular libraries for their electronic properties, including HOMO values.

\begin{table*}[t]
\caption{The results of the ablation study of our GSO on the seven scientific optimization tasks, using GPT-J 6B as the backbone model. The symbol N/A indicates that the model is unable to provide a feasible solution for the current task.}\label{tab:abla_gpt6b}

\resizebox{\linewidth}{!}{
\begin{tabular}{l|ccccccc}
\toprule
\multirow{2}{*}{\textbf{Method}} & \textbf{Linear System} & \textbf{Travel Salesman} & \multicolumn{2}{c}{\textbf{Constitutive Law}} & \multicolumn{3}{c}{\textbf{Molecule Property}} \\ \cmidrule(lr){2-2} \cmidrule(lr){3-3} \cmidrule(lr){4-5} \cmidrule(lr){6-8}
 & \textbf{(a)} $\downarrow$ & \textbf{(b)}  $\downarrow$ & \textbf{(c)} $\downarrow$ & \textbf{(d)} $\downarrow$ & \textbf{(e)} $\downarrow$ & \textbf{(f)} $\downarrow$ & \textbf{(g)} $\downarrow$ \\
\midrule
\textbf{GSO$_{w/o \hspace{1mm} edit}$} & N/A & 0.5 $\pm$ 0.1 & 37.1 $\pm$ 10.0 & 185.3 $\pm$ 36.0 & 180.2 $\pm$ 19.1 & 497.1 $\pm$ 98.3  & 25.1 $\pm$ 4.1 \\
\textbf{GSO$_{w/o \hspace{1mm} dynamic}$} & 15.1 $\pm$ 3.2  & 0.0 $\pm$ 0.1 & 19.9 $\pm$ 6.1  & 131.7 $\pm$ 26.0 & 120.9 $\pm$ 64.4 & 135.9 $\pm$ 27.1 & 15.4 $\pm$ 4.1 \\
\midrule
\textbf{GSO (ours)} &
\textbf{12.3 $\pm$ 4.8} & \textbf{0.0 $\pm$ 0.0} & \textbf{15.7 $\pm$ 5.0} & {59.9 $\pm$ 10.3} & \textbf{70.1 $\pm$ 17.7} & \textbf{95.5 $\pm$ 10.0} & \textbf{8.5 $\pm$ 2.0} \\
\bottomrule
\end{tabular}
}
\end{table*}

\begin{table*}[t]
\caption{The results of the ablation study of our GSO on the seven scientific optimization tasks, using Llama2 13B as the backbone model.}\label{tab:abla_llama2}

\resizebox{\linewidth}{!}{
\begin{tabular}{l|ccccccc}
\toprule
\multirow{2}{*}{\textbf{Method}} & \textbf{Linear System} & \textbf{Travel Salesman} & \multicolumn{2}{c}{\textbf{Constitutive Law}} & \multicolumn{3}{c}{\textbf{Molecule Property}} \\ \cmidrule(lr){2-2} \cmidrule(lr){3-3} \cmidrule(lr){4-5} \cmidrule(lr){6-8}
 & \textbf{(a)} $\downarrow$ & \textbf{(b)}  $\downarrow$ & \textbf{(c)} $\downarrow$ & \textbf{(d)} $\downarrow$ & \textbf{(e)} $\downarrow$ & \textbf{(f)} $\downarrow$ & \textbf{(g)} $\downarrow$ \\
\midrule
\textbf{GSO$_{w/o \hspace{1mm} edit}$} & 22.3 $\pm$ 3.9  & 0.2 $\pm$ 0.0 & 73.3 $\pm$ 16.0 & 113.3 $\pm$ 20.1 & 540.2 $\pm$ 38.3 & 613.1 $\pm$ 49.9 & 332.1 $\pm$ 61.2 \\
\textbf{GSO$_{w/o \hspace{1mm} dynamic}$} & 13.7 $\pm$ 2.9  & \textbf{0.1 $\pm$ 0.1} & 6.3 $\pm$ 2.0 & 20.5 $\pm$ 4.0 & \textbf{35.3 $\pm$ 5.0} & 87.6 $\pm$ 14.1 & 19.3 $\pm$ 6.0 \\
\midrule
\textbf{GSO (ours)} &
\textbf{5.0 $\pm$ 1.6} & \textbf{0.1 $\pm$ 0.1} & \textbf{2.9 $\pm$ 1.3} & \textbf{7.4 $\pm$ 3.1} & {48.1 $\pm$ 12.0} & \textbf{79.1 $\pm$ 13.3} & \textbf{8.3 $\pm$ 3.2} \\

\bottomrule
\end{tabular}
}
\end{table*}

\begin{table*}[t]
\caption{The results of the ablation study of our GSO on the seven scientific optimization tasks, using Yi 9B as the backbone model.}\label{tab:abla_yi9b}

\resizebox{\linewidth}{!}{
\begin{tabular}{l|ccccccc}
\toprule
\multirow{2}{*}{\textbf{Method}} & \textbf{Linear System} & \textbf{Travel Salesman} & \multicolumn{2}{c}{\textbf{Constitutive Law}} & \multicolumn{3}{c}{\textbf{Molecule Property}} \\ \cmidrule(lr){2-2} \cmidrule(lr){3-3} \cmidrule(lr){4-5} \cmidrule(lr){6-8}
 & \textbf{(a)} $\downarrow$ & \textbf{(b)}  $\downarrow$ & \textbf{(c)} $\downarrow$ & \textbf{(d)} $\downarrow$ & \textbf{(e)} $\downarrow$ & \textbf{(f)} $\downarrow$ & \textbf{(g)} $\downarrow$ \\
\midrule
\textbf{GSO$_{w/o \hspace{1mm} edit}$} & 9.0 $\pm$ 1.5  & 0.4 $\pm$ 0.1 & 55.4 $\pm$ 10.1 & 408.1 $\pm$ 39.0 & 821.0 $\pm$ 112.0 & 982.1 $\pm$ 101.5 & 155.0 $\pm$ 30.1 \\
\textbf{GSO$_{w/o \hspace{1mm} dynamic}$} & 5.5 $\pm$ 1.7  &\textbf{ 0.0 $\pm$ 0.0 }& 19.1 $\pm$ 3.0 & \textbf{80.9 $\pm$ 5.0} & 73.0 $\pm$ 5.3 & 194.0 $\pm$ 20.1 & 9.0 $\pm$ 2.3 \\
\midrule
\textbf{GSO (ours)} &
\textbf{3.0 $\pm$ 0.8} & \textbf{0.0 $\pm$ 0.0} & \textbf{5.9 $\pm$ 2.1} & {89.1 $\pm$ 33.9} & \textbf{67.9 $\pm$ 23.1} & \textbf{{172.9 $\pm$ 43.1}} & \textbf{5.5 $\pm$ 2.1} \\
\bottomrule
\end{tabular}
}
\end{table*}

\begin{table*}[t]
\caption{The results of the ablation study of our GSO on the seven scientific optimization tasks, using Internlm 7B as the backbone model.}\label{tab:abla_internlm}

\resizebox{\linewidth}{!}{
\begin{tabular}{l|ccccccc}
\toprule
\multirow{2}{*}{\textbf{Method}} & \textbf{Linear System} & \textbf{Travel Salesman} & \multicolumn{2}{c}{\textbf{Constitutive Law}} & \multicolumn{3}{c}{\textbf{Molecule Property}} \\ \cmidrule(lr){2-2} \cmidrule(lr){3-3} \cmidrule(lr){4-5} \cmidrule(lr){6-8}
 & \textbf{(a)} $\downarrow$ & \textbf{(b)}  $\downarrow$ & \textbf{(c)} $\downarrow$ & \textbf{(d)} $\downarrow$ & \textbf{(e)} $\downarrow$ & \textbf{(f)} $\downarrow$ & \textbf{(g)} $\downarrow$ \\
\midrule
\textbf{GSO$_{w/o \hspace{1mm} edit}$} & 12.0 $\pm$ 3.1  & 0.3 $\pm$ 0.1 & 65.1 $\pm$ 10.5& 420 $\pm$ 82.1 & 490.0 $\pm$ 55.0 & 782.5 $\pm$ 133.2 & 150.2 $\pm$ 10.5 \\
\textbf{GSO$_{w/o \hspace{1mm} dynamic}$} &\textbf{ 9.0 $\pm$ 1.1}  & 0.1 $\pm$ 0.0 & 17.5 $\pm$ 4.3 & 79.0 $\pm$ 13.1 & 78.5 $\pm$ 20.0 & 159.3 $\pm$ 23.0 & 19.3 $\pm$ 2.5 \\
\midrule
\textbf{GSO (ours)} &
{10.0 $\pm$ 2.1} & \textbf{0.0 $\pm$ 0.0} & \textbf{10.4 $\pm$ 3.1} & \textbf{37.0 $\pm$ 12.9} & \textbf{47.2 $\pm$ 10.8} & \textbf{73.5 $\pm$ 19.6} & \textbf{4.2 $\pm$ 2.3} \\

\bottomrule
\end{tabular}
}
\end{table*}

\begin{table*}[t]
\caption{The results of the ablation study of our GSO on the seven scientific optimization tasks, using mistral 7B as the backbone model.}\label{tab:abla_mistral}

\resizebox{\linewidth}{!}{
\begin{tabular}{l|ccccccc}
\toprule
\multirow{2}{*}{\textbf{Method}} & \textbf{Linear System} & \textbf{Travel Salesman} & \multicolumn{2}{c}{\textbf{Constitutive Law}} & \multicolumn{3}{c}{\textbf{Molecule Property}} \\ \cmidrule(lr){2-2} \cmidrule(lr){3-3} \cmidrule(lr){4-5} \cmidrule(lr){6-8}
 & \textbf{(a)} $\downarrow$ & \textbf{(b)}  $\downarrow$ & \textbf{(c)} $\downarrow$ & \textbf{(d)} $\downarrow$ & \textbf{(e)} $\downarrow$ & \textbf{(f)} $\downarrow$ & \textbf{(g)} $\downarrow$ \\
\midrule
\textbf{GSO$_{w/o \hspace{1mm} edit}$} & 13.0 $\pm$ 2.5  & 0.3 $\pm$ 0.1 & 55.4 $\pm$ 13.5 & 31.7 $\pm$ 10.4 & 151.1 $\pm$ 20.0 & 45.4 $\pm$ 10.1 & 115.3 $\pm$ 20.1 \\
\textbf{GSO$_{w/o \hspace{1mm} dynamic}$} & 9.4 $\pm$ 2.0  & \textbf{0.0 $\pm$ 0.0} & 14.1 $\pm$ 3.0 & 16.0 $\pm$ 3.0 & 47.3 $\pm$ 6.8 & 10.5 $\pm$ 1.9 & 16.8 $\pm$ 4.0 \\
\midrule
\textbf{GSO (ours)} &
\textbf{5.6 $\pm$ 3.0} & {0.1 $\pm$ 0.0} & \textbf{3.2 $\pm$ 1.7} & \textbf{13.0 $\pm$ 2.1} & \textbf{23.9 $\pm$ 5.1} & \textbf{2.8 $\pm$ 1.5} & \textbf{1.7 $\pm$ 0.4} \\
\bottomrule
\end{tabular}
}
\end{table*}

\section{More Results of Ablation Study} \label{app:abla}
In Section \ref{sec:abla}, we report the results of the ablation study using Llama3 8B as the backbone model. In this section, we will further present the results using GPT-J-6B, Llama2 13B, Yi 9B, InternLM 7B, and Mistral 7B as backbone models to obtain more insights into the individual components constituting GSO across various backbone models. As illustrated in Tables \ref{tab:abla_gpt6b}, \ref{tab:abla_llama2}, \ref{tab:abla_yi9b}, \ref{tab:abla_internlm}, and \ref{tab:abla_mistral}, we still observe that the absence of each component within GSO leads to a decline in performance across diverse domains for almost all applied backbone models in the seven tested scientific optimization tasks, which further demonstrates that GSO organically integrates the two-level optimization into a unified framework as well. We also observe that the absence of model editing leads to a more significant decline in accuracy, which further demonstrates the importance of effectively utilizing observational feedback to adaptively adjust the optimization direction.

\section{More Case Study on Generalization or Memorization}

To investigate whether the improvements brought by our method are solely due to the LLM having seen solutions during the training phase, we designed an experiment aimed at eliminating this factor by having the model invent an imaginary constitutive law that does not exist on Earth. We combined the von Mises plasticity constitutive law, granular material, and weakly compressible fluid in proportions of $50\%$, $30\%$, and $20\%$, respectively, creating a new constitutive law that represents an exceedingly complex hypothetical material. As shown in Table \ref{tab:memory}, our method was still able to discover the constitutive law with minimal quantitative loss compared to other existing baselines. These results also indicate that GSO does not simply rely on memorization to achieve results but instead indeed effectively performs optimization. 
% From our observations, the visual differences between the ground truth material and the optimized constitutive law were minimal.

\section{Inference Time Comparisons}
We note that GSO requires an additional process to utilize the optimization feedback by applying model editing techniques to update LLM's parameters. Compared to the vanilla model, this process could potentially introduce extra inference time. However, we also observe that GSO effectively reduces the input prompt length through this approach, thereby decreasing the LLM’s inference time. Hence, to explore the interplay between these two effects, we record and compare the average inference time per sample of different methods using a single Nvidia A100 GPU (80GB), including vanilla, CoT, OPRO, Eureka, Funsearch, SGA, and GSO on each mentioned scientific optimization task to explore the influence of additional inference time and provide more insight of our GSO. 

We report the average inference time per sample as a metric, as shown in Table \ref{tab:inference_time}. We observe from the table that although the incorporation of model editing within the GSO framework introduces just marginal additional inference time costs. On average, the increase in inference time cost per sample due to the introduction of GSO, compared to the vanilla LLM, is $7.75$s. Notably, when utilizing Mistral 7B as the backbone model, this additional inference time is only $5.8$s, yet it offers a maximum precision improvement of 32.6$\times$ in the HOMO-LUMO gap prediction task compared to all other methods. Furthermore, GSO exhibits higher inference efficiency compared to existing methods. We may focus on exploring ways to further reduce the time cost of the GSO, including lightweighting LLMs to get a faster calculation \cite{qinglianghua, qinglianghua2} or adopting more efficient \cite{eff, eff3, eff2} and rational large model inference strategies such as speculative decoding \cite{sd, sd2} as for future works.

\section{Scientific Artifacts}

The data we collect in specialized domains is publicly available and viewable online. The data owners have indicated that the data can be used for scientific research or have not indicated that the data cannot be used for scientific research, and our collection process is also in compliance with regulations. Moreover, there is no unique identification of individuals or offensive content in these data.

\section{More Discussions On GSO} \label{app:more_discuss}

\subsection{What is the key advantage of using LLMs to optimize over some traditional optimization algorithms,
especially on classical optimization problems?}

The key advantage is that one can use natural language to describe the optimization problem. Instead of formally defining the optimization problem and deriving the update step with a programmed solver. This makes optimization more accessible to general users who may not have extensive domain knowledge of the specific types of optimization tasks in question, and it may also enhance the productivity of optimization experts who work on these tasks daily.
Moreover, employing LLMs for optimization is both domain-agnostic and task-agnostic, offering excellent universality and generalization. It enables rapid understanding and optimization of relatively unknown problems and can assist in finding more refined solutions. With the rapid advancement of LLMs, their abilities are improving swiftly, indicating the great potential of LLM-based optimization methods to address more complex and even unknown optimization tasks.

\subsection{What specific impact does the adoption of the exploitation/exploration phase actually have, given that sometimes satisfactory results can be achieved without introducing it? }
To investigate how the dynamic 
exploitation/exploration strategy works, we designed an ablation experiment in Section \ref{sec:abla} and Appendix \ref{app:abla}. Statistically, the strategy has a relatively positive effect, especially for tasks like predicting the HOMO-LUMO gap of molecules. We find that:
\textbf{ Without exploitation}, GSO sometimes finds it difficult to effectively provide consistent high-quality parameter hypotheses.
\textbf{Without exploration}, GSO sometimes becomes trapped in local optima, unable to further optimize the task to gain additional insights.
Consequently, using a balanced exploitation/exploration strategy enables one to escape local optima and obtain consistent loss decay. As shown in Figure \ref{fig:loss_curve}, GSO with this strategy is the only one that continues to progress toward better results, whereas others exhibit plateaued stagnation curves.

\subsection{Is GSO still effective when dealing with high-dimensional data or extremely intricate optimization tasks? }
Our experimental design includes both simple, low-dimensional ideal, and high-dimensional complex scientific optimization problems aimed at minimizing the sim-to-real gap to the greatest extent possible. This ensures that GSO generates meaningful and practical optimizations rather than questionable toys. We believe this is a strength rather than a weakness for GSO. Moreover, the success of GSO stems from the effective utilization of optimization feedback at each step and a dynamic exploration/exploitation strategy, which is domain-agnostic and can be applied to other domains. We also acknowledge that when the underlying physics of the optimization problem is too complex, this represents a limitation of our approach. A potential extension could focus on data compression by pre-training an auto-encoder to project high-dimensional data into a latent space.

\subsection{What is the difference between updating LLM parameters using model editing and employing methods like fine-tuning?}

The differences are two-folds:\textbf{(i)} 
the key difference lies in their\textbf{ scope and flexibility}. Fine-tuning typically involves adjusting a large portion of the model's parameters over a dataset, often requiring extensive computational resources and time \cite{lora,modify}, and it can result in overfitting or catastrophic forgetting of prior knowledge \cite{emp}. Model editing is a more targeted approach, allowing for precise modifications to specific parts of the model in response to new information or tasks without altering the broader knowledge encoded within the model. (ii)
Another difference is its \textbf{efficiency}; it enables rapid updates to the model without the need for extensive retraining. This makes it particularly suitable for dynamic optimization tasks where quick adjustments are necessary. Moreover, model editing preserves the generalization abilities of the LLM while fine-tuning may risk degrading its performance on unrelated tasks. In this sense, model editing offers a more controlled and adaptive method for refining LLM behavior, especially in domain-agnostic and task-agnostic scenarios, making it an ideal tool for iterative optimization processes.
% \section{Appendix}
% \label{sec:appendix}


\begin{thebibliography}{138}
\providecommand{\natexlab}[1]{#1}

\bibitem[{Aho and Ullman(1972)}]{Aho:72}
Alfred~V. Aho and Jeffrey~D. Ullman. 1972.
\newblock \emph{The Theory of Parsing, Translation and Compiling}, volume~1.
\newblock Prentice-Hall, Englewood Cliffs, NJ.

\bibitem[{AI4Science and Quantum(2023)}]{ai4sci}
Microsoft~Research AI4Science and Microsoft~Azure Quantum. 2023.
\newblock The impact of large language models on scientific discovery: a
  preliminary study using gpt-4.
\newblock \emph{arXiv preprint arXiv:2311.07361}.

\bibitem[{Amari(1993)}]{back_sgd}
Shun-ichi Amari. 1993.
\newblock Backpropagation and stochastic gradient descent method.
\newblock \emph{Neurocomputing}, 5(4-5):185--196.

\bibitem[{{American Psychological Association}(1983)}]{APA:83}
{American Psychological Association}. 1983.
\newblock \emph{Publications Manual}.
\newblock American Psychological Association, Washington, DC.

\bibitem[{Anderson(1972)}]{store2}
James~A Anderson. 1972.
\newblock A simple neural network generating an interactive memory.
\newblock \emph{Mathematical biosciences}, 14(3-4):197--220.

\bibitem[{Ando and Zhang(2005)}]{Ando2005}
Rie~Kubota Ando and Tong Zhang. 2005.
\newblock A framework for learning predictive structures from multiple tasks
  and unlabeled data.
\newblock \emph{Journal of Machine Learning Research}, 6:1817--1853.

\bibitem[{Andrew and Gao(2007)}]{andrew2007scalable}
Galen Andrew and Jianfeng Gao. 2007.
\newblock Scalable training of {L1}-regularized log-linear models.
\newblock In \emph{Proceedings of the 24th International Conference on Machine
  Learning}, pages 33--40.

\bibitem[{Bau et~al.(2020)Bau, Liu, Wang, Zhu, and Torralba}]{rewrite}
David Bau, Steven Liu, Tongzhou Wang, Jun-Yan Zhu, and Antonio Torralba. 2020.
\newblock Rewriting a deep generative model.
\newblock In \emph{Computer Vision--ECCV 2020: 16th European Conference,
  Glasgow, UK, August 23--28, 2020, Proceedings, Part I 16}, pages 351--369.
  Springer.

\bibitem[{Bengio and LeCun(2007)}]{Bengio+chapter2007}
Yoshua Bengio and Yann LeCun. 2007.
\newblock Scaling learning algorithms towards {AI}.
\newblock In \emph{Large Scale Kernel Machines}. MIT Press.

\bibitem[{Boiko et~al.(2023)Boiko, MacKnight, Kline, and Gomes}]{nature_chem}
Daniil~A Boiko, Robert MacKnight, Ben Kline, and Gabe Gomes. 2023.
\newblock Autonomous chemical research with large language models.
\newblock \emph{Nature}, 624(7992):570--578.

\bibitem[{Brown et~al.(2020)Brown, Mann, Ryder, Subbiah, Kaplan, Dhariwal,
  Neelakantan, Shyam, Sastry, Askell et~al.}]{few-shot}
Tom Brown, Benjamin Mann, Nick Ryder, Melanie Subbiah, Jared~D Kaplan, Prafulla
  Dhariwal, Arvind Neelakantan, Pranav Shyam, Girish Sastry, Amanda Askell,
  et~al. 2020.
\newblock Language models are few-shot learners.
\newblock \emph{Advances in neural information processing systems},
  33:1877--1901.

\bibitem[{Chandra et~al.(1981)Chandra, Kozen, and Stockmeyer}]{Chandra:81}
Ashok~K. Chandra, Dexter~C. Kozen, and Larry~J. Stockmeyer. 1981.
\newblock \href {https://doi.org/10.1145/322234.322243} {Alternation}.
\newblock \emph{Journal of the Association for Computing Machinery},
  28(1):114--133.

\bibitem[{Chen et~al.(2024{\natexlab{a}})Chen, Dohan, and So}]{EvoPrompting}
Angelica Chen, David Dohan, and David So. 2024{\natexlab{a}}.
\newblock Evoprompting: language models for code-level neural architecture
  search.
\newblock \emph{Advances in Neural Information Processing Systems}, 36.

\bibitem[{Chen et~al.(2023{\natexlab{a}})Chen, Scheurer, Korbak, Campos, Chan,
  Bowman, Cho, and Perez}]{code3}
Angelica Chen, J{\'e}r{\'e}my Scheurer, Tomasz Korbak, Jon~Ander Campos,
  Jun~Shern Chan, Samuel~R Bowman, Kyunghyun Cho, and Ethan Perez.
  2023{\natexlab{a}}.
\newblock Improving code generation by training with natural language feedback.
\newblock \emph{arXiv preprint arXiv:2303.16749}.

\bibitem[{Chen et~al.(2024{\natexlab{b}})Chen, Huang, Li, Chen, Lai, Xu, Gu,
  Gu, Yao, Xiao et~al.}]{hurt2}
Canyu Chen, Baixiang Huang, Zekun Li, Zhaorun Chen, Shiyang Lai, Xiongxiao Xu,
  Jia-Chen Gu, Jindong Gu, Huaxiu Yao, Chaowei Xiao, et~al. 2024{\natexlab{b}}.
\newblock Can editing llms inject harm?
\newblock \emph{arXiv preprint arXiv:2407.20224}.

\bibitem[{Chen et~al.(2024{\natexlab{c}})Chen, Shen, Lv, Wang, Ni, and
  Ye}]{sac-kg}
Hanzhu Chen, Xu~Shen, Qitan Lv, Jie Wang, Xiaoqi Ni, and Jieping Ye.
  2024{\natexlab{c}}.
\newblock Sac-kg: Exploiting large language models as skilled automatic
  constructors for domain knowledge graph.
\newblock In \emph{Proceedings of the 62nd Annual Meeting of the Association
  for Computational Linguistics (Volume 1: Long Papers)}, pages 4345--4360.

\bibitem[{Chen et~al.(2025)Chen, Shen, Wang, Wang, Lv, He, Wu, Ye, and
  Wu}]{chenknowledge}
Hanzhu Chen, Xu~Shen, Jie Wang, Zehao Wang, Qitan Lv, Junjie He, Rong Wu,
  Jieping Ye, and Feng Wu. 2025.
\newblock Knowledge graph finetuning enhances knowledge manipulation in large
  language models.
\newblock In \emph{The Thirteenth International Conference on Learning
  Representations}.

\bibitem[{Chen et~al.(2023{\natexlab{b}})Chen, Chen, Huang, and Zhou}]{spec2}
Jiuhai Chen, Lichang Chen, Heng Huang, and Tianyi Zhou. 2023{\natexlab{b}}.
\newblock When do you need chain-of-thought prompting for chatgpt?
\newblock \emph{arXiv preprint arXiv:2304.03262}.

\bibitem[{Chen et~al.()Chen, Lin, Sch{\"a}rli, and Zhou}]{code1}
Xinyun Chen, Maxwell Lin, Nathanael Sch{\"a}rli, and Denny Zhou.
\newblock Teaching large language models to self-debug.
\newblock In \emph{The Twelfth International Conference on Learning
  Representations}.

\bibitem[{Chen and Tian(2019)}]{chen2019learning}
Xinyun Chen and Yuandong Tian. 2019.
\newblock Learning to perform local rewriting for combinatorial optimization.
\newblock \emph{Advances in neural information processing systems}, 32.

\bibitem[{Chen et~al.(2022)Chen, Song, Lee, Wang, Zhang, Dohan, Kawakami,
  Kochanski, Doucet, Ranzato et~al.}]{OptFormer}
Yutian Chen, Xingyou Song, Chansoo Lee, Zi~Wang, Richard Zhang, David Dohan,
  Kazuya Kawakami, Greg Kochanski, Arnaud Doucet, Marc'aurelio Ranzato, et~al.
  2022.
\newblock Towards learning universal hyperparameter optimizers with
  transformers.
\newblock \emph{Advances in Neural Information Processing Systems},
  35:32053--32068.

\bibitem[{Chithrananda et~al.(2020)Chithrananda, Grand, and
  Ramsundar}]{chem_berta}
Seyone Chithrananda, Gabriel Grand, and Bharath Ramsundar. 2020.
\newblock Chemberta: large-scale self-supervised pretraining for molecular
  property prediction.
\newblock \emph{arXiv preprint arXiv:2010.09885}.

\bibitem[{Cova and Pais(2019)}]{op_chem}
T{\^a}nia~FGG Cova and Alberto~ACC Pais. 2019.
\newblock Deep learning for deep chemistry: optimizing the prediction of
  chemical patterns.
\newblock \emph{Frontiers in chemistry}, 7:809.

\bibitem[{Da et~al.(2021)Da, Bras, Lu, Choi, and Bosselut}]{ana1}
Jeff Da, Ronan~Le Bras, Ximing Lu, Yejin Choi, and Antoine Bosselut. 2021.
\newblock Analyzing commonsense emergence in few-shot knowledge models.
\newblock \emph{arXiv preprint arXiv:2101.00297}.

\bibitem[{De~Cao et~al.(2021)De~Cao, Aziz, and Titov}]{ke}
Nicola De~Cao, Wilker Aziz, and Ivan Titov. 2021.
\newblock Editing factual knowledge in language models.
\newblock \emph{arXiv preprint arXiv:2104.08164}.

\bibitem[{Deudon et~al.(2018)Deudon, Cournut, Lacoste, Adulyasak, and
  Rousseau}]{deudon2018learning}
Michel Deudon, Pierre Cournut, Alexandre Lacoste, Yossiri Adulyasak, and
  Louis-Martin Rousseau. 2018.
\newblock Learning heuristics for the tsp by policy gradient.
\newblock In \emph{Integration of Constraint Programming, Artificial
  Intelligence, and Operations Research: 15th International Conference, CPAIOR
  2018, Delft, The Netherlands, June 26--29, 2018, Proceedings 15}, pages
  170--181. Springer.

\bibitem[{Dhingra et~al.()Dhingra, Zaheer, Balachandran, Neubig, Salakhutdinov,
  and Cohen}]{qa1}
Bhuwan Dhingra, Manzil Zaheer, Vidhisha Balachandran, Graham Neubig, Ruslan
  Salakhutdinov, and William~W Cohen.
\newblock Differentiable reasoning over a virtual knowledge base.
\newblock In \emph{International Conference on Learning Representations}.

\bibitem[{Didelez and Pigeot(2001)}]{cg_01}
Vanessa Didelez and Iris Pigeot. 2001.
\newblock Causality: models, reasoning, and inference.

\bibitem[{Fang et~al.(2022)Fang, Liu, Lei, He, Zhang, Zhou, Wang, Wu, and
  Wang}]{data_set}
Xiaomin Fang, Lihang Liu, Jieqiong Lei, Donglong He, Shanzhuo Zhang, Jingbo
  Zhou, Fan Wang, Hua Wu, and Haifeng Wang. 2022.
\newblock Geometry-enhanced molecular representation learning for property
  prediction.
\newblock \emph{Nature Machine Intelligence}, 4(2):127--134.

\bibitem[{Fisher(1922)}]{fisher1922}
R.~A. Fisher. 1922.
\newblock On the mathematical foundations of theoretical statistics.
\newblock \emph{Philosophical Transactions of the Royal Society of London.
  Series A}, 222:309--368.

\bibitem[{Fortunato et~al.(2018)Fortunato, Bergstrom, B{\"o}rner, Evans,
  Helbing, Milojevi{\'c}, Petersen, Radicchi, Sinatra, Uzzi et~al.}]{sos}
Santo Fortunato, Carl~T Bergstrom, Katy B{\"o}rner, James~A Evans, Dirk
  Helbing, Sta{\v{s}}a Milojevi{\'c}, Alexander~M Petersen, Filippo Radicchi,
  Roberta Sinatra, Brian Uzzi, et~al. 2018.
\newblock Science of science.
\newblock \emph{Science}, 359(6379):eaao0185.

\bibitem[{Frantar and Alistarh(2023)}]{eff2}
Elias Frantar and Dan Alistarh. 2023.
\newblock Qmoe: Practical sub-1-bit compression of trillion-parameter models.
\newblock \emph{arXiv preprint arXiv:2310.16795}.

\bibitem[{Geva et~al.(2022)Geva, Caciularu, Wang, and Goldberg}]{ana3}
Mor Geva, Avi Caciularu, Kevin~Ro Wang, and Yoav Goldberg. 2022.
\newblock Transformer feed-forward layers build predictions by promoting
  concepts in the vocabulary space.
\newblock \emph{arXiv preprint arXiv:2203.14680}.

\bibitem[{Geva et~al.(2020)Geva, Schuster, Berant, and Levy}]{ana2}
Mor Geva, Roei Schuster, Jonathan Berant, and Omer Levy. 2020.
\newblock Transformer feed-forward layers are key-value memories.
\newblock \emph{arXiv preprint arXiv:2012.14913}.

\bibitem[{Golden et~al.(1980)Golden, Bodin, Doyle, and
  Stewart~Jr}]{golden1980approximate}
Bruce Golden, Lawrence Bodin, T~Doyle, and W~Stewart~Jr. 1980.
\newblock Approximate traveling salesman algorithms.
\newblock \emph{Operations research}, 28(3-part-ii):694--711.

\bibitem[{Goodfellow et~al.(2016)Goodfellow, Bengio, Courville, and
  Bengio}]{goodfellow2016deep}
Ian Goodfellow, Yoshua Bengio, Aaron Courville, and Yoshua Bengio. 2016.
\newblock \emph{Deep learning}, volume~1.
\newblock MIT Press.

\bibitem[{Gu et~al.(2024)Gu, Xu, Ma, Lu, Ling, Chang, and Peng}]{hurt}
Jia-Chen Gu, Hao-Xiang Xu, Jun-Yu Ma, Pan Lu, Zhen-Hua Ling, Kai-Wei Chang, and
  Nanyun Peng. 2024.
\newblock Model editing can hurt general abilities of large language models.
\newblock \emph{arXiv preprint arXiv:2401.04700}.

\bibitem[{Gurnee et~al.()Gurnee, Nanda, Pauly, Harvey, Troitskii, and
  Bertsimas}]{finding}
Wes Gurnee, Neel Nanda, Matthew Pauly, Katherine Harvey, Dmitrii Troitskii, and
  Dimitris Bertsimas.
\newblock Finding neurons in a haystack: Case studies with sparse probing.
\newblock \emph{Transactions on Machine Learning Research}.

\bibitem[{{Gurobi Optimization, LLC}(2024)}]{gurobi}
{Gurobi Optimization, LLC}. 2024.
\newblock \href {https://www.gurobi.com} {{Gurobi Optimizer Reference Manual}}.

\bibitem[{Gusfield(1997)}]{Gusfield:97}
Dan Gusfield. 1997.
\newblock \emph{Algorithms on Strings, Trees and Sequences}.
\newblock Cambridge University Press, Cambridge, UK.

\bibitem[{Gutin and Punnen(2006)}]{gutin2006traveling}
Gregory Gutin and Abraham~P Punnen. 2006.
\newblock \emph{The traveling salesman problem and its variations}, volume~12.
\newblock Springer Science \& Business Media.

\bibitem[{Guu et~al.(2020)Guu, Lee, Tung, Pasupat, and Chang}]{qa2}
Kelvin Guu, Kenton Lee, Zora Tung, Panupong Pasupat, and Mingwei Chang. 2020.
\newblock Retrieval augmented language model pre-training.
\newblock In \emph{International conference on machine learning}, pages
  3929--3938. PMLR.

\bibitem[{Hanna et~al.(2024)Hanna, Liu, and Variengien}]{greater}
Michael Hanna, Ollie Liu, and Alexandre Variengien. 2024.
\newblock How does gpt-2 compute greater-than?: Interpreting mathematical
  abilities in a pre-trained language model.
\newblock \emph{Advances in Neural Information Processing Systems}, 36.

\bibitem[{Hao et~al.(2023)Hao, Tan, Tang, Ni, Shao, Zhang, Xing, and
  Hu}]{bertnet}
Shibo Hao, Bowen Tan, Kaiwen Tang, Bin Ni, Xiyan Shao, Hengzhe Zhang, Eric
  Xing, and Zhiting Hu. 2023.
\newblock Bertnet: Harvesting knowledge graphs with arbitrary relations from
  pretrained language models.
\newblock In \emph{Findings of the Association for Computational Linguistics:
  ACL 2023}, pages 5000--5015.

\bibitem[{Hartvigsen et~al.(2024)Hartvigsen, Sankaranarayanan, Palangi, Kim,
  and Ghassemi}]{pre_edit1}
Tom Hartvigsen, Swami Sankaranarayanan, Hamid Palangi, Yoon Kim, and Marzyeh
  Ghassemi. 2024.
\newblock Aging with grace: Lifelong model editing with discrete key-value
  adaptors.
\newblock \emph{Advances in Neural Information Processing Systems}, 36.

\bibitem[{Helsgaun(2017)}]{helsgaun2017extension}
Keld Helsgaun. 2017.
\newblock An extension of the lin-kernighan-helsgaun tsp solver for constrained
  traveling salesman and vehicle routing problems.
\newblock \emph{Roskilde: Roskilde University}, 12:966--980.

\bibitem[{Hinton et~al.(2006)Hinton, Osindero, and Teh}]{Hinton06}
Geoffrey~E. Hinton, Simon Osindero, and Yee~Whye Teh. 2006.
\newblock A fast learning algorithm for deep belief nets.
\newblock \emph{Neural Computation}, 18:1527--1554.

\bibitem[{Hsieh et~al.(2023)Hsieh, Li, Yeh, Nakhost, Fujii, Ratner, Krishna,
  Lee, and Pfister}]{qinglianghua2}
Cheng-Yu Hsieh, Chun-Liang Li, Chih-Kuan Yeh, Hootan Nakhost, Yasuhisa Fujii,
  Alexander Ratner, Ranjay Krishna, Chen-Yu Lee, and Tomas Pfister. 2023.
\newblock Distilling step-by-step! outperforming larger language models with
  less training data and smaller model sizes.
\newblock \emph{arXiv preprint arXiv:2305.02301}.

\bibitem[{Hu et~al.(2021)Hu, Shen, Wallis, Allen-Zhu, Li, Wang, Wang, and
  Chen}]{lora}
Edward~J Hu, Yelong Shen, Phillip Wallis, Zeyuan Allen-Zhu, Yuanzhi Li, Shean
  Wang, Lu~Wang, and Weizhu Chen. 2021.
\newblock Lora: Low-rank adaptation of large language models.
\newblock \emph{arXiv preprint arXiv:2106.09685}.

\bibitem[{Intriligator(2002)}]{math_opt}
Michael~D Intriligator. 2002.
\newblock \emph{Mathematical optimization and economic theory}.
\newblock SIAM.

\bibitem[{Jiang et~al.(2023)Jiang, Sablayrolles, Mensch, Bamford, Chaplot,
  Casas, Bressand, Lengyel, Lample, Saulnier et~al.}]{mistral}
Albert~Q Jiang, Alexandre Sablayrolles, Arthur Mensch, Chris Bamford,
  Devendra~Singh Chaplot, Diego de~las Casas, Florian Bressand, Gianna Lengyel,
  Guillaume Lample, Lucile Saulnier, et~al. 2023.
\newblock Mistral 7b.
\newblock \emph{arXiv preprint arXiv:2310.06825}.

\bibitem[{Jiang et~al.(2016)Jiang, Schroeder, Teran, Stomakhin, and
  Selle}]{jiang2016material}
Chenfanfu Jiang, Craig Schroeder, Joseph Teran, Alexey Stomakhin, and Andrew
  Selle. 2016.
\newblock The material point method for simulating continuum materials.
\newblock In \emph{Acm siggraph 2016 courses}, pages 1--52.

\bibitem[{J{\"u}nger et~al.(1995)J{\"u}nger, Reinelt, and
  Rinaldi}]{junger1995traveling}
Michael J{\"u}nger, Gerhard Reinelt, and Giovanni Rinaldi. 1995.
\newblock The traveling salesman problem.
\newblock \emph{Handbooks in operations research and management science},
  7:225--330.

\bibitem[{Kenton and Toutanova(2019)}]{bert}
Jacob Devlin Ming-Wei~Chang Kenton and Lee~Kristina Toutanova. 2019.
\newblock Bert: Pre-training of deep bidirectional transformers for language
  understanding.
\newblock In \emph{Proceedings of naacL-HLT}, volume~1, page~2. Minneapolis,
  Minnesota.

\bibitem[{Kim et~al.(2024)Kim, Baldi, and McAleer}]{soon2}
Geunwoo Kim, Pierre Baldi, and Stephen McAleer. 2024.
\newblock Language models can solve computer tasks.
\newblock \emph{Advances in Neural Information Processing Systems}, 36.

\bibitem[{Kohonen(1972)}]{store1}
Teuvo Kohonen. 1972.
\newblock Correlation matrix memories.
\newblock \emph{IEEE transactions on computers}, 100(4):353--359.

\bibitem[{Kojima et~al.(2022)Kojima, Gu, Reid, Matsuo, and Iwasawa}]{dra2}
Takeshi Kojima, Shixiang~Shane Gu, Machel Reid, Yutaka Matsuo, and Yusuke
  Iwasawa. 2022.
\newblock Large language models are zero-shot reasoners.
\newblock \emph{Advances in neural information processing systems},
  35:22199--22213.

\bibitem[{Kool et~al.(2018)Kool, Van~Hoof, and Welling}]{kool2018attention}
Wouter Kool, Herke Van~Hoof, and Max Welling. 2018.
\newblock Attention, learn to solve routing problems!
\newblock \emph{arXiv preprint arXiv:1803.08475}.

\bibitem[{Lehman et~al.(2022)Lehman, Gordon, Jain, Ndousse, Yeh, and
  Stanley}]{ELM}
Joel Lehman, Jonathan Gordon, Shawn Jain, Kamal Ndousse, Cathy Yeh, and
  Kenneth~O. Stanley. 2022.
\newblock \href {https://arxiv.org/abs/2206.08896} {Evolution through large
  models}.
\newblock \emph{Preprint}, arXiv:2206.08896.

\bibitem[{Leviathan et~al.(2023)Leviathan, Kalman, and Matias}]{sd}
Yaniv Leviathan, Matan Kalman, and Yossi Matias. 2023.
\newblock Fast inference from transformers via speculative decoding.
\newblock In \emph{International Conference on Machine Learning}, pages
  19274--19286. PMLR.

\bibitem[{Li et~al.(2024{\natexlab{a}})Li, Liu, Fan, Wei, Liu, Tang, and
  Li}]{module_chem}
Jiatong Li, Yunqing Liu, Wenqi Fan, Xiao-Yong Wei, Hui Liu, Jiliang Tang, and
  Qing Li. 2024{\natexlab{a}}.
\newblock Empowering molecule discovery for molecule-caption translation with
  large language models: A chatgpt perspective.
\newblock \emph{IEEE Transactions on Knowledge and Data Engineering}.

\bibitem[{Li et~al.(2024{\natexlab{b}})Li, Zhang, Wang, Hao, Lei, Tan, Zhou,
  Liu, Yang, Xiong, Wang, Chen, Wang, Li, Zhang, Su, Ouyang, Li, and
  Zhou}]{chem_2024}
Junxian Li, Di~Zhang, Xunzhi Wang, Zeying Hao, Jingdi Lei, Qian Tan, Cai Zhou,
  Wei Liu, Yaotian Yang, Xinrui Xiong, Weiyun Wang, Zhe Chen, Wenhai Wang, Wei
  Li, Shufei Zhang, Mao Su, Wanli Ouyang, Yuqiang Li, and Dongzhan Zhou.
  2024{\natexlab{b}}.
\newblock \href {https://arxiv.org/abs/2408.07246} {Chemvlm: Exploring the
  power of multimodal large language models in chemistry area}.
\newblock \emph{Preprint}, arXiv:2408.07246.

\bibitem[{Li et~al.(2024{\natexlab{c}})Li, Zhang, Wang, Hao, Lei, Tan, Zhou,
  Liu, Yang, Xiong, Wang, Chen, Wang, Li, Zhang, Su, Ouyang, Li, and
  Zhou}]{tune3}
Junxian Li, Di~Zhang, Xunzhi Wang, Zeying Hao, Jingdi Lei, Qian Tan, Cai Zhou,
  Wei Liu, Yaotian Yang, Xinrui Xiong, Weiyun Wang, Zhe Chen, Wenhai Wang, Wei
  Li, Shufei Zhang, Mao Su, Wanli Ouyang, Yuqiang Li, and Dongzhan Zhou.
  2024{\natexlab{c}}.
\newblock \href {https://arxiv.org/abs/2408.07246} {Chemvlm: Exploring the
  power of multimodal large language models in chemistry area}.
\newblock \emph{Preprint}, arXiv:2408.07246.

\bibitem[{Liu et~al.(2024{\natexlab{a}})Liu, Lin, Hewitt, Paranjape,
  Bevilacqua, Petroni, and Liang}]{lost}
Nelson~F Liu, Kevin Lin, John Hewitt, Ashwin Paranjape, Michele Bevilacqua,
  Fabio Petroni, and Percy Liang. 2024{\natexlab{a}}.
\newblock Lost in the middle: How language models use long contexts.
\newblock \emph{Transactions of the Association for Computational Linguistics},
  11:157--173.

\bibitem[{Liu et~al.(2024{\natexlab{b}})Liu, Li, Lv, Liu, Zhu, and Hu}]{sd2}
Tianyu Liu, Yun Li, Qitan Lv, Kai Liu, Jianchen Zhu, and Winston Hu.
  2024{\natexlab{b}}.
\newblock \href {https://arxiv.org/abs/2408.11850} {Parallel speculative
  decoding with adaptive draft length}.
\newblock \emph{Preprint}, arXiv:2408.11850.

\bibitem[{Liu et~al.(2024{\natexlab{c}})Liu, Lv, Wang, Yang, and
  Chen}]{restliu}
Tianyu Liu, Qitan Lv, Jie Wang, Shuling Yang, and Hanzhu Chen.
  2024{\natexlab{c}}.
\newblock Learning rule-induced subgraph representations for inductive relation
  prediction.
\newblock \emph{Advances in Neural Information Processing Systems}, 36.

\bibitem[{Liu et~al.(2019)Liu, Ott, Goyal, Du, Joshi, Chen, Levy, Lewis,
  Zettlemoyer, and Stoyanov}]{roberta}
Yinhan Liu, Myle Ott, Naman Goyal, Jingfei Du, Mandar Joshi, Danqi Chen, Omer
  Levy, Mike Lewis, Luke Zettlemoyer, and Veselin Stoyanov. 2019.
\newblock \href {https://arxiv.org/abs/1907.11692} {Roberta: A robustly
  optimized bert pretraining approach}.
\newblock \emph{Preprint}, arXiv:1907.11692.

\bibitem[{Liu et~al.(2021)Liu, Roberts, Lal-Nag, Chen, Huang, and
  Tong}]{ai_drug_2021}
Zhichao Liu, Ruth~A Roberts, Madhu Lal-Nag, Xi~Chen, Ruili Huang, and Weida
  Tong. 2021.
\newblock Ai-based language models powering drug discovery and development.
\newblock \emph{Drug Discovery Today}, 26(11):2593--2607.

\bibitem[{Lu et~al.(2022)Lu, Bartolo, Moore, Riedel, and Stenetorp}]{sen_1}
Yao Lu, Max Bartolo, Alastair Moore, Sebastian Riedel, and Pontus Stenetorp.
  2022.
\newblock Fantastically ordered prompts and where to find them: Overcoming
  few-shot prompt order sensitivity.
\newblock In \emph{Proceedings of the 60th Annual Meeting of the Association
  for Computational Linguistics (Volume 1: Long Papers)}, pages 8086--8098.

\bibitem[{Luo et~al.(2023)Luo, Yang, Meng, Li, Zhou, and Zhang}]{emp}
Yun Luo, Zhen Yang, Fandong Meng, Yafu Li, Jie Zhou, and Yue Zhang. 2023.
\newblock An empirical study of catastrophic forgetting in large language
  models during continual fine-tuning.
\newblock \emph{arXiv preprint arXiv:2308.08747}.

\bibitem[{Lv et~al.(2024)Lv, Wang, Chen, Li, Zhang, and Wu}]{coft_2024}
Qitan Lv, Jie Wang, Hanzhu Chen, Bin Li, Yongdong Zhang, and Feng Wu. 2024.
\newblock Coarse-to-fine highlighting: Reducing knowledge hallucination in
  large language models.
\newblock In \emph{Forty-first International Conference on Machine Learning}.

\bibitem[{Ma et~al.(2024{\natexlab{a}})Ma, Ling, Zhang, and Gu}]{edit1}
Jun-Yu Ma, Zhen-Hua Ling, Ningyu Zhang, and Jia-Chen Gu. 2024{\natexlab{a}}.
\newblock Neighboring perturbations of knowledge editing on large language
  models.
\newblock \emph{arXiv preprint arXiv:2401.17623}.

\bibitem[{Ma et~al.(2024{\natexlab{b}})Ma, Wang, Xu, Ling, and Gu}]{ma}
Jun-Yu Ma, Hong Wang, Hao-Xiang Xu, Zhen-Hua Ling, and Jia-Chen Gu.
  2024{\natexlab{b}}.
\newblock Perturbation-restrained sequential model editing.
\newblock \emph{arXiv preprint arXiv:2405.16821}.

\bibitem[{Ma et~al.(2023{\natexlab{a}})Ma, Chen, Deng, Tenenbaum, Du, Gan, and
  Matusik}]{ma2023learning}
Pingchuan Ma, Peter~Yichen Chen, Bolei Deng, Joshua~B Tenenbaum, Tao Du, Chuang
  Gan, and Wojciech Matusik. 2023{\natexlab{a}}.
\newblock Learning neural constitutive laws from motion observations for
  generalizable pde dynamics.
\newblock In \emph{International Conference on Machine Learning}, pages
  23279--23300. PMLR.

\bibitem[{Ma et~al.(2024{\natexlab{c}})Ma, Wang, Guo, Sun, Tenenbaum, Rus, Gan,
  and Matusik}]{bi_level}
Pingchuan Ma, Tsun-Hsuan Wang, Minghao Guo, Zhiqing Sun, Joshua~B Tenenbaum,
  Daniela Rus, Chuang Gan, and Wojciech Matusik. 2024{\natexlab{c}}.
\newblock Llm and simulation as bilevel optimizers: A new paradigm to advance
  physical scientific discovery.
\newblock \emph{arXiv preprint arXiv:2405.09783}.

\bibitem[{Ma et~al.(2023{\natexlab{b}})Ma, Mishra, Beirami, Beutel, and
  Chen}]{spec1}
Xiao Ma, Swaroop Mishra, Ahmad Beirami, Alex Beutel, and Jilin Chen.
  2023{\natexlab{b}}.
\newblock Let's do a thought experiment: Using counterfactuals to improve moral
  reasoning.
\newblock \emph{arXiv preprint arXiv:2306.14308}.

\bibitem[{Ma et~al.(2023{\natexlab{c}})Ma, Liang, Wang, Huang, Bastani,
  Jayaraman, Zhu, Fan, and Anandkumar}]{eureka}
Yecheng~Jason Ma, William Liang, Guanzhi Wang, De-An Huang, Osbert Bastani,
  Dinesh Jayaraman, Yuke Zhu, Linxi Fan, and Anima Anandkumar.
  2023{\natexlab{c}}.
\newblock Eureka: Human-level reward design via coding large language models.
\newblock In \emph{The Twelfth International Conference on Learning
  Representations}.

\bibitem[{Madaan et~al.(2024)Madaan, Tandon, Gupta, Hallinan, Gao, Wiegreffe,
  Alon, Dziri, Prabhumoye, Yang et~al.}]{reason2}
Aman Madaan, Niket Tandon, Prakhar Gupta, Skyler Hallinan, Luyu Gao, Sarah
  Wiegreffe, Uri Alon, Nouha Dziri, Shrimai Prabhumoye, Yiming Yang, et~al.
  2024.
\newblock Self-refine: Iterative refinement with self-feedback.
\newblock \emph{Advances in Neural Information Processing Systems}, 36.

\bibitem[{Madaan and Yazdanbakhsh(2022)}]{sen_3}
Aman Madaan and Amir Yazdanbakhsh. 2022.
\newblock Text and patterns: For effective chain of thought, it takes two to
  tango.
\newblock \emph{arXiv preprint arXiv:2209.07686}.

\bibitem[{Meng et~al.(2022)Meng, Bau, Andonian, and Belinkov}]{edit2}
Kevin Meng, David Bau, Alex Andonian, and Yonatan Belinkov. 2022.
\newblock Locating and editing factual associations in gpt.
\newblock \emph{Advances in Neural Information Processing Systems},
  35:17359--17372.

\bibitem[{Meng et~al.()Meng, Sharma, Andonian, Belinkov, and Bau}]{edit3}
Kevin Meng, Arnab~Sen Sharma, Alex~J Andonian, Yonatan Belinkov, and David Bau.
\newblock Mass-editing memory in a transformer.
\newblock In \emph{The Eleventh International Conference on Learning
  Representations}.

\bibitem[{Meyerson et~al.(2023{\natexlab{a}})Meyerson, Nelson, Bradley, Gaier,
  Moradi, Hoover, and Lehman}]{tune1}
Elliot Meyerson, Mark~J Nelson, Herbie Bradley, Adam Gaier, Arash Moradi, Amy~K
  Hoover, and Joel Lehman. 2023{\natexlab{a}}.
\newblock Language model crossover: Variation through few-shot prompting.
\newblock \emph{arXiv preprint arXiv:2302.12170}.

\bibitem[{Meyerson et~al.(2023{\natexlab{b}})Meyerson, Nelson, Bradley, Gaier,
  Moradi, Hoover, and Lehman}]{LMX}
Elliot Meyerson, Mark~J Nelson, Herbie Bradley, Adam Gaier, Arash Moradi, Amy~K
  Hoover, and Joel Lehman. 2023{\natexlab{b}}.
\newblock Language model crossover: Variation through few-shot prompting.
\newblock \emph{arXiv preprint arXiv:2302.12170}.

\bibitem[{Mitchell et~al.()Mitchell, Lin, Bosselut, Finn, and Manning}]{edit4}
Eric Mitchell, Charles Lin, Antoine Bosselut, Chelsea Finn, and Christopher~D
  Manning.
\newblock Fast model editing at scale.
\newblock In \emph{International Conference on Learning Representations}.

\bibitem[{Mitchell et~al.(2022)Mitchell, Lin, Bosselut, Manning, and
  Finn}]{edit5}
Eric Mitchell, Charles Lin, Antoine Bosselut, Christopher~D Manning, and
  Chelsea Finn. 2022.
\newblock Memory-based model editing at scale.
\newblock In \emph{International Conference on Machine Learning}, pages
  15817--15831. PMLR.

\bibitem[{Montgomery et~al.(2021)Montgomery, Peck, and
  Vining}]{montgomery2021introduction}
Douglas~C. Montgomery, Elizabeth~A. Peck, and Geoffrey~G. Vining. 2021.
\newblock \emph{Introduction to Linear Regression Analysis}, 6th edition
  edition.
\newblock Wiley.

\bibitem[{Nair et~al.(2023)Nair, Schumacher, Tso, and Kannan}]{diag1}
Varun Nair, Elliot Schumacher, Geoffrey Tso, and Anitha Kannan. 2023.
\newblock Dera: enhancing large language model completions with dialog-enabled
  resolving agents.
\newblock \emph{arXiv preprint arXiv:2303.17071}.

\bibitem[{Nazari et~al.(2018)Nazari, Oroojlooy, Snyder, and
  Tak{\'a}c}]{nazari2018reinforcement}
Mohammadreza Nazari, Afshin Oroojlooy, Lawrence Snyder, and Martin Tak{\'a}c.
  2018.
\newblock Reinforcement learning for solving the vehicle routing problem.
\newblock \emph{Advances in neural information processing systems}, 31.

\bibitem[{Neuberg(2003)}]{casual_graph}
Leland~Gerson Neuberg. 2003.
\newblock Causality: models, reasoning, and inference, by judea pearl,
  cambridge university press, 2000.
\newblock \emph{Econometric Theory}, 19(4):675--685.

\bibitem[{Olausson et~al.(2023)Olausson, Inala, Wang, Gao, and
  Solar-Lezama}]{code4}
Theo~X Olausson, Jeevana~Priya Inala, Chenglong Wang, Jianfeng Gao, and Armando
  Solar-Lezama. 2023.
\newblock Demystifying gpt self-repair for code generation.
\newblock \emph{arXiv preprint arXiv:2306.09896}.

\bibitem[{OpenAI(2020)}]{chatgpt}
OpenAI. 2020.
\newblock Chatgpt: A large-scale generative model for conversation.

\bibitem[{OpenAI(2023)}]{gpt4}
OpenAI. 2023.
\newblock \href {https://arxiv.org/abs/2303.08774} {Gpt-4 technical report}.
\newblock \emph{Preprint}, arXiv:2303.08774.

\bibitem[{Popper(2005)}]{logic_sci}
Karl Popper. 2005.
\newblock \emph{The logic of scientific discovery}.
\newblock Routledge.

\bibitem[{Radford et~al.(2018)Radford, Narasimhan, Salimans, Sutskever
  et~al.}]{gpt1}
Alec Radford, Karthik Narasimhan, Tim Salimans, Ilya Sutskever, et~al. 2018.
\newblock Improving language understanding by generative pre-training.

\bibitem[{Radford et~al.(2019)Radford, Wu, Child, Luan, Amodei, Sutskever
  et~al.}]{gpt6b}
Alec Radford, Jeffrey Wu, Rewon Child, David Luan, Dario Amodei, Ilya
  Sutskever, et~al. 2019.
\newblock Language models are unsupervised multitask learners.
\newblock \emph{OpenAI blog}, 1(8):9.

\bibitem[{Ramakrishnan et~al.(2014)Ramakrishnan, Dral, Rupp, and
  Von~Lilienfeld}]{qm9}
Raghunathan Ramakrishnan, Pavlo~O Dral, Matthias Rupp, and O~Anatole
  Von~Lilienfeld. 2014.
\newblock Quantum chemistry structures and properties of 134 kilo molecules.
\newblock \emph{Scientific data}, 1(1):1--7.

\bibitem[{Rasooli and Tetreault(2015)}]{rasooli-tetrault-2015}
Mohammad~Sadegh Rasooli and Joel~R. Tetreault. 2015.
\newblock \href {http://arxiv.org/abs/1503.06733} {Yara parser: {A} fast and
  accurate dependency parser}.
\newblock \emph{Computing Research Repository}, arXiv:1503.06733.
\newblock Version 2.

\bibitem[{Romera-Paredes et~al.(2024{\natexlab{a}})Romera-Paredes, Barekatain,
  Novikov, Balog, Kumar, Dupont, Ruiz, Ellenberg, Wang, Fawzi
  et~al.}]{nature_math}
Bernardino Romera-Paredes, Mohammadamin Barekatain, Alexander Novikov, Matej
  Balog, M~Pawan Kumar, Emilien Dupont, Francisco~JR Ruiz, Jordan~S Ellenberg,
  Pengming Wang, Omar Fawzi, et~al. 2024{\natexlab{a}}.
\newblock Mathematical discoveries from program search with large language
  models.
\newblock \emph{Nature}, 625(7995):468--475.

\bibitem[{Romera-Paredes et~al.(2024{\natexlab{b}})Romera-Paredes, Barekatain,
  Novikov, Balog, Kumar, Dupont, Ruiz, Ellenberg, Wang, Fawzi
  et~al.}]{funsearch}
Bernardino Romera-Paredes, Mohammadamin Barekatain, Alexander Novikov, Matej
  Balog, M~Pawan Kumar, Emilien Dupont, Francisco~JR Ruiz, Jordan~S Ellenberg,
  Pengming Wang, Omar Fawzi, et~al. 2024{\natexlab{b}}.
\newblock Mathematical discoveries from program search with large language
  models.
\newblock \emph{Nature}, 625(7995):468--475.

\bibitem[{Romera-Paredes et~al.(2024{\natexlab{c}})Romera-Paredes, Barekatain,
  Novikov, Balog, Kumar, Dupont, Ruiz, Ellenberg, Wang, Fawzi et~al.}]{many_na}
Bernardino Romera-Paredes, Mohammadamin Barekatain, Alexander Novikov, Matej
  Balog, M~Pawan Kumar, Emilien Dupont, Francisco~JR Ruiz, Jordan~S Ellenberg,
  Pengming Wang, Omar Fawzi, et~al. 2024{\natexlab{c}}.
\newblock Mathematical discoveries from program search with large language
  models.
\newblock \emph{Nature}, 625(7995):468--475.

\bibitem[{Rosenkrantz et~al.(1977)Rosenkrantz, Stearns, and
  Lewis}]{rosenkrantz1977analysis}
Daniel~J Rosenkrantz, Richard~E Stearns, and Philip~M Lewis, II. 1977.
\newblock An analysis of several heuristics for the traveling salesman problem.
\newblock \emph{SIAM journal on computing}, 6(3):563--581.

\bibitem[{Seber and Lee(2012)}]{seber2012linear}
George A.~F. Seber and Alan~J. Lee. 2012.
\newblock \emph{Linear Regression Analysis}, 2nd edition edition.
\newblock Wiley.

\bibitem[{Sharma and Thakur(2023)}]{gpt_drug}
Gaurav Sharma and Abhishek Thakur. 2023.
\newblock Chatgpt in drug discovery.

\bibitem[{Shi et~al.(2023)Shi, Chen, Misra, Scales, Dohan, Chi, Sch{\"a}rli,
  and Zhou}]{distract}
Freda Shi, Xinyun Chen, Kanishka Misra, Nathan Scales, David Dohan, Ed~H Chi,
  Nathanael Sch{\"a}rli, and Denny Zhou. 2023.
\newblock Large language models can be easily distracted by irrelevant context.
\newblock In \emph{International Conference on Machine Learning}, pages
  31210--31227. PMLR.

\bibitem[{Shinn et~al.(2023)Shinn, Labash, and Gopinath}]{reason1}
Noah Shinn, Beck Labash, and Ashwin Gopinath. 2023.
\newblock Reflexion: an autonomous agent with dynamic memory and
  self-reflection.
\newblock \emph{arXiv preprint arXiv:2303.11366}, 2(5):9.

\bibitem[{Shoeybi et~al.(2019)Shoeybi, Patwary, Puri, LeGresley, Casper, and
  Catanzaro}]{megantron}
Mohammad Shoeybi, Md.MostofaAli Patwary, Raul Puri, Patrick LeGresley, Jared
  Casper, and Bryan Catanzaro. 2019.
\newblock Megatron-lm: Training multi-billion parameter language models using
  model parallelism.
\newblock \emph{Cornell University - arXiv,Cornell University - arXiv}.

\bibitem[{Sulsky et~al.(1995)Sulsky, Zhou, and
  Schreyer}]{sulsky1995application}
Deborah Sulsky, Shi-Jian Zhou, and Howard~L Schreyer. 1995.
\newblock Application of a particle-in-cell method to solid mechanics.
\newblock \emph{Computer physics communications}, 87(1-2):236--252.

\bibitem[{Sumers et~al.()Sumers, Yao, Narasimhan, and Griffiths}]{LATM}
Theodore Sumers, Shunyu Yao, Karthik Narasimhan, and Thomas Griffiths.
\newblock Cognitive architectures for language agents.
\newblock \emph{Transactions on Machine Learning Research}.

\bibitem[{Team(2024{\natexlab{a}})}]{claude3}
Anthropic Team. 2024{\natexlab{a}}.
\newblock \href {https://api.semanticscholar.org/CorpusID:268232499} {The
  claude 3 model family: Opus, sonnet, haiku}.

\bibitem[{Team(2023{\natexlab{a}})}]{gemini}
Gemini Team. 2023{\natexlab{a}}.
\newblock \href {https://arxiv.org/abs/2312.11805} {Gemini: A family of highly
  capable multimodal models}.
\newblock \emph{Preprint}, arXiv:2312.11805.

\bibitem[{Team(2023{\natexlab{b}})}]{internlm}
InternLM Team. 2023{\natexlab{b}}.
\newblock Internlm: A multilingual language model with progressively enhanced
  capabilities.

\bibitem[{Team(2023{\natexlab{c}})}]{llama2}
Llama2 Team. 2023{\natexlab{c}}.
\newblock \href {https://arxiv.org/abs/2307.09288} {Llama 2: Open foundation
  and fine-tuned chat models}.
\newblock \emph{Preprint}, arXiv:2307.09288.

\bibitem[{Team(2024{\natexlab{b}})}]{llama3}
Llama3 Team. 2024{\natexlab{b}}.
\newblock The llama 3 herd of models.
\newblock \emph{arXiv preprint arXiv:2407.21783}.

\bibitem[{Team(2022)}]{palm}
PaLM Team. 2022.
\newblock \href {https://arxiv.org/abs/2204.02311} {Palm: Scaling language
  modeling with pathways}.
\newblock \emph{Preprint}, arXiv:2204.02311.

\bibitem[{Trinh et~al.(2024)Trinh, Wu, Le, He, and Luong}]{geo}
Trieu~H Trinh, Yuhuai Wu, Quoc~V Le, He~He, and Thang Luong. 2024.
\newblock Solving olympiad geometry without human demonstrations.
\newblock \emph{Nature}, 625(7995):476--482.

\bibitem[{Udrescu et~al.(2020)Udrescu, Tan, Feng, Neto, Wu, and
  Tegmark}]{limi1}
Silviu-Marian Udrescu, Andrew Tan, Jiahai Feng, Orisvaldo Neto, Tailin Wu, and
  Max Tegmark. 2020.
\newblock Ai feynman 2.0: Pareto-optimal symbolic regression exploiting graph
  modularity.
\newblock \emph{Advances in Neural Information Processing Systems},
  33:4860--4871.

\bibitem[{Wang et~al.(2023{\natexlab{a}})Wang, Xie, Jiang, Mandlekar, Xiao,
  Zhu, Fan, and Anandkumar}]{soon1}
Guanzhi Wang, Yuqi Xie, Yunfan Jiang, Ajay Mandlekar, Chaowei Xiao, Yuke Zhu,
  Linxi Fan, and Anima Anandkumar. 2023{\natexlab{a}}.
\newblock Voyager: An open-ended embodied agent with large language models.
\newblock \emph{arXiv preprint arXiv:2305.16291}.

\bibitem[{Wang et~al.(2023{\natexlab{b}})Wang, Fu, Du, Gao, Huang, Liu,
  Chandak, Liu, Van~Katwyk, Deac et~al.}]{ai_age}
Hanchen Wang, Tianfan Fu, Yuanqi Du, Wenhao Gao, Kexin Huang, Ziming Liu, Payal
  Chandak, Shengchao Liu, Peter Van~Katwyk, Andreea Deac, et~al.
  2023{\natexlab{b}}.
\newblock Scientific discovery in the age of artificial intelligence.
\newblock \emph{Nature}, 620(7972):47--60.

\bibitem[{Wang et~al.(2023{\natexlab{c}})Wang, Ma, Dong, Huang, Wang, Ma, Yang,
  Wang, Wu, and Wei}]{eff}
Hongyu Wang, Shuming Ma, Li~Dong, Shaohan Huang, Huaijie Wang, Lingxiao Ma, Fan
  Yang, Ruiping Wang, Yi~Wu, and Furu Wei. 2023{\natexlab{c}}.
\newblock Bitnet: Scaling 1-bit transformers for large language models.
\newblock \emph{arXiv preprint arXiv:2310.11453}.

\bibitem[{Wang et~al.()Wang, Variengien, Conmy, Shlegeris, and
  Steinhardt}]{104interpretability}
Kevin~Ro Wang, Alexandre Variengien, Arthur Conmy, Buck Shlegeris, and Jacob
  Steinhardt.
\newblock Interpretability in the wild: a circuit for indirect object
  identification in gpt-2 small.
\newblock In \emph{The Eleventh International Conference on Learning
  Representations}.

\bibitem[{Wei et~al.(2022)Wei, Wang, Schuurmans, Bosma, Xia, Chi, Le, Zhou
  et~al.}]{wei2022chain}
Jason Wei, Xuezhi Wang, Dale Schuurmans, Maarten Bosma, Fei Xia, Ed~Chi, Quoc~V
  Le, Denny Zhou, et~al. 2022.
\newblock Chain-of-thought prompting elicits reasoning in large language
  models.
\newblock \emph{Advances in neural information processing systems},
  35:24824--24837.

\bibitem[{Wei et~al.(2023)Wei, Wei, Tay, Tran, Webson, Lu, Chen, Liu, Huang,
  Zhou et~al.}]{sen_2}
Jerry Wei, Jason Wei, Yi~Tay, Dustin Tran, Albert Webson, Yifeng Lu, Xinyun
  Chen, Hanxiao Liu, Da~Huang, Denny Zhou, et~al. 2023.
\newblock Larger language models do in-context learning differently.
\newblock \emph{arXiv preprint arXiv:2303.03846}.

\bibitem[{Wuestman et~al.(2020)Wuestman, Hoekman, and Frenken}]{inspired}
Mignon Wuestman, Jarno Hoekman, and Koen Frenken. 2020.
\newblock A typology of scientific breakthroughs.
\newblock \emph{Quantitative Science Studies}, 1(3):1203--1222.

\bibitem[{Yang et~al.(2023)Yang, Wang, Lu, Liu, Le, Zhou, and Chen}]{opro}
Chengrun Yang, Xuezhi Wang, Yifeng Lu, Hanxiao Liu, Quoc~V. Le, Denny Zhou, and
  Xinyun Chen. 2023.
\newblock Large language models as optimizers.

\bibitem[{Young et~al.(2024)Young, Chen, Li, Huang, Zhang, Zhang, Li, Zhu,
  Chen, Chang et~al.}]{yi9b}
Alex Young, Bei Chen, Chao Li, Chengen Huang, Ge~Zhang, Guanwei Zhang, Heng Li,
  Jiangcheng Zhu, Jianqun Chen, Jing Chang, et~al. 2024.
\newblock Yi: Open foundation models by 01. ai.
\newblock \emph{arXiv preprint arXiv:2403.04652}.

\bibitem[{Yu et~al.(2024)Yu, Chen, Zhou, and He}]{pre_edit2}
Lang Yu, Qin Chen, Jie Zhou, and Liang He. 2024.
\newblock Melo: Enhancing model editing with neuron-indexed dynamic lora.
\newblock In \emph{Proceedings of the AAAI Conference on Artificial
  Intelligence}, volume~38, pages 19449--19457.

\bibitem[{Yu et~al.(2023)Yu, Yang, Liu, and Ananiadou}]{inter_analy}
Zeping Yu, Kailai Yang, Zhiwei Liu, and Sophia Ananiadou. 2023.
\newblock Exploring the residual stream of transformers.
\newblock \emph{arXiv preprint arXiv:2312.12141}.

\bibitem[{Yuan et~al.(2024)Yuan, Cho, and Weston}]{diag2}
Weizhe Yuan, Kyunghyun Cho, and Jason Weston. 2024.
\newblock System-level natural language feedback.
\newblock In \emph{Proceedings of the 18th Conference of the European Chapter
  of the Association for Computational Linguistics (Volume 1: Long Papers)},
  pages 2773--2789.

\bibitem[{Zhang et~al.(2024{\natexlab{a}})Zhang, Liu, Tan, Chen, Yan, Yan, Li,
  Huang, Yue, Zhou et~al.}]{tune2}
Di~Zhang, Wei Liu, Qian Tan, Jingdan Chen, Hang Yan, Yuliang Yan, Jiatong Li,
  Weiran Huang, Xiangyu Yue, Dongzhan Zhou, et~al. 2024{\natexlab{a}}.
\newblock Chemllm: A chemical large language model.
\newblock \emph{arXiv preprint arXiv:2402.06852}.

\bibitem[{Zhang et~al.(2024{\natexlab{b}})Zhang, Qi, and Zhou}]{robots2}
Kaiyan Zhang, Biqing Qi, and Bowen Zhou. 2024{\natexlab{b}}.
\newblock Towards building specialized generalist ai with system 1 and system 2
  fusion.
\newblock \emph{arXiv preprint arXiv:2407.08642}.

\bibitem[{Zhang et~al.(2024{\natexlab{c}})Zhang, Yao, Tian, Wang, Deng, Wang,
  Xi, Mao, Zhang, Ni et~al.}]{ME_survey}
Ningyu Zhang, Yunzhi Yao, Bozhong Tian, Peng Wang, Shumin Deng, Mengru Wang,
  Zekun Xi, Shengyu Mao, Jintian Zhang, Yuansheng Ni, et~al.
  2024{\natexlab{c}}.
\newblock A comprehensive study of knowledge editing for large language models.
\newblock \emph{arXiv preprint arXiv:2401.01286}.

\bibitem[{Zhang et~al.()Zhang, Pearson, and Wang}]{robots1}
Starkson Zhang, Alfredo Pearson, and Zhenting Wang.
\newblock Autonomous generalist scientist: Towards and beyond human-level
  automatic research using foundation model-based ai agents and robots (a
  position).

\bibitem[{Zhao et~al.(2024)Zhao, Lin, Zhu, Ye, Chen, Zheng, Ceze,
  Krishnamurthy, Chen, and Kasikci}]{eff3}
Yilong Zhao, Chien-Yu Lin, Kan Zhu, Zihao Ye, Lequn Chen, Size Zheng, Luis
  Ceze, Arvind Krishnamurthy, Tianqi Chen, and Baris Kasikci. 2024.
\newblock Atom: Low-bit quantization for efficient and accurate llm serving.
\newblock \emph{Proceedings of Machine Learning and Systems}, 6:196--209.

\bibitem[{Zheng et~al.(2023)Zheng, Yin, Xie, Sun, Huang, Yu, Cao, Kozyrakis,
  Stoica, Gonzalez et~al.}]{limi2}
Lianmin Zheng, Liangsheng Yin, Zhiqiang Xie, Chuyue Sun, Jeff Huang, Cody~Hao
  Yu, Shiyi Cao, Christos Kozyrakis, Ion Stoica, Joseph~E Gonzalez, et~al.
  2023.
\newblock Sglang: Efficient execution of structured language model programs.

\bibitem[{Zhong et~al.(2023)Zhong, Wu, Manning, Potts, and Chen}]{assess}
Zexuan Zhong, Zhengxuan Wu, Christopher~D Manning, Christopher Potts, and Danqi
  Chen. 2023.
\newblock Mquake: Assessing knowledge editing in language models via multi-hop
  questions.
\newblock \emph{arXiv preprint arXiv:2305.14795}.

\bibitem[{Zhou et~al.()Zhou, Muresanu, Han, Paster, Pitis, Chan, and Ba}]{dra3}
Yongchao Zhou, Andrei~Ioan Muresanu, Ziwen Han, Keiran Paster, Silviu Pitis,
  Harris Chan, and Jimmy Ba.
\newblock Large language models are human-level prompt engineers.
\newblock In \emph{The Eleventh International Conference on Learning
  Representations}.

\bibitem[{Zhu et~al.(2020)Zhu, Rawat, Zaheer, Bhojanapalli, Li, Yu, and
  Kumar}]{modify}
Chen Zhu, Ankit~Singh Rawat, Manzil Zaheer, Srinadh Bhojanapalli, Daliang Li,
  Felix Yu, and Sanjiv Kumar. 2020.
\newblock Modifying memories in transformer models.
\newblock \emph{arXiv preprint arXiv:2012.00363}.

\bibitem[{Zhu et~al.(2023)Zhu, Li, Liu, Ma, and Wang}]{qinglianghua}
Xunyu Zhu, Jian Li, Yong Liu, Can Ma, and Weiping Wang. 2023.
\newblock A survey on model compression for large language models.
\newblock \emph{arXiv preprint arXiv:2308.07633}.

\end{thebibliography}
\end{document}